\documentclass[11pt,a4paper]{article}

\usepackage[english]{babel}
\usepackage[utf8]{inputenc}

\usepackage{amsmath,amssymb,amsbsy,amsfonts}
\usepackage{bm}
\usepackage{stmaryrd}

\usepackage{ragged2e}
\justifying

\usepackage[paperwidth=192mm,
		    paperheight=262mm,
		    vmargin={19mm,19mm},
		    hmargin={13.7mm,13.7mm},
		    headsep=12pt,
		    footskip=12pt]{geometry}

\usepackage{hyperref}

\usepackage{siunitx}
\usepackage{graphicx}
\graphicspath{{./figures/}}
\usepackage{float}
\usepackage{caption}
\usepackage{subcaption}

\usepackage{tabularx}
\usepackage{multirow}
\usepackage{booktabs}

\usepackage[backend=bibtex,style=numeric,sorting=nyt,defernumbers=true]{biblatex}
\renewbibmacro{in:}{}
\DeclareBibliographyCategory{cited}
\AtEveryCitekey{\addtocategory{cited}{\thefield{entrykey}}}
\DeclareMathOperator{\grad}{\nabla}
\DeclareMathOperator{\dive}{\nabla\cdot}

\newcommand{\vel}{\bm{u}}
\newcommand{\pmix}{\overline{p}}

\newcommand{\ave}[1]{\left\{\!\!\left\{#1\right\}\!\!\right\}}
\newcommand{\jump}[1]{\left\llbracket#1\right\rrbracket}
\newcommand{\odif}[2]{\frac{\mathrm{d}{#1}}{\mathrm{d}{#2}}}
\newcommand{\pad}[2]{\frac{\partial{#1}}{\partial{#2}}}

\newcommand{\rpth}[1]{\left(#1\right)}
\newcommand{\spth}[1]{\left[#1\right]}
\newcommand{\cpth}[1]{\left\{#1\right\}}

\newcommand{\Crouzet}{\textit{NC-centered-2013}}

\begin{document}

\title{A robust computational framework for the mixture-energy-consistent six-equation two-phase model with instantaneous mechanical relaxation terms}

\author{Giuseppe Orlando$^{(1)}$, Ward Haegeman$^{(2,1)}$ \\ Marica Pelanti$^{(3)}$, Marc Massot$^{(1)}$}

\date{}

\maketitle

\begin{center}
{
    \small
    $^{(1)}$
    CMAP, CNRS, \'{E}cole polytechnique, Institut Polytechnique de Paris \\
    Route de Saclay, 91120 Palaiseau, France \\
    {\tt giuseppe.orlando@polytechnique.edu, ward.haegeman@polytechnique.edu, marc.massot@polytechnique.edu} \\
    \ \\
    $^{(2)}$
    Department of Multi-Physics for Energetics, ONERA, Université Paris-Saclay \\ 
    F-91123 Palaiseau, France \\
    \ \\
    $^{(3)}$
    IMSIA, UMR 9219 ENSTA-CNRS-EDF, \\
    ENSTA Paris, Institut Polytechnique de Paris,\\
    91120 Palaiseau, France \\
    {\tt marica.pelanti@ensta.fr} \\
}
\end{center}

\noindent

{\bf Keywords}: Two-phase flows, 6-equation model, Finite Volume schemes, Wave-propagation scheme, Non-conservative terms discretization, Riemann solvers, Relaxation source terms.

\pagebreak

\begin{abstract}
	We present a robust computational framework for the numerical solution of a hyperbolic 6-equation single-velocity two-phase system. The system's main interest is that, when combined with instantaneous mechanical relaxation, it recovers the solution of the 5-equation model of Kapila. Several numerical methods based on this strategy have been developed over the years. However, neither the 5- nor 6-equation model admits a complete set of jump conditions because they involve non-conservative products. Different discretizations of these terms in the 6-equation model exist. The precise impact of these discretizations on the numerical solutions of the 5-equation model, in particular for shocks, is still an open question to which this work provides new insights. We consider the phasic total energies as prognostic variables to naturally enforce discrete conservation of total energy and compare the accuracy and robustness of different discretizations for the hyperbolic operator. Namely, we discuss the construction of an HLLC approximate Riemann solver in relation to jump conditions. We then compare an HLLC wave-propagation scheme which includes the non-conservative terms, with Rusanov and HLLC solvers for the conservative part in combination with suitable approaches for the non-conservative terms. We show that some approaches for the discretization of non-conservative terms fit within the framework of path-conservative schemes for hyperbolic problems. We then analyze the use of various numerical strategies on several relevant test cases, showing both the impact of the theoretical shortcomings of the models as well as the importance of the choice of a robust framework for the global numerical strategy.
\end{abstract}
\setlength{\leftskip}{0pt}
\setlength{\rightskip}{0pt}

\pagebreak

\section{Introduction}
\label{sec:intro} \indent

Two-phase flows are relevant in several areas of engineering, such as naval engineering or aerospace engineering, and are characterized by the simultaneous presence of two phases with different properties. A widely established model for compressible two-phase flows is the full non-equilibrium 7-equation Baer-Nunziato (BN) model \cite{baer:1986}, originally proposed for detonation waves in granular explosives, and then extended to liquid-gas flows in \cite{saurel:1999}. The model assumes that each phase is compressible and evolves with its own pressure, temperature, and velocity, together with an evolution equation for the volume fraction of one of the two phases. Taking into account the compressibility of both phases is fundamental to correctly capture wave-propagation phenomena and acoustic perturbations, and it is crucial in the case of liquid-vapour transition \cite{saurel:2008}.

Reduced models have been derived by means of asymptotic expansions of the BN model with assumptions of infinitely fast relaxation towards the equilibrium for the velocity, pressure or temperature, see, e.g., \cite{kapila:2001, pelanti:2022}. Among them, a widely employed model is the 5-equation single-velocity single-pressure model, also known as Kapila model \cite{kapila:2001}. The main difficulty in the design of robust numerical schemes for the 5-equation model stems from the non-conservative term in the volume fraction equation that depends on the divergence of the flow velocity and on the phasic compressibilities \cite{pelanti:2014, saurel:2009}. The variation of the volume fraction across shock waves associated with this term makes the construction of approximate Riemann solvers more challenging with respect to the 6-equation model. Moreover, the presence of this non-conservative term in the equation of the volume fraction makes it difficult to preserve the bounds of the volume fraction when dealing with shocks and strong rarefaction waves \cite{saurel:2009}.

Alternative strategies for the numerical solution of the 5-equation model have been proposed in the literature relying on the relaxation of non-equilibrium models. Following the formal derivation performed at the continuous level, one could start from the numerical solution of the BN model and then perform instantaneous velocity and pressure relaxations. The BN model is also characterized by non-conservative terms, but, under specific choices of the interfacial pressure and velocity, it is well posed, in the sense that non-conservative products are uniquely defined and the model admits a mixture entropy inequality to select physically relevant weak solutions \cite{gallouet:2004}. Its numerical solution is rather challenging, even though numerical methods which are positivity-preserving and entropy-consistent are nowadays available \cite{coquel:2017, saleh:2019}. However, the wave structure of the BN model is definitely more complex than that of the reduced Kapila model, as it exhibits a larger number of acoustic waves and involves multiple velocities \cite{tokareva:2010}.

An alternative and popular strategy consists in resolving the 6-equation single-velocity model with instantaneous mechanical relaxation \cite{pelanti:2014, petitpas:2007, saurel:2009}. This model assumes velocity equilibrium between the two phases, keeping mechanical, thermal, and chemical non-equilibrium effects. More specifically, we consider here the 6-equation single-velocity model as presented in \cite{pelanti:2014}. Hence, the model consists of an advection equation for the volume fraction of one phase, together with the continuity equation and the total energy balance of each phase, and a mixture momentum equation. Alternative formulations of the 6-equation model have been proposed in the literature using the equations of the phasic internal energies instead of those of the total energies, see \cite{saurel:2009, zein:2010}. However, as discussed in \cite{pelanti:2014}, while the two systems are mathematically equivalent for smooth solutions, the use of phasic total energies provides several advantages in terms of numerical robustness. In particular, it ensures automatically conservation of the mixture total energy. Formulations based on phasic internal energies require instead the use of an additional conservation law to correct the thermodynamic state and to guarantee the conservation of the total energy \cite{saurel:2009, schmidmayer:2023, zein:2010}. The main advantage of relaxing the 6-equation model instead of the 7-equation model is that its wave structure is similar to that of the 5-equation model. However, the 6-equation model contains non-conservative terms for which, unlike the BN model, a complete set of jump conditions cannot be uniquely defined, making the design of Riemann solvers and the choice of suitable numerical methods challenging. The choice of starting from the 7- or 6-equation model reflects the trade-off between the good mathematical properties of the 7-equation model and the simpler wave structure of the 6-equation model and its reduced number of equations, in particular for multi-dimensional computations. In this work, we focus on the latter of these two strategies.

The aim of the present paper is to provide a quantitative comparison among different numerical discretization strategies for the aforementioned 6-equation model so as to identify a robust computational framework that can be employed in combination with the instantaneous mechanical relaxation. In most standard contributions \cite{pelanti:2014, saurel:2009, saurel:2007}, the 6-equation model has been employed only as an auxiliary tool to solve the 5-equation model, thus combining its numerical solution with instantaneous mechanical relaxation. As a consequence, most existing discretizations have been specifically designed and assessed only in the pressure-equilibrium regime, and little attention has been paid to their behaviour in the pressure non-equilibrium 6-equation setting. However, the accurate and robust numerical solution of the 6-equation model is a crucial step, and it is therefore important to assess the properties of the numerical method employed to solve the 6-equation model. Furthermore, some recent contributions equipped the 6-equation model with finite-rate mechanical relaxation \cite{pelanti:2022,schmidmayer:2023}. When doing so, the model is no longer an auxiliary system used for the numerical resolution of the Kapila model, but acquires the status of a standalone model; thus investigating its discretization is even more important, and this is one of the goals of the present work.

The spatial discretization is based on the finite volume (FV) method as implemented in the framework of \texttt{samurai}\footnote{\url{https://github.com/hpc-maths/samurai}} and approximate Riemann solvers, such as the Rusanov flux \cite{rusanov:1962} and the HLLC \cite{leveque:2002}, are employed. Several numerical methods for the 6-equation model have been developed over the years, see, among many others \cite{delorenzo:2018, pelanti:2014, petitpas:2007, saurel:2009} and the references therein. In this work, we compare in terms of robustness and accuracy two numerical strategies; in the first one, we employ the approximate Riemann solver for the conservative part in combination with suitable approaches for the non-conservative terms \cite{bassi:1997a, crouzet:2013, orlando:2023, tumolo:2015}, whereas in the second one we adopt the so-called wave-propagation method as described in \cite{leveque:2002}. This second approach has been employed for the solution of the 6-equation model with instantaneous mechanical relaxation in \cite{pelanti:2022, pelanti:2014}.

Instantaneous mechanical relaxation may be applied after each hyperbolic step. We have already mentioned the advantages in resolving the 6-equation model with instantaneous relaxation compared to a direct discretization of the 5-equation model. However, we show here that, when combined with instantaneous relaxation, all numerical schemes yield virtually identical results for the 5-equation model, except under extreme conditions, where different numerical schemes converge to different solutions. Since jump conditions are not uniquely defined, different discretization techniques of the hyperbolic operator, in particular the non-conservative products, may lead to different partitions of the phasic energies and ultimately lead to different shock profiles, even though global conservation is ensured. The analysis of the discretization strategies for the 6-equation model carried out in this contribution is complementary to a recent work of two of the authors in which a complete analysis of the pressure relaxation operator is provided \cite{haegeman:2024}. In their work, it has been shown how, for a given spatial discretization of the 6-equation model, the choice of the pressure relaxation operator may impact the numerical results. In particular, the interfacial pressure plays a role only in the relaxation process and does not appear in the 5-equation model. The choice of the interfacial pressure in the relaxation process, as well as the relaxation techniques, may impact the robustness of the numerical strategy or affect the computed shock profiles. The impact of the spatial discretization scheme for the hyperbolic operator of the 6-equation model was not part of the study presented in \cite{haegeman:2024}. This is what the present work focuses on as we consider various treatments for the spatial discretization of the non-conservative products present in the model and consider only a single relaxation strategy. Combined, the present contribution and the one in \cite{haegeman:2024}, aim at providing a complete and in-depth analysis of the overall numerical strategy. In conclusion, the goal of the paper is to show a robust computational framework for the solution of the 6-equation model that can then be employed in combination with the instantaneous mechanical relaxation to solve numerically the 5-equation model, limiting the impact of the shortcomings of the models. In order to achieve the goal, we present a novel analysis of the jump conditions and construction of approximate Riemann solvers as well as a novel perspective on the numerical treatment of non-conservative terms.

The paper is structured as follows. The 6-equation single-velocity two-phase model with instantaneous mechanical relaxation is presented in Section \ref{sec:model}. The different numerical strategies for the discretization of the hyperbolic operator are reported in Section \ref{sec:num_hyp}. Here we also briefly recall the wave-propagation method. Since the 6-equation model does not admit a complete set of jump conditions, a detailed discussion on the hypotheses employed to build the HLLC approximate Riemann solver is provided. Numerical results on a number of benchmarks to assess the robustness and the accuracy of the different approaches presented in Section \ref{sec:num_hyp} for the 6-equation model are reported in Section \ref{sec:num_res_homogeneous}. We stress the fact that, to the best of our knowledge, a detailed numerical assessment of the 6-equation model as independent system has never been performed in the literature. This is because the 6-equation model has traditionally been employed merely as an auxiliary tool for solving the Kapila model. However, when finite-rate relaxation is adopted, as in \cite{pelanti:2022, schmidmayer:2023}, this is no longer the case, and the 6-equation model assumes the status of a standalone system. Following the outcomes of these results, we show that some of the approaches employed for the  discretization of non-conservative terms, originally designed for diffusion operators, can be reformulated in the framework of the well-known path conservative schemes for hyperbolic problems \cite{pares:2006}, so as to justify why certain approaches are more robust than other ones. The analysis of the instantaneous mechanical relaxation is presented in Section \ref{sec:relax}. The quantitative numerical assessment of the mechanical relaxation is presented in Section \ref{sec:num_res_relax}. Finally, some conclusions and perspectives for future work are presented in Section \ref{sec:conclu}.

\section{Two-phase models}
\label{sec:model}

In this section, the 6-equation model followed by the 5-equation model and their properties are described. Endowed with instantaneous mechanical relaxation source terms, the 6-equation model may be used as an auxiliary model for the resolution of the 5-equation model of Kapila. Using the formulation proposed in \cite{pelanti:2014}, the model reads as follows
\begin{subequations}
	\begin{align}
		\pad{\alpha_{1}}{t} + \vel \cdot \grad\alpha_{1} &= \mu\rpth{p_{1} - p_{2}}, \label{eq:volume_fraction_advection} \\
		\pad{\alpha_{1}\rho_{1}}{t} + \dive\rpth{\alpha_{1}\rho_{1}\vel} &= 0, \\
		\pad{\alpha_{2}\rho_{2}}{t} + \dive\rpth{\alpha_{2}\rho_{2}\vel} &= 0, \\
		\pad{\rho\vel}{t} + \dive\rpth{\rho\vel \otimes \vel} + \grad\pmix &= \bm{0}, \label{eq:mixture_momentum} \\
		\pad{\alpha_{1}\rho_{1}E_{1}}{t} + \dive\rpth{\alpha_{1}\rho_{1}E_{1}\vel + \alpha_{1}p_{1}\vel} + \vel \cdot \boldsymbol{\Sigma} &= -\mu p_{I}\rpth{p_{1} - p_{2}}, \\
		\pad{\alpha_{2}\rho_{2}E_{2}}{t} + \dive\rpth{\alpha_{2}\rho_{2}E_{2}\vel + \alpha_{2}p_{2}\vel} - \vel \cdot \boldsymbol{\Sigma} &= \mu p_{I}\rpth{p_{1} - p_{2}},
	\end{align}
\end{subequations}
where the non-conservative term $\boldsymbol{\Sigma}$ is
\begin{equation}\label{eq:non_conservative_term}
	\boldsymbol{\Sigma} = -\spth{Y_{2}\grad\rpth{\alpha_{1}p_{1}} - Y_{1}\grad\rpth{\alpha_{2}p_{2}}}.
\end{equation}
Here $\alpha_{k}$ is the volume fraction of phase $k, k=1,2$, $\rho_{k}$ is the density of phase $k$, $p_{k}$ is the pressure of phase $k$, $E_{k}$ is the total energy per unit of mass of phase $k$, and $\vel$ is the velocity field. Moreover, the mixture density is defined as $\rho = \alpha_{1}\rho_{1} + \alpha_{2}\rho_{2},$ and the mixture pressure $\pmix$ is given by $\pmix = \alpha_{1}p_{1} + \alpha_{2}p_{2}.$ The notation $Y_{k}$ denotes the mass fraction of phase $k$, defined as $Y_{k} = \frac{\alpha_{k}\rho_{k}}{\rho}.$ We recall that the volume fractions satisfy the saturation constraint $\alpha_{1} + \alpha_{2} = 1.$

The previous set of equations has to be closed by specifying an equation of state (EOS) for each phase. A complete EOS for each phase corresponds to specifying the entropies per unit of mass $s_{k}$ as a function of the specific volumes $\rho_{k}^{-1}$ and the internal energies per unit of mass $e_{k} = E_{k} - \left|\vel\right|^{2}/2$, such that $\rpth{\rho_{k}^{-1},e_{k}} \mapsto s_{k}$ are concave. The pressure laws, $p_{k} = p_{k}\rpth{e_{k},\rho_{k}}$, are then defined through the Gibbs identity
$$T_{k}\mathrm{d}s_{k} = \mathrm{d}e_{k} + p_{k}\mathrm{d}\rho_{k}^{-1},$$
with $T_{k} > 0$, the temperature of each phase. In this work, we restrict our attention to the stiffened gas equation of state (SG-EOS) \cite{metayer:2016}, for which the pressure law writes
\begin{equation}
	p_{k}\rpth{e_{k}, \rho_{k}} = \rpth{\gamma_{k} - 1}\rho_{k}e_{k} - \gamma_{k}\pi_{k} - \rpth{\gamma_{k} - 1}\eta_{k}\rho_{k}, \label{eq:SG_EOS}
\end{equation}
where $\gamma_{k}, \pi_{k}$, and $\eta_{k}$ are constant material-dependent parameters. The speed of sound $c_{k}$ of phase $k$ reads therefore \cite{pelanti:2014}
\begin{equation}
	c_{k} = \sqrt{\gamma_{k}\frac{p_{k} + \pi_{k}}{\rho_{k}}}.
\end{equation}
When values of the flow state variables are defined in the region of thermodynamic stability, the model is hyperbolic, i.e. it has real eigenvalues and a complete set of eigenvectors \cite{pelanti:2014, zein:2010}. The linearly degenerate --- or convective --- eigenvalues write $\lambda = u$ with multiplicity $3+d$, $d$ being the number of space dimensions, while the genuinely nonlinear --- or acoustic --- eigenvalues write $\lambda = u \pm c_{\text{f}}$ and have multiplicity one. Here, the mixture speed of sound $c_{\text{f}}$ of the 6-equation model, also known as frozen speed of sound, is defined through
\begin{equation}\label{eq:speed_sound_mixture}
	c^{2}_{\text{f}} = Y_{1}c_{1}^{2} + Y_{2}c_{2}^{2}.
\end{equation}

The main advantage of expressing the model in terms of phasic total energies is that the sum of the two phasic total energy equations recovers immediately the conservation equation of the mixture total energy
\begin{equation}\label{eq:conservation_total_energy}
	\pad{\rho E}{t} + \dive\rpth{\rho E\vel + \pmix\vel} = 0,
\end{equation}
where $\rho E = \alpha_{1}\rho_{1}E_{1} + \alpha_{2}\rho_{2}E_{2}.$ This avoids the need of any type of energy correction to ensure numerically total energy conservation \cite{saurel:2009, schmidmayer:2023} and this property is guaranteed for the corresponding discrete equations, provided that symmetric discretizations of the non-conservative terms are used.

Finally, $\mu > 0$ represents the pressure relaxation rate and $p_{I}$ denotes the interfacial pressure.  When an infinite-rate pressure relaxation is considered, i.e. $\mu \to \infty$, mechanical equilibrium is reached instantaneously and one recovers the solutions to the 5-equation model of Kapila \cite{kapila:2001, pelanti:2014, saurel:2009} which writes
\begin{subequations}\label{eq:Kapila_system}
	\begin{align}
		\pad{\alpha_{1}}{t} + \vel \cdot \grad\alpha_{1} + \frac{\rho_1 c_1^2 - \rho_2 c_2^2}{\frac{\rho_1 c_1^2}{\alpha_1}+\frac{\rho_2 c_2^2}{\alpha_2}}\dive\vel &= 0, \\
		\pad{\alpha_{1}\rho_{1}}{t} + \dive\rpth{\alpha_{1}\rho_{1}\vel} &= 0, \\
		\pad{\alpha_{2}\rho_{2}}{t} + \dive\rpth{\alpha_{2}\rho_{2}\vel} &= 0, \\
		\pad{\rho\vel}{t} + \dive\rpth{\rho\vel \otimes \vel} + \grad p &= \bm{0}, \\
		\pad{\rho E}{t} + \dive\rpth{\rho E\vel + p\vel}  &= 0.
	\end{align}
\end{subequations}
The notations are identical to those defined previously. Let us underline that $p$ denotes the common pressure as the phases remain at pressure equilibrium at all times. This model is also hyperbolic, with a similar eigenstructure as the previous one, i.e. simple acoustic eigenvalues and convective eigenvalues with multiplicity $2+d$. However, its mixture speed of sound $c_{\text{w}}$ is the so-called Wood speed of sound \cite{wood:1930} and is defined through
$$\frac{1}{\rho c_{\text{w}}^2} = \frac{\alpha_{1}}{\rho_{1}c_{1}^{2}} + \frac{\alpha_{2}}{\rho_{2}c_{2}^{2}}.$$
Notice that, as discussed in \cite{haegeman:2024} and in Section \ref{sec:intro}, the relaxation source terms of the 6-equation model relies on an interfacial pressure which is not present in the 5-equation model. We refer to the complementary work \cite{haegeman:2024} for a detailed analysis of the impact of the choice of the interfacial pressure $p_{I}$.

For later reference, we rewrite the 6-equation model with mechanical relaxation in a compact form as
\begin{equation}\label{eq:system_compact}
	\pad{\mathbf{q}}{t} + \dive\mathbf{F}(\mathbf{q}) + \boldsymbol{\sigma}\rpth{\mathbf{q}, \grad\mathbf{q}} = \boldsymbol{\Psi}_{\mu}(\mathbf{q}), 
\end{equation}
where
\begin{subequations}
\begin{align}
	\mathbf{q} &=
	\begin{bmatrix}
		\alpha_{1} \\
		\alpha_{1}\rho_{1} \\
		\alpha_{2}\rho_{2} \\
		\rho\vel \\
		\alpha_{1}\rho_{1}E_{1} \\
		\alpha_{2}\rho_{2}E_{2} 
	\end{bmatrix},
	\quad
	\mathbf{F}(\mathbf{q}) =
	\begin{bmatrix}
		0 \\
		\alpha_{1}\rho_{1}\vel^{\top} \\
		\alpha_{2}\rho_{2}\vel^{\top} \\
		\rho\vel \otimes \vel + \pmix\mathbf{I} \\
		\alpha_{1}\rpth{\rho_{1}E_{1} + p_{1}}\vel^{\top} \\
		\alpha_{2}\rpth{\rho_{2}E_{2} + p_{2}}\vel^{\top}
	\end{bmatrix}, \\
	\boldsymbol{\sigma}\rpth{\mathbf{q}, \grad\mathbf{q}} &=
	\begin{bmatrix}
		\vel \cdot \grad\alpha_{1} \\
		0 \\
		0 \\
		\bm{0} \\
		\vel \cdot \boldsymbol{\Sigma} \\
		-\vel \cdot \boldsymbol{\Sigma},
	\end{bmatrix},
	\quad
	\boldsymbol{\Psi}_{\mu} =
	\begin{bmatrix}
		\mu\rpth{p_{1} - p_{2}} \\
		0 \\
		0 \\
		\bm{0} \\
		-\mu p_{I}\rpth{p_{1} - p_{2}} \\
		\mu p_{I}\rpth{p_{1} - p_{2}}
	\end{bmatrix},
\end{align}
\end{subequations}
with $\boldsymbol{\Sigma}$ as in \eqref{eq:non_conservative_term} and $\mathbf{I}$ denoting the identity tensor.

The numerical discretization of system \eqref{eq:system_compact} is based on an operator splitting approach, similarly to \cite{pelanti:2014, saurel:2009, zein:2010}. More specifically, we first solve over a time interval $\Delta t$ the homogeneous hyperbolic problem
\begin{equation}\label{eq:system_hyperbolic}
	\pad{\mathbf{q}}{t} + \dive\mathbf{F}(\mathbf{q}) + \boldsymbol{\sigma}\rpth{\mathbf{q}, \grad\mathbf{q}} = \mathbf{0}.
\end{equation}
Next, when solving the 5-equation model, the hyperbolic step is followed by the relaxation operator  
\begin{equation}
	\odif{\mathbf{q}}{t} = \boldsymbol{\Psi}_{\mu}(\mathbf{q}),
\end{equation}
in the limit of instantaneous relaxation $\mu\to\infty$. First, we focus on \eqref{eq:system_hyperbolic}, and as such, we set $\mu=0$ until Section \ref{sec:relax}. The goal of the first part of the paper is to compare different numerical strategies for the numerical solution of \eqref{eq:system_hyperbolic} so as to identify robust numerical schemes before considering the instantaneous mechanical relaxation. We recall indeed that neither 6- nor 5- equation model admits a complete set of jump conditions. Hence, we restrict our attention to simple first order schemes both in space and time, so as to avoid any further dependence of higher order methods related to reconstruction, limiters, and large stencils. Moreover, we employ a number of cells that is unfeasible for practical applications with the ultimate goal to highlight the possible limitations of the models when mesh convergence is really established. Then, in the second part of the paper, we consider suitable strategies for the solution of the instantaneous mechanical relaxation and we present the numerical results where instantaneous relaxation is taken into account.

\section{Numerical discretization of the hyperbolic operator}
\label{sec:num_hyp}

This section is dedicated to the numerical discretization of system \eqref{eq:system_hyperbolic}. Since solving system \eqref{eq:system_hyperbolic} is a prerequisite for the numerical solution of the 5-equation model, it is therefore important to assess also the properties of the numerical method for the homogeneous system \eqref{eq:system_hyperbolic} as well as to analyze its interaction with the instantaneous relaxation terms. A comparative study of different numerical methods for the relaxed 6-equation single-velocity two-phase model \eqref{eq:system_compact} has been proposed in \cite{delorenzo:2018}. However, all test cases considered in \cite{delorenzo:2018} have been combined with instantaneous mechanical relaxation, except for one with limited results, thus preventing a clear comparison of the various discretization strategies for the hyperbolic operator.

The discretization of \eqref{eq:system_hyperbolic} employs the finite volume method \cite{leveque:2002, toro:2009}. In this section, we first briefly recall two classical approximate Riemann solvers employed to compute the numerical fluxes. Since the 6-equation model does not admit a complete set of jump conditions, the construction of approximate Riemann solvers is not standard. A detailed discussion on the underlying hypotheses employed to build the approximate Riemann solver is provided. Next, we depict two strategies to obtain a first order method, i.e. the Godunov method \cite{toro:2009} with a suitable treatment of the non-conservative terms, and the wave-propagation method \cite{leveque:2002}.

\subsection{Approximate Riemann solvers}
\label{ssec:Riemann_solvers}

Consider a face separating a left state $\mathbf{q}_{L}$ and a right state $\mathbf{q}_{R}$. We consider simple approximate Riemann solvers \cite{bouchut:2004} for the conservative part, i.e. a sequence of $M + 1$ constant states between $\mathbf{q}_{L}$ and $\mathbf{q}_{R}$ separated by $M$ waves that propagate at velocity $s^{(1)} < s^{(2)} < \dots < s^{(M)}$. We recall that for a conservative system of the form
$$\pad{\mathbf{q}}{t} + \dive\mathbf{F}(\mathbf{q}) = \mathbf{0},$$
the consistency condition is \cite{bouchut:2004, pelanti:2022, toro:2009}
\begin{equation}\label{eq:consistency_flux}
	\sum_{l=1}^{M}s^{(l)}\rpth{\mathbf{q}^{(l+1)} - \mathbf{q}^{(l)}} = \spth{\mathbf{F}(\mathbf{q}_{R}) - \mathbf{F}(\mathbf{q}_{L})}\bm{n}_{L},
\end{equation}
with $\mathbf{q}^{(1)} = \mathbf{q}_{L}$ and $\mathbf{q}^{(M+1)} = \mathbf{q}_{R}$. 

First, we consider the Rusanov flux that consists of two waves moving at speed $s^{(1)} = -s^{(2)} = s > 0$. Hence, it reads as follows \cite{rusanov:1962, toro:2009}
\begin{equation}\label{eq:Rusanov_flux}
	\hat{\mathbf{F}}\rpth{\mathbf{q}_{L}, \mathbf{q}_{R}} = \frac{1}{2}\spth{\mathbf{F}(\mathbf{q}_{L}) + \mathbf{F}(\mathbf{q}_{R})} -
	s\rpth{\mathbf{q}_{L} \otimes \bm{n}_{L} + \mathbf{q}_{R} \otimes \bm{n}_{R}},
\end{equation}
where $\bm{n}_{L}$ and $\bm{n}_{R}$ are the outward unit normal vectors from the element with state $\mathbf{q}_{L}$ and $\mathbf{q}_{R}$, respectively, and
\begin{equation}\label{eq:lambda}
	s = \max\rpth{\left|\vel_{L} \cdot \bm{n}_{L}\right| + c_{\text{f},L}, \left|\vel_{R} \cdot \bm{n}_{R}\right| + c_{\text{f},R}}.
\end{equation}

Next, we consider the HLLC-type solver \cite{toro:2009} presented in \cite{pelanti:2014} and employed also in \cite{delorenzo:2018, pelanti:2022}. The HLLC approximate Riemann solver \cite{toro:2009} consists of three waves moving at speed
\begin{equation}\label{eq:s_HLLC}
	s^{(1)} = s_{L}, \qquad s^{(2)} = s_{\star}, \qquad s^{(3)} = s_{R},
\end{equation}
that separate four constant states $\mathbf{q}_{L}, \mathbf{q}_{\star,L}, \mathbf{q}_{\star,R}, \mathbf{q}_{R}$. In the HLLC approximate Riemann solver, the two external waves are assumed to be shock waves, while the middle wave is a contact discontinuity. Following \cite{davies:1988}, we define
\begin{equation}\label{eq:sL_sR_HLLC}
	s_{L} = \min\rpth{\vel_{L} \cdot \bm{n}_{L} - c_{\text{f},L}, \vel_{R} \cdot \bm{n}_{L} - c_{\text{f},R}}, \qquad s_{R} = \max\rpth{\vel_{L} \cdot \bm{n}_{L} + c_{\text{f},L}, \vel_{R} \cdot \bm{n}_{L} + c_{\text{f},R}}.
\end{equation}
As the mixture quantities satisfy equations similar to the Euler equations, the intermediate speed is determined as \cite{pelanti:2022, toro:2009}
\begin{equation}\label{eq:s_star_HLLC}
	s_{\star} = \frac{\pmix_{R} - \pmix_{L} + \rho_{L}\rpth{\vel_{L} \cdot \bm{n}_{L}}\rpth{s_{L} - \vel_{L} \cdot \bm{n}_{L}} - \rho_{R}\rpth{\vel_{R} \cdot \bm{n}_{L}}\rpth{s_{R} - \vel_{R} \cdot \bm{n}_{L}}}{\rho_{L}\rpth{s_{L} - \vel_{L} \cdot \bm{n}_{L}} - \rho_{R}\rpth{s_{R} - \vel_{R} \cdot \bm{n}_{L}}}.
\end{equation}
Hence, the HLLC flux reads as follows
\begin{equation}\label{eq:HLLC_flux}
	\hat{\mathbf{F}}\rpth{\mathbf{q}_{L}, \mathbf{q}_{R}} =
	\begin{cases}
		\mathbf{F}(\mathbf{q}_{L}) \qquad &\text{if } 0 \le s_{L}, \\
		\mathbf{F}(\mathbf{q}_{L}) + s_{L}\rpth{\mathbf{q}_{\star,L} - \mathbf{q}_{L}} \otimes \bm{n}_{L} \qquad &\text{if } s_{L} < 0 \le s_{\star}, \\
		\mathbf{F}(\mathbf{q}_{R}) + s_{R}\rpth{\mathbf{q}_{\star,R} - \mathbf{q}_{R}} \otimes \bm{n}_{L} \qquad &\text{if } s_{\star} < 0 \le s_{R}, \\
		\mathbf{F}(\mathbf{q}_{R}) \qquad &\text{if } 0 > s_{R}.
	\end{cases}
\end{equation}
For the conserved variables, the determination of the intermediate state is given by the associated jump conditions. Moreover, since the volume fraction is constant along the fluid trajectories, one immediately obtains that it is constant through shocks and therefore
\begin{equation}
	\alpha_{1_{\star,L}} = \alpha_{1_{L}}, \qquad \alpha_{1_{\star,R}} = \alpha_{1_{R}}. 
\end{equation}
Approximate jump conditions must be provided to determine the remaining variables in the intermediate region. We refer to the upcoming Section \ref{ssec:approximate_jump_relations} for a detailed discussion on how the remaining variables of $\mathbf{q}_{\star,L}$ and $\mathbf{q}_{\star,R}$ are determined. It results that the middle states read as follows \cite{pelanti:2022, pelanti:2014}
\begin{equation}\label{eq:q_star_HLLC}
	\mathbf{q}_{\star,\#} =
	\begin{bmatrix}
		\alpha_{1_{\#}} \\
		\alpha_{1_{\#}}\rho_{1_{\#}}\frac{s_{\#} - \vel_{\#} \cdot \bm{n}_{L}}{s_{\#} - s_{\star}} \\
		\alpha_{2_{\#}}\rho_{2_{\#}}\frac{s_{\#} - \vel_{\#} \cdot \bm{n}_{L}}{s_{\#} - s_{\star}} \\
		\rho_{\#}\frac{s_{\#} - \vel_{\#} \cdot \bm{n}_{L}}{s_{\#} - s_{\star}}\rpth{\vel_{\#} + \rpth{s_{\star} - \vel_{\#} \cdot \bm{n}_{L}}\bm{n}_{L}} \\
		\alpha_{1_{\#}}\rho_{1_{\#}}\frac{s_{\#} - \vel_{\#} \cdot \bm{n}_{L}}{s_{\#} - s_{\star}}E_{1_{\star,\#}} \\
		\alpha_{2_{\#}}\rho_{2_{\#}}\frac{s_{\#} - \vel_{\#} \cdot \bm{n}_{L}}{s_{\#} - s_{\star}}E_{2_{\star,\#}}
	\end{bmatrix}, \qquad \# = L,R.
\end{equation}
Notice that, in the above expression \eqref{eq:q_star_HLLC}, $\rho_{\#}\frac{s_{\#} - \vel_{\#} \cdot \bm{n}_{L}}{s_{\#} - s_{\star}}\rpth{\vel_{\#} + \rpth{s_{\star} - \vel_{\#} \cdot \bm{n}_{L}}\bm{n}_{L}}$ reduces to $\rho_{\#}\frac{s_{\#} - \vel_{\#} \cdot \bm{n}_{L}}{s_{\#} - s_{\star}}s_{\star}$ in the direction normal to the interface, and to $\rho_{\#}\frac{s_{\#} - \vel_{\#} \cdot \bm{n}_{L}}{s_{\#} - s_{\star}}\rpth{\vel_{\#} - \rpth{\vel_{\#} \cdot \bm{n}_{L}}\bm{n}_{L}}$ for the tangential components of the velocity. In expression \eqref{eq:q_star_HLLC}, all terms are closed except for $E_{k_{\star,\#}}, k=1,2$ and $\#=L,R$, as no jump conditions are available for the phasic energy equations. We will now discuss the approximate jump conditions that are provided to close the Riemann solver. These conditions relate the energy jumps to pressure and density jumps, leading to the determination of both the internal energies $e_{k_{\star,\#}}$ and the phasic pressures $p_{k_{\star,\#}}$, and ultimately the remaining unknowns $E_{k_{\star,\#}}$.

\subsection{Discussion on approximate jump relations for the mixture-energy-consistent 6-equation model}
\label{ssec:approximate_jump_relations}

In the previous section, the intermediate states of the HLLC Riemann solver are given by \eqref{eq:q_star_HLLC}. However, their expressions are not fully closed since the system lacks a complete set of jump conditions. This section is devoted to the closures of the remaining terms through approximate jump conditions relating the energy jumps to those of pressures and densities. Several approaches can be found in the literature, in this section we will show that they can be seen as different numerical discretizations of the same approximate jump conditions, thus bridging the various approaches. From now on, for a generic variable $\varphi$, the symbol $\ave{\varphi}$ denotes the arithmetic average between the left value and the right value, i.e.
$$\ave{\varphi} = \frac{1}{2}\rpth{\varphi_{L} + \varphi_{R}},$$
whereas the symbol $\jump{\varphi}$ denotes the jump between the left value and the right value, i.e.
$$\jump{\varphi} = \varphi_{L} - \varphi_{R}.$$
We focus on the one-dimensional setting for the sake of simplicity. For system \eqref{eq:system_hyperbolic}, the jump relations through a shock are only mathematically defined for the conserved equations and the volume fraction equation, and they can be written as ($k = 1,2$)
\begin{subequations}
	\begin{align}
		\jump{\alpha_{k}} &= 0, \label{eq:jump_relation_vol_fraction} \\
		\jump{\alpha_{k}\rho_{k}\rpth{u-S}} &= 0, \label{eq:jump_relation_densities} \\
		\jump{\rho u\rpth{u-S}} + \jump{\pmix} &= 0, \label{eq:jump_relation_momentum} 
	\end{align}
\end{subequations}
with $S$ denoting the speed of the shock. In particular, these relations lead to closed expressions of the intermediate state in the HLLC approximate Riemann solver \eqref{eq:q_star_HLLC} for the corresponding variables. The variables for which the intermediate state \eqref{eq:q_star_HLLC} is still unclosed are those for which no jump conditions are readily available, i.e. the phasic energies, and thus jump conditions must be provided.

The following approximate jump conditions for the 5-equation model \eqref{eq:Kapila_system} have been proposed \cite{saurel:2007}
\begin{equation}\label{eq:jump_relation_5-eq}
	\jump{e_{k}} + \ave{p}\jump{\frac{1}{\rho_{k}}} = 0, \qquad k=1,2. 
\end{equation}
We stress the fact that in the 5-equation model \eqref{eq:Kapila_system}, the phases are at pressure equilibrium at all times and thus $p$ denotes their common pressure. Keeping in mind that one of the main uses of the 6-equation model is to serve as an auxiliary system for the resolution of the 5-equation model \eqref{eq:Kapila_system}, a first approach, that has been proposed in \cite{saurel:2009}, is to extend the jump condition \eqref{eq:jump_relation_5-eq} to the 6-equation model as follows
\begin{equation}\label{eq:jump_relation_eint_Saurel}
	\jump{e_{k}} + \ave{p_{k}}\jump{\frac{1}{\rho_{k}}} = 0, \qquad k=1,2. 
\end{equation}
These relations are compatible with the volume fraction jump conditions \eqref{eq:jump_relation_vol_fraction} as well as the mixture energy's jump condition, which can be deduced from \eqref{eq:conservation_total_energy}. It also recovers \eqref{eq:jump_relation_5-eq} when pressure equilibrium is imposed. However, this approach raises the question of the validity of the proposed jump conditions \eqref{eq:jump_relation_eint_Saurel} far from pressure equilibrium, namely in the case where finite-rate pressure relaxations are considered. The link between this approach and another approach for the numerical solution of 6-equation model that we are going to present will provide an insight on this question. Relation \eqref{eq:jump_relation_eint_Saurel} is not yet sufficient to close the intermediate states in the HLLC approximate Riemann solver as it provides a single equation for each phase, for two unknowns, e.g., $e_{k_{\star,L}}$ and $p_{k_{\star,L}}$ for the left intermediate state. For this approach, it is proposed in \cite{saurel:2009} to use the phasic equation of state and to solve
\begin{equation}\label{eq:closure_thermo_Saurel}
	e_{k_{\star,L}} - e_{k,L} + \frac{1}{2}\rpth{p_{k_{\star,L}}\rpth{e_{k_{\star,L}}, \rho_{k_{\star,L}}} + p_{k,L}}\rpth{\frac{1}{\rho_{k_{\star,L}}} - \frac{1}{\rho_{k,L}}} = 0.
\end{equation}
However, this results in nonlinear thermodynamic equations to be solved for each face of each computational cell. Although it admits an analytical solution in the case of stiffened gas laws \eqref{eq:SG_EOS}, it is not the case in general. Moreover, when doing so, in general we do not have that $\alpha_{1_{\star,L}} p_{1_{\star,L}} + \alpha_{2_{\star,L}} p_{2_{\star,L}} = \alpha_{1_{\star,R}} p_{1_{\star,R}} + \alpha_{2_{\star,R}} p_{2_{\star,R}}$, and neither of these expressions are in general equal to $\pmix_{\star}$ which can be determined using the mixture momentum equation. It then results that the mixture total energy and the mixture momentum rely on different expressions of the mixture pressure in the intermediate regions.

A second approach has been proposed in \cite{pelanti:2014}. It consists in neglecting the contribution of the non-conservative products in the phasic total energy equations through shocks. Hence, it results in the following approximate jump condition through shock waves
\begin{equation}\label{eq:jump_relation_phasic_energies}
	\jump{\alpha_{k}\rho_{k}E_{k}\rpth{u-S}} + \jump{\alpha_{k}p_{k}u} = 0, \qquad k=1,2.
\end{equation}
However, once again, this is not sufficient to yield a closed solution for the HLLC's intermediate states, since, just as before, it introduces two unknowns, e.g., $E_{k_{\star,L}}$ and $p_{k_{\star,L}}$ for a single equation. Following \cite{pelanti:2014}, instead of calling upon the equation of state, the phasic momentum equations are considered
\begin{subequations}
	\begin{align}
		\pad{\alpha_{1}\rho_{1}u}{t} + \pad{\rpth{\alpha_{1}\rho_{1}u^{2} + \alpha_{1}p_{1}}}{x} + u\Sigma &= 0, \\
		\pad{\alpha_{2}\rho_{2}u}{t} + \pad{\rpth{\alpha_{2}\rho_{2}u^{2} + \alpha_{2}p_{2}}}{x} - u\Sigma &= 0.
	\end{align}
\end{subequations}
Neglecting also their non-conservative contributions through shocks yields
\begin{equation}\label{eq:jump_relations_phasic_momenta}
	\jump{\rho_{k}u\rpth{u-S}} + \jump{p_{k}} = 0, \qquad k=1,2,
\end{equation}
where we have simplified the jump relation by $\alpha_k$ owing to \eqref{eq:jump_relation_vol_fraction}. These are of course compatible with the jump relation \eqref{eq:jump_relation_momentum} for the mixture momentum and provide expressions for the phasic pressures in the intermediate regions. Together with \eqref{eq:jump_relation_phasic_energies}, it allows to close the intermediate states of the HLLC approximate Riemann solver \eqref{eq:q_star_HLLC} as it results in
\begin{equation}\label{fermeture_E_kstar_Pelanti}
	E_{k_{\star,\#}} = E_{k_{\#}} + \rpth{s_{\star} - \vel_{\#} \cdot \bm{n}_{L}}\rpth{s_{\star} + \frac{p_{k_{\#}}}{\rho_{k_{\#}}\rpth{s_{\#} - \vel_{\#} \cdot \bm{n}_{L}}}},\qquad k=1,2,\quad \#=L,R.
\end{equation}
Contrary to the previous approach, here a closed analytical expression is obtained. Moreover, it is compatible with mixture pressure in the intermediate region as the condition
$$ 
\alpha_{1_{\star,L}} p_{1_{\star,L}} + \alpha_{2_{\star,L}} p_{2_{\star,L}} = \pmix_\star =\alpha_{1_{\star,R}} p_{1_{\star,R}} + \alpha_{2_{\star,R}} p_{2_{\star,R}},
$$
is now automatically satisfied, thus also ensuring compatibility with the mixture's total energy equation.

It is important to note that in this approach, the non-conservative products have only been neglected through the shocks, i.e. the external waves of the HLLC Riemann solver. Moreover, the resulting intermediate states result in a consistent discretization through contact discontinuities of the system, including its non-conservative products as shown in \cite{pelanti:2022}. Indeed, in 1D, the non-conservative product for the phasic total energies writes
$$u \Sigma = -u \rpth{Y_{2}\pad{\alpha_{1}p_{1}}{x} - Y_{1}\pad{\alpha_{2}p_{2}}{x}} = -u \pad{\alpha_{1}p_{1}}{x} + Y_{1} u \pad{\pmix}{x}.$$
Through a contact discontinuity, the mixture pressure $\pmix$ and velocity $u$ are both constant with the velocity being equal to $u = s_{\star}$, the speed of the moving discontinuity. It then results that, through the contact discontinuity, the non-conservative term $u \Sigma$ writes

$$u \Sigma = -s_\star \rpth{ \rpth{\alpha_{1}p_{1}}_{\star,R} - \rpth{\alpha_{1}p_{1}}_{\star,L}}.$$

The intermediate total energy of phase 1 then writes
\begin{align*}
	&s_{\star}\rpth{\rpth{\alpha_{1}\rho_{1}E_{1}}_{\star,R} - \rpth{\alpha_{1}\rho_{1}E_{1}}_{\star,L}} = \\ &\spth{\rpth{\alpha_{1}\rho_{1}E_{1}u  + \alpha_{1}p_{1}u}_{\star,R} - \rpth{\alpha_{1}\rho_{1}E_{1}u  + \alpha_{1}p_{1}u}_{\star,L}} - s_{\star}\rpth{ \rpth{\alpha_{1}p_{1}}_{\star,R} - \rpth{\alpha_{1}p_{1}}_{\star,L}}.
\end{align*}
The left-hand side is a consistent discretization of the time-derivative of the phasic total energy through the contact discontinuity, and on the right-hand side we obtain a consistent discretization of the conservative flux in between the square brackets while the remaining term is exactly the non-conservative product that we have just determined. The same result holds for the second phase up to the correct change of signs. \\

The two approaches we have presented have both been used in the literature, however a link between them has not been established yet. We now show how they are related and can be seen as two different discretization strategies of the same approximate jump conditions, namely \eqref{eq:jump_relation_eint_Saurel}. Indeed, starting from the second approach with the assumption that the non-conservative products in the phasic total energy equations can be neglected through shocks, we have obtained \eqref{eq:jump_relation_phasic_energies}. We develop this relation using the identity $\jump{\phi\psi} = \ave{\phi}\jump{\psi} + \ave{\psi}\jump{\phi}$ so as to obtain
$$\tilde{m}_{k}\jump{e_{k}} + \jump{\tilde{m}_{k}u}\ave{u} + \jump{\alpha_{k}p_{k}}\ave{u} + \ave{\alpha_{k}p_{k}}\jump{u}=0.$$
Here we have denoted $\tilde{m}_{k} = \alpha_{k}\rho_{k}\rpth{u - S}$ the phasic mass flux, which owing to the mass conservation \eqref{eq:jump_relation_densities} is continuous through the shock. Regrouping the two middle terms and the two remaining terms, we have, owing to \eqref{eq:jump_relation_vol_fraction} and dividing by $\alpha_{k}$, that
$$\rpth{\rho_{k}\rpth{u-S}\jump{e_{k}} + \ave{p_{k}}\jump{u}} + \ave{u}\rpth{\jump{\rho_{k}u\rpth{u-S}} + \jump{p_{k}}}=0.$$
Since $u = S + \tilde{m}_{k}/\rpth{\alpha_{k}\rho_{k}}$, we have $\jump{u} = \tilde{m}_{k}\jump{1/\rpth{\alpha_{k}\rho_{k}}} = \rho_{k}\rpth{u-S}\jump{1/\rho_{k}}$, and thus we obtain
\begin{equation}\label{eq:jump_relation_internal_energies_6eqs}
	\rpth{\rho_{k}\rpth{u-S}}\rpth{\jump{e_{k}} + \ave{p_{k}}\jump{\frac{1}{\rho_{k}}}} + \ave{u}\rpth{\jump{\rho_{k}u\rpth{u-S}}+\jump{p_{k}}} = 0.
\end{equation}
In order to close the intermediate states of the HLLC solver, in the second approach we have also assumed \eqref{eq:jump_relations_phasic_momenta}, which owing to \eqref{eq:jump_relation_internal_energies_6eqs} is now shown to also yield \eqref{eq:jump_relation_eint_Saurel}. It then results that both approaches can be described in the following unified manner. First, at the continuous level, we provide the additional jump conditions \eqref{eq:jump_relation_eint_Saurel}
$$\jump{e_{k}} + \ave{p_{k}}\jump{\frac{1}{\rho_{k}}} = 0, \qquad k=1,2. $$
At the discrete level, these are not sufficient to close the intermediate states of the HLLC since, for each phase, it provides a single equation with two unknowns. In the first approach, proposed in \cite{saurel:2009}, one closes the HLLC states using the equation of state leading to \eqref{eq:closure_thermo_Saurel}, while in the second approach, one assumes the additional relation \eqref{eq:jump_relations_phasic_momenta} resulting in the analytical expressions \eqref{fermeture_E_kstar_Pelanti}. Considering the previous discussions, we favour the second approach which we use to close \eqref{eq:q_star_HLLC}. In particular, it imposes that the mixture pressure is continuous through a material interface, which is a property of the exact Riemann solution of the 6-equation model. Moreover, this approach ensures that
$$\pmix_{\star,\#} = \sum\limits_{k=1,2}\alpha_{k_{\star,\#}}p_{k_{\star,\#}},$$
which allows for a coherent discretization of the phasic energy equations with respect to the momentum equation. Note also that, for the classical HLLC scheme applied to the Euler equations, the equation of state is not used to determine the intermediate states \cite{toro:2009}, as the computed flux is a mechanical response to the imposed shock velocities, which ensures integral consistency with the equations. In particular, the pressure in the intermediate regions are, in general, not equal to thermodynamic pressures obtained from the densities and internal energies of those regions.

\subsection{Treatment of non-conservative terms}
\label{ssec:non_conservative}

In this section, we analyze the numerical treatment of the non-conservative terms. Since the system is non-conservative, the numerical flux consists in general of a conservative part, $\hat{\mathbf{F}}$, and an additional term for non-conservative terms, which is face-dependent and will be denoted by $\boldsymbol{\mathcal{T}^{\pm}}$. Non-conservative hyperbolic systems arise in a wide range of applications, e.g. two-phase flows \cite{baer:1986, kapila:2001, murrone:2005}, shallow water equations \cite{bouchut:2008, rosatti:2011}, and plasma physics \cite{wargnier:2020}. Several approaches have been proposed in the literature to address the discretization of non-conservative terms. A popular strategy is based on the Dal Maso, Le Floch and Murat (DLM) theory \cite{dalmaso:1995} (see also \cite{bianchini:2005}). Matching at the discrete level the DLM theory, the concept of path-conservative schemes for the numerical solution of non-conservative hyperbolic systems has been proposed in \cite{pares:2006}. Several path-conservative Riemann solvers for non-conservative hyperbolic systems have been designed \cite{delorenzo:2018, dumbser:2016, dumbser:2011}. The main issue with path-conservative schemes is that, even for conservative systems, when rewritten in a non-conservative form, the path-conservative scheme may not converge to the specified entropic weak solution in the case of discontinuities \cite{castro:2008, hou:1994}, even when the path is chosen such that, at the continuous level, the path-integral yields the exact jump conditions \cite{abgrall:2010}. However, this is a common issue of all non-conservative schemes when applied to conservation laws \cite{hou:1994}, unless some additional corrections are included \cite{abgrall:2018, chalons:2017}. As such, the path-conservative approach is mostly formal, still it remains quite popular and yields reasonable results for several non-conservative systems of equations \cite{busto:2021, dumbser:2011, nguyen:2015, tokareva:2010}

Other approaches to treat numerically non-conservative terms have been proposed in the literature. We refer in particular to \cite{abgrall:2018}, where a high-order approach for non-conservative systems based on the so-called residual distribution formulation is proposed (see also \cite{sirianni:2025, sirianni:2024}). In the present work, we consider and we compare two simpler strategies for the approximation of non-conservative terms, which do not rely on Riemann solvers. The first strategy is inspired by a Bassi-Rebay (BR)-type approach \cite{bassi:1997a}, originally developed for the discretization of diffusion operators and later extended to the discretization of non-conservative terms in hyperbolic systems, see, e.g., \cite{orlando:2023, tumolo:2015}. The term BR-type is used here to refer to a specific algebraic structure based on two integrations by parts (see Appendix \ref{app:BR_non_cons} for further details), rather than to claim any priority over the broad finite volume literature on consistent gradient discretizations. We refer the reader to, e.g., \cite{eymard:2000} for a comprehensive overview of gradient discretizations in finite volume methods, including a detailed comparison with mixed formulations, which share ideas with those presented in \cite{bassi:1997a}. Consider the generic non-conservative product $\mathbf{B}(w)\grad w$, where $w$ is a scalar field. Assuming $K$ is the cell containing the left state $L$ such that the out-going normal at each interface points towards the right state $R$, we obtain \cite{orlando:2023}
\begin{align}\label{eq:BR}
	\int_{K} \mathbf{B}(w)\grad w\mathrm{d}\Omega &= \int_{K} \dive\rpth{\mathbf{B}(w)w}\mathrm{d}\Omega - \int_{K} w \dive\mathbf{B}(w) \mathrm{d}\Omega \nonumber \\
	&\approx \sum_{\Gamma \in \partial K}\int_{\Gamma} \rpth{\widehat{\mathbf{B}w}(w_{L}, w_{R}) - \widehat{\mathbf{B}}(w_{L}, w_{R})w_{L}}\bm{n}\mathrm{d}\Gamma,
\end{align}
where the approximations $\widehat{\mathbf{B}w}(w_{L}, w_{R})$ and $\widehat{\mathbf{B}}(w_{L}, w_{R})w_{L}$ have still to be defined. We consider here the approximations adopted in \cite{orlando:2023} (we refer to it as \textit{BR-2023})
\begin{equation}\label{eq:fluxes_Orlando}
	\widehat{\mathbf{B}w}(w_{L}, w_{R}) = \frac{1}{2}\rpth{\mathbf{B}(w_{L})w_{L} + \mathbf{B}(w_{R})w_{R}}, \quad \widehat{\mathbf{B}}(w_{L}, w_{R}) = \frac{1}{2}\rpth{\mathbf{B}(w_{L}) + \mathbf{B}(w_{R})},
\end{equation}
and those employed in \cite{tumolo:2015} (we refer to it as \textit{BR-2015})
\begin{equation}\label{eq:fluxes_Tumolo}
	\widehat{\mathbf{B}w}(w_{L}, w_{R}) = \frac{1}{4}\rpth{\mathbf{B}(w_{L}) + \mathbf{B}(w_{R})}\rpth{w_{L} + w_{R}}, \quad \widehat{\mathbf{B}}(w_{L}, w_{R}) = \frac{1}{2}\rpth{\mathbf{B}(w_{L}) + \mathbf{B}(w_{R})}.
\end{equation}
Hence, we consider two different discretizations of the non-conservative terms of system \eqref{eq:system_hyperbolic}. First, employing \eqref{eq:fluxes_Orlando}, we set for future reference
\begin{equation}\label{eq:BR_Orlando_non_cons}
	\boldsymbol{\mathcal{T}}^{-}\rpth{\mathbf{q}_{L}, \mathbf{q}_{R}} =
	\begin{bmatrix}
		\rpth{\ave{\vel\alpha_{1}} - \ave{\vel}\alpha_{1,L}} \cdot \bm{n}_{L} \\
		0 \\
		0 \\
		0 \\
		+\hat{\boldsymbol{\Sigma}}^{-} \cdot \bm{n}_{L} \\
		-\hat{\boldsymbol{\Sigma}}^{-} \cdot \bm{n}_{L}
	\end{bmatrix}, \\
	\boldsymbol{\mathcal{T}}^{+}\rpth{\mathbf{q}_{L}, \mathbf{q}_{R}} =
	\begin{bmatrix}
		\rpth{\ave{\vel\alpha_{1}} - \ave{\vel}\alpha_{1,R}} \cdot \bm{n}_{L} \\
		0 \\
		0 \\
		0 \\
		+\hat{\boldsymbol{\Sigma}}^{+} \cdot \bm{n}_{L} \\
		-\hat{\boldsymbol{\Sigma}}^{+} \cdot \bm{n}_{L}
	\end{bmatrix}.
\end{equation}
Moreover, we set
\begin{subequations}
	\begin{align}
		\hat{\boldsymbol{\Sigma}}^{-} &= -\rpth{\ave{\vel Y_{2}\alpha_{1}p_{1}} - \ave{\vel Y_{2}}\alpha_{1,L}p_{1,L}} + \rpth{\ave{\vel Y_{1}\alpha_{2}p_{2}} - \ave{\vel Y_{1}}\alpha_{2,L}p_{2,L}}, \\
		\hat{\boldsymbol{\Sigma}}^{+} &= -\rpth{\ave{\vel Y_{2}\alpha_{1}p_{1}} - \ave{\vel Y_{2}}\alpha_{1,R}p_{1,R}} + \rpth{\ave{\vel Y_{1}\alpha_{2}p_{2}} - \ave{\vel Y_{1}}\alpha_{2,R}p_{2,R}}. 
	\end{align}
\end{subequations}
Second, employing \eqref{eq:fluxes_Tumolo}, we define
\begin{equation}\label{eq:BR_Tumolo_non_cons}
	\boldsymbol{\mathcal{T}}^{-}\rpth{\mathbf{q}_{L}, \mathbf{q}_{R}} =
	\begin{bmatrix}
		\rpth{\ave{\vel}\ave{\alpha_{1}} - \ave{\vel}\alpha_{1,L}} \cdot \bm{n}_{L} \\
		0 \\
		0 \\
		0 \\
		+\bar{\boldsymbol{\Sigma}}^{-} \cdot \bm{n}_{L} \\
		-\bar{\boldsymbol{\Sigma}}^{-} \cdot \bm{n}_{L}
	\end{bmatrix},
	\boldsymbol{\mathcal{T}}^{+}\rpth{\mathbf{q}_{L}, \mathbf{q}_{R}} =
	\begin{bmatrix}
		\rpth{\ave{\vel}\ave{\alpha_{1}} - \ave{\vel}\alpha_{1,R}} \cdot \bm{n}_{L} \\
		0 \\
		0 \\
		0 \\
		+\bar{\boldsymbol{\Sigma}}^{+} \cdot \bm{n}_{L} \\
		-\bar{\boldsymbol{\Sigma}}^{+} \cdot \bm{n}_{L}
	\end{bmatrix},
\end{equation}
where
\begin{subequations}
	\begin{align}
		\bar{\boldsymbol{\Sigma}}^{-} = &-\rpth{\ave{\vel Y_{2}}\ave{\alpha_{1}p_{1}} - \ave{\vel Y_{2}}\alpha_{1,L}p_{1,L}} + \rpth{\ave{\vel Y_{1}}\ave{\alpha_{2}p_{2}} - \ave{\vel Y_{1}}\alpha_{2,L}p_{2,L}}, \\
		\bar{\boldsymbol{\Sigma}}^{+} = &-\rpth{\ave{\vel Y_{2}}\ave{\alpha_{1}p_{1}} - \ave{\vel Y_{2}}\alpha_{1,R}p_{1,R}} + \rpth{\ave{\vel Y_{1}}\ave{\alpha_{2}p_{2}} - \ave{\vel Y_{1}}\alpha_{2,R}p_{2,R}}.
	\end{align}
\end{subequations}

The second strategy that we consider is the one presented in \cite{crouzet:2013}, so that
\begin{equation}
	\boldsymbol{\mathcal{T}}^{-}\rpth{\mathbf{q}_{L}, \mathbf{q}_{R}} =
	\begin{bmatrix}
		\vel_{L}\ave{\alpha_{1}} \cdot \bm{n}_{L} \\
		0 \\
		0 \\
		0 \\
		\tilde{\boldsymbol{\Sigma}}^{-} \cdot \bm{n}_{L} \\
		-\tilde{\boldsymbol{\Sigma}}^{-} \cdot \bm{n}_{L}
	\end{bmatrix}, \quad
	\boldsymbol{\mathcal{T}}^{+}\rpth{\mathbf{q}_{L}, \mathbf{q}_{R}} =
	\begin{bmatrix}
		\vel_{R}\ave{\alpha_{1}} \bm{n}_{L} \\
		0 \\
		0 \\
		0 \\
		\tilde{\boldsymbol{\Sigma}}^{+} \cdot \bm{n}_{L} \\
		-\tilde{\boldsymbol{\Sigma}}^{+} \cdot \bm{n}_{L}
	\end{bmatrix}, \label{eq:Crouzet_non_cons}
\end{equation}
where
\begin{subequations}
	\begin{align}
		\tilde{\boldsymbol{\Sigma}}^{-} &= -\vel_{L}\rpth{Y_{2,L}\ave{\alpha_{1}p_{1}} - Y_{1,L}\ave{\alpha_{2}p_{2}}}, \\
		\tilde{\boldsymbol{\Sigma}}^{+} &= -\vel_{R}\rpth{Y_{2,R}\ave{\alpha_{1}p_{1}} - Y_{1,R}\ave{\alpha_{2}p_{2}}}.
	\end{align}
\end{subequations}

The strategies discussed so far treat the non-conservative terms separately from the conservative contribution. However, there exist solvers that can naturally be applied to non-conservative systems. Examples are Suliciu-type relaxation schemes, as the one developed in \cite{delorenzo:2018}, the Roe method \cite{pelanti:2012, roe:1981} or the VFRoe-ncv scheme \cite{murrone:2005}, or the wave-propagation method \cite{ketcheson:2013, leveque:2002}, that we will also employ in our numerical experiments and that we will briefly recall in the one-dimensional setting in Section \ref{ssec:disc_scheme}.

Finally, the treatment of the term $\vel \cdot \grad\alpha_{1}$ in the equation of the volume fraction deserves some comments. In the classical Rusanov formulation \eqref{eq:Rusanov_flux}, the entire vector of variables $\mathbf{q}$ is used, which gives rise to the stabilizing term
$$s\rpth{\mathbf{q}_{L} \otimes \bm{n}_{L} + \mathbf{q}_{R} \otimes \bm{n}_{R}},$$
where $s$ is defined as in \eqref{eq:lambda}. This term therefore also accounts for jumps in the volume fraction and is essential to obtain physically meaningful solutions when discretizations \eqref{eq:BR_Orlando_non_cons}, \eqref{eq:BR_Tumolo_non_cons}, and \eqref{eq:Crouzet_non_cons} are combined with the Rusanov flux \eqref{eq:Rusanov_flux}. The discretization of the non-conservative terms alone does not provide this stabilization. However, the diffusive term of the Rusanov flux is sufficient to stabilize the scheme. If the same discretizations are instead coupled with the HLLC-type conservative flux \eqref{eq:HLLC_flux}, the solution becomes severely corrupted, which further highlights the critical role of this stabilizing contribution. Although HLLC is more accurate, it does not provide sufficient numerical diffusion to ensure stability in this context. Consequently, when using HLLC-type approximate Riemann solvers, we apply the aforementioned discretization framework only to the non-conservative terms in the energy equations, while we use upwinding for the volume fraction equation. Since, \eqref{eq:volume_fraction_advection} is simply an advection equation, jump conditions are well-defined and the volume fraction can be discontinuous only along a contact discontinuity \cite{toro:2009}. Hence, following, e.g., \cite{lochon:2016}, when using the HLLC flux for the conservative portion, we apply the following discretization of the non-conservative contribution that acts as a stabilization for the HLLC flux
\begin{equation}
	\boldsymbol{\mathcal{T}}^{-}\rpth{\mathbf{q}_{L}, \mathbf{q}_{R}} =
	\begin{bmatrix}
		\frac{1 - \text{sign}(s_{\star})}{2}\rpth{\alpha_{1,R} - \alpha_{1,L}} \\
		0 \\
		0 \\
		0 \\
		\square^{-} \cdot \bm{n}_{L} \\
		-\square^{-} \cdot \bm{n}_{L}
	\end{bmatrix}, \quad
	\boldsymbol{\mathcal{T}}^{+}\rpth{\mathbf{q}_{L}, \mathbf{q}_{R}} =
	\begin{bmatrix}
		-\frac{1 + \text{sign}(s_{\star})}{2}\rpth{\alpha_{1,R} - \alpha_{1,L}} \\
		0 \\
		0 \\
		0 \\
		\square^{+} \cdot \bm{n}_{L} \\
		-\square^{+} \cdot \bm{n}_{L}
	\end{bmatrix}, \label{eq:HLLC_non_cons}
\end{equation}
where $\square \in \cpth{\tilde{\boldsymbol{\Sigma}}, \bar{\boldsymbol{\Sigma}}, \hat{\boldsymbol{\Sigma}}}$. We have presented here several, among many, possible approaches for the direct treatment of the non-conservative products. More specifically, we have recalled the \textit{BR-2015} scheme, for which the numerical approximations at the two faces of the boundary are given by \eqref{eq:BR_Orlando_non_cons}; the \textit{BR-2023} approach, for which the numerical approximations at the two faces of the boundary are given by \eqref{eq:BR_Tumolo_non_cons}; and the discretization \eqref{eq:Crouzet_non_cons} proposed in \cite{crouzet:2013}, which we will denote as \Crouzet.\ These discretizations are appropriately modified according to \eqref{eq:HLLC_non_cons} when coupled with the HLLC-type conservative flux \eqref{eq:HLLC_flux}. We will show in Section \ref{sec:BR_path_conservative} that the approaches recalled here can be framed within the unified formalism of path-conservative schemes.

\subsection{Fully discrete scheme}
\label{ssec:disc_scheme}

In this section, we present the fully discrete numerical method for \eqref{eq:system_hyperbolic}. For the sake of simplicity in the notation, we consider the one-dimensional setting. We assume therefore a spatial discretization on a mesh with cells of uniform size $\Delta x$. We denote by $\mathbf{q}_{j}$ the approximate solution of the system at the cell $j$. Following the classical Godunov method \cite{toro:2009}, we obtain
\begin{equation}\label{eq:discrete_scheme}
	\mathbf{q}_{j}^{n+1} = \mathbf{q}_{j}^{n} - \frac{\Delta t}{\Delta x}\rpth{\boldsymbol{\mathcal{H}}^{-}_{j + 1/2}\rpth{\mathbf{q}_{j}^{n}, \mathbf{q}_{j+1}^{n}} - \boldsymbol{\mathcal{H}}^{+}_{j - 1/2}\rpth{\mathbf{q}_{j-1}^{n}, \mathbf{q}_{j}^{n}}}.
\end{equation}
Since system \eqref{eq:system_hyperbolic} is not conservative, $\boldsymbol{\mathcal{H}}^{-} \neq \boldsymbol{\mathcal{H}}^{+}$ in general. More specifically, we consider the numerical fluxes described in Section \ref{ssec:Riemann_solvers} for the conservative part, whereas the strategies depicted in Section \ref{ssec:non_conservative} are employed for the discretization of the non-conservative terms. Hence, we obtain
\begin{align}
	\boldsymbol{\mathcal{H}}^{-} &= \hat{\mathbf{F}} + \boldsymbol{\mathcal{T}}^{-}, \\
	\boldsymbol{\mathcal{H}}^{+} &= \hat{\mathbf{F}} + \boldsymbol{\mathcal{T}}^{+},
\end{align}
where we recall that $\hat{\mathbf{F}}$ denotes the discretization of the conservative portion of system \eqref{eq:system_hyperbolic}, whereas $\boldsymbol{\mathcal{T}}$ denotes the discretization of the non-conservative terms.

We also consider the wave-propagation method \cite{leveque:2002}, which provides an elegant and effective scheme for the discretization of systems including (possibly) non-conservative terms. The one-dimensional first order wave-propagation method reads \cite{pelanti:2014}
\begin{equation}\label{eq:wave_propagation}
	\mathbf{q}_{j}^{n+1} = \mathbf{q}_{j}^{n} - \frac{\Delta t}{\Delta x}\rpth{\boldsymbol{\mathcal{A}}^{-}\Delta\mathbf{q}_{j+1/2} - \rpth{-\boldsymbol{\mathcal{A}}^{+}\Delta\mathbf{q}_{j-1/2}}},
\end{equation}
where $\boldsymbol{\mathcal{A}}^{\pm}\Delta\mathbf{q}_{j\mp1/2}$ are the so-called fluctuations at the interface between cells $j$ with state variables $\mathbf{q}_{j}$ and $j+1$ with state variables $\mathbf{q}_{j+1}$. Note that the formulation \eqref{eq:wave_propagation} fits into the path-conservative formalism \cite{delorenzo:2018}, which we will briefly recall in Section \ref{sec:BR_path_conservative}. The fluctuations are computed solving local Riemann problems by means of suitable Riemann solvers. Following the discussion at the beginning of Section \ref{ssec:Riemann_solvers}, the structure of the solution defined by the Riemann solver can be expressed in general by a set of $M$ jumps $\boldsymbol{\mathcal{W}}^{(l)}, l=1 \dots M$ associated with the waves and corresponding speed $s^{(l)}$. The sum of the waves must be equal to the initial jump in the vector $\mathbf{q}$ of the system variables. Moreover, for any variable of the model system governed by a conservative equation (including total mixture energy), the solver must guarantee that the initial jump in the associated ﬂux function is recovered by the sum of the waves multiplied by the corresponding speeds. The fluctuations are then computed as \cite{pelanti:2014}
\begin{equation}\label{eq:wave_propagation_waves}
	\boldsymbol{\mathcal{A}}^{\pm}\Delta\mathbf{q}_{j+1/2} = \sum_{k=1}^{M}\rpth{s^{(k)}_{j+1/2}}^{\pm}\boldsymbol{\mathcal{W}}^{(k)}_{j+1/2},
\end{equation}
where
$$\rpth{s^{(k)}_{j+1/2}}^{+} = \max\rpth{s^{(k)}_{j+1/2}, 0}, \qquad \rpth{s^{(k)}_{j+1/2}}^{-} = \max\rpth{-s^{(k)}_{j+1/2}, 0}.$$
For conservative systems, \eqref{eq:wave_propagation}-\eqref{eq:wave_propagation_waves} is equivalent to the classical Godunov method for simple approximate Riemann solvers \cite{pelanti:2018}. Indeed, conservative schemes in their usual form can be rewritten as
\begin{equation}
	\mathbf{q}_{j}^{n+1} = \mathbf{q}_{j}^{n} - \frac{\Delta t}{\Delta x}\rpth{\cpth{\hat{\mathbf{F}}_{j+1/2}\rpth{\mathbf{q}_{j}^{n}, \mathbf{q}_{j+1}^{n}} - \mathbf{F}(\mathbf{q}_{j})} + \cpth{\mathbf{F}(\mathbf{q}_{j}) -\hat{\mathbf{F}}_{j-1/2}\rpth{\mathbf{q}_{j-1}^{n}, \mathbf{q}_{j}^{n}}} }.
\end{equation}
The terms in curly brackets can then be expressed as jumps multiplied by the corresponding wave speeds via the Rankine--Hugoniot conditions, thus defining the fluctuations $\mathcal{A}^{\mp}\Delta\mathbf{q}_{j\pm 1/2}$ \cite{leveque:1997}. Notice that, in such case, one can easily verify that
\begin{equation}
	\boldsymbol{\mathcal{A}}^{-}\Delta\mathbf{q}_{j+1/2} + \boldsymbol{\mathcal{A}}^{+}\Delta\mathbf{q}_{j+1/2} = \sum_{l=1}^{M}s^{(l)}\boldsymbol{\mathcal{W}}^{(l)} = \mathbf{F}(\mathbf{q}_{j+1}) - \mathbf{F}(\mathbf{q}_{j}).
\end{equation}
In the case of HLLC, the speed of waves are defined as in \eqref{eq:s_HLLC} and the jumps associated with the waves are
\begin{equation}
	\boldsymbol{\mathcal{W}}^{(1)} = \mathbf{q}_{\star,L} - \mathbf{q}_{L}, \qquad \boldsymbol{\mathcal{W}}^{(2)} = \mathbf{q}_{\star,R} - \mathbf{q}_{\star,L}, \qquad \boldsymbol{\mathcal{W}}^{(3)} = \mathbf{q}_{R} - \mathbf{q}_{\star,R}.
\end{equation}
Following \cite{pelanti:2014}, the states $\mathbf{q}_{\star,L}$ and $\mathbf{q}_{\star,R}$ coincide with \eqref{eq:q_star_HLLC}. Hence, the non-conservative contribution $\vel \cdot \boldsymbol{\Sigma}$ in \eqref{eq:system_hyperbolic} is neglected through shocks and discretized only at contact discontinuities (see the discussion in the previous Section \ref{ssec:Riemann_solvers}).

For the sake of clarity, we conclude this section by providing an overview of the numerical schemes that will be compared through numerical experiments in Section \ref{sec:num_res_homogeneous} and \ref{sec:num_res_relax}. We consider either a Rusanov or HLLC flux for the conservative part in combination with either the \textit{BR-2023} \eqref{eq:BR_Orlando_non_cons}, the \textit{BR-2015} \eqref{eq:BR_Tumolo_non_cons} or \Crouzet\ \eqref{eq:Crouzet_non_cons} for the discretization of non-conservative terms. Additionally, we consider an HLLC-type wave-propagation scheme that naturally incorporates non-conservative terms. Finally, as extensively discussed in Section \ref{ssec:approximate_jump_relations}, the closures proposed in \cite{pelanti:2014} are adopted for the HLLC scheme.

\section{Numerical results without mechanical relaxation}
\label{sec:num_res_homogeneous}

The techniques outlined in Section \ref{sec:num_hyp} are now employed in a number of benchmarks consisting of Riemann problem tests where initially a discontinuity separates two uniform states. In this section, we will focus on the results obtained solving the homogeneous problem \eqref{eq:system_hyperbolic}, while we will discuss benchmarks with the instantaneous mechanical relaxation in Section \ref{sec:num_res_relax}. As already mentioned in Section \ref{sec:intro}, this stage is particularly relevant because recent contributions have equipped the 6-equation model with finite-rate mechanical relaxation \cite{pelanti:2022, schmidmayer:2023}. When doing so, the 6-equation model is no longer merely an auxiliary system for the numerical solution of the Kapila model, but becomes a standalone model in its own right. For the sake of compactness, we report here the parameters of all the one-dimensional test cases analyzed in this work. The computational domain is $\Omega = \rpth{0,1}$. The parameters of the EOS, the location of the initial discontinuity, and the final time are listed in Table \ref{tab:parameters}. The initial conditions for phase $1$ and phase $2$ are listed in Table \ref{tab:phase1_init} and \ref{tab:phase2_init}, respectively. We consider two computational meshes, a coarse mesh composed by $N_{\text{cells}} = 1024$ cells and a fine one composed by $N_{\text{cells}} = 65536$ cells. Discrete parameter choices for the numerical simulations are associated with the so-called Courant number
\begin{equation}\label{eq:Courant}
	C = \rpth{\left|\vel\right|_{\infty} + c_{\text{f}}}\Delta t/{\Delta x},
\end{equation}
where $\left|\vel\right|_{\infty}$ is the $l^{\infty}$ norm of the vector velocity and $c_{\text{f}}$ denotes the mixture speed of sound \eqref{eq:speed_sound_mixture}. Unless differently stated, we consider a time step such that the Courant number \eqref{eq:Courant} is equal to $0.9$.

\begin{table}[h!]
	\centering
	\footnotesize
	\begin{tabularx}{0.97\columnwidth}{lrrrrrrrr}
		\toprule
		\textbf{Test case} & $T_{f}$ & $x_{0}$ & $\gamma_{1}$ & $\pi_{1}$ & $\eta_{1}$ & $\gamma_{2}$ & $\pi_{2}$ & $\eta_{2}$ \\
		\midrule
		Sonic rarefaction & $0.15$ & $0.5$ & $1.4$ & $0$ & $0$ & $1.4$ & $0$ & $0$ \\
		\midrule
		Low-density flow & $0.15$ & $0.5$ & $1.4$ & $0$ & $0$ & $1.4$ & $0$ & $0$ \\ 
		\midrule
		Water-air shock tube & $2.4 \times 10^{-4}$ & $0.7$ & $4.4$ & $6 \times 10^{8}$ & $0$ & $1.4$ & $0$ & $0$ \\
		\midrule
		Epoxy-spinel shock & $2.9 \times 10^{-5}$ & $0.6$ & $2.43$ & $5.3 \times 10^{9}$ & $0$ & $1.62$ & $141 \times 10^{9}$ & $0$ \\
		\midrule
		Water cavitation tube & $3.2 \times 10^{-3}$ & $0.5$ & $2.35$ & $10^{9}$ & $-1.167 \times 10^{6}$ & $1.43$ & $0$ & $2.030 \times 10^{6}$ \\
		\bottomrule
	\end{tabularx}
	\caption{Final time $T_{f}$, location of the initial discontinuity $x_{0}$, and EOS parameters for the one-dimensional test cases, see main text for details. All the parameters are in SI units. Here, and in the following tables of this section, Sonic rarefaction: test case in Section \ref{ssec:sonic_rarefaction}. Low-density flow: test case in Section \ref{ssec:low_density}. Water-air shock tube: test case in Section \ref{ssec:water_air_shock_homogeneous} and in Section \ref{ssec:water_air_shock_relax}. Epoxy-spinel shock: test case in Section \ref{ssec:epoxy_spinel_homogeneous} and in Section \ref{ssec:epoxy_spinel_relax}. Water cavitation tube: test case in Section \ref{ssec:cavitation}.}
	\label{tab:parameters}
\end{table}

\begin{table}[h!]
	\centering
	\footnotesize
	\begin{tabularx}{0.8\columnwidth}{lrrrrrrrr}
		\toprule
		\textbf{Test case} & $\alpha_{1_{L}}$ & $\rho_{1_{L}}$ & $p_{1_{L}}$ & $u_{L}$ & $\alpha_{1_{R}}$ & $\rho_{1_{R}}$ & $p_{1_{R}}$ & $u_{R}$ \\
		\midrule
		Sonic rarefaction & $0.8$ & $1$ & $1$ & $0.75$ & $0.3$ & $0.125$ & $0.1$ & $0$ \\
		\midrule
		Low-density flow & $0.8$ & $1$ & $0.4$ & $-2$ & $0.5$ & $1$ & $0.4$ & $2$ \\ 
		\midrule
		Water-air shock tube & $1 - 10^{-6}$ & $10^{3}$ & $10^{9}$ & $0$ & $10^{-6}$ & $10^{3}$ & $10^{5}$ & $0$ \\
		\midrule
		Epoxy-spinel shock & $0.5954$ & $1185$ & $2 \times 10^{11}$ & $0$ & $0.5954$ & $1185$ & $10^{5}$ & $0$ \\
		\midrule
		Water cavitation tube & $0.99$ & $1150$ & $10^{5}$ & $-2$ & $0.99$ & $1150$ & $10^{5}$ & $2$ \\
		\bottomrule
	\end{tabularx}
	\caption{Initial conditions of phase $1$ variables for the one-dimensional test cases, see main text for details. All the parameters are in SI units. The subscripts $L$ and $R$ denote the values at the left and at the right, respectively, of the discontinuity $x_{0}$ whose location is reported in Table \ref{tab:parameters}.}
	\label{tab:phase1_init}
\end{table}

\begin{table}[h!]
	\centering
	\footnotesize
	\begin{tabularx}{0.8\columnwidth}{lrrrrrrrr}
		\toprule
		\textbf{Test case} & $\alpha_{2_{L}}$ & $\rho_{2_{L}}$ & $p_{2_{L}}$ & $u_{L}$ & $\alpha_{2_{R}}$ & $\rho_{2_{R}}$ & $p_{2_{R}}$ & $u_{R}$ \\
		\midrule
		Sonic rarefaction & $0.2$ & $1$ & $1$ & $0.75$ & $0.7$ & $0.125$ & $0.1$ & $0$ \\
		\midrule
		Low-density flow & $0.2$ & $1$ & $0.4$ & $-2$ & $0.5$ & $1$ & $0.4$ & $2$ \\ 
		\midrule
		Water-air shock tube & $10^{-6}$ & $1$ & $10^{9}$ & $0$ & $1 - 10^{-6}$ & $1$ & $10^{5}$ & $0$ \\
		\midrule
		Epoxy-spinel shock & $0.4046$ & $3622$ & $2 \times 10^{11}$ & $0$ & $0.4046$ & $3622$ & $10^{5}$ & $0$ \\
		\midrule
		Water cavitation tube & $0.01$ & $0.63$ & $10^{5}$ & $-2$ & $0.01$ & $0.63$ & $10^{5}$ & $2$ \\
		\bottomrule
	\end{tabularx}
	\caption{Initial conditions of phase $2$ variables for the one-dimensional test cases, see main text for details. All the parameters are in SI units. The subscripts $L$ and $R$ denote the values at the left and at the right, respectively, of the discontinuity $x_{0}$ whose location is reported in Table \ref{tab:parameters}.}
	\label{tab:phase2_init}
\end{table}

\subsection{Sonic rarefaction test case}
\label{ssec:sonic_rarefaction}

\begin{figure}[h!]
	\centering
	\begin{subfigure}{0.475\textwidth}
		\centering
		\includegraphics[width = 0.95\textwidth]{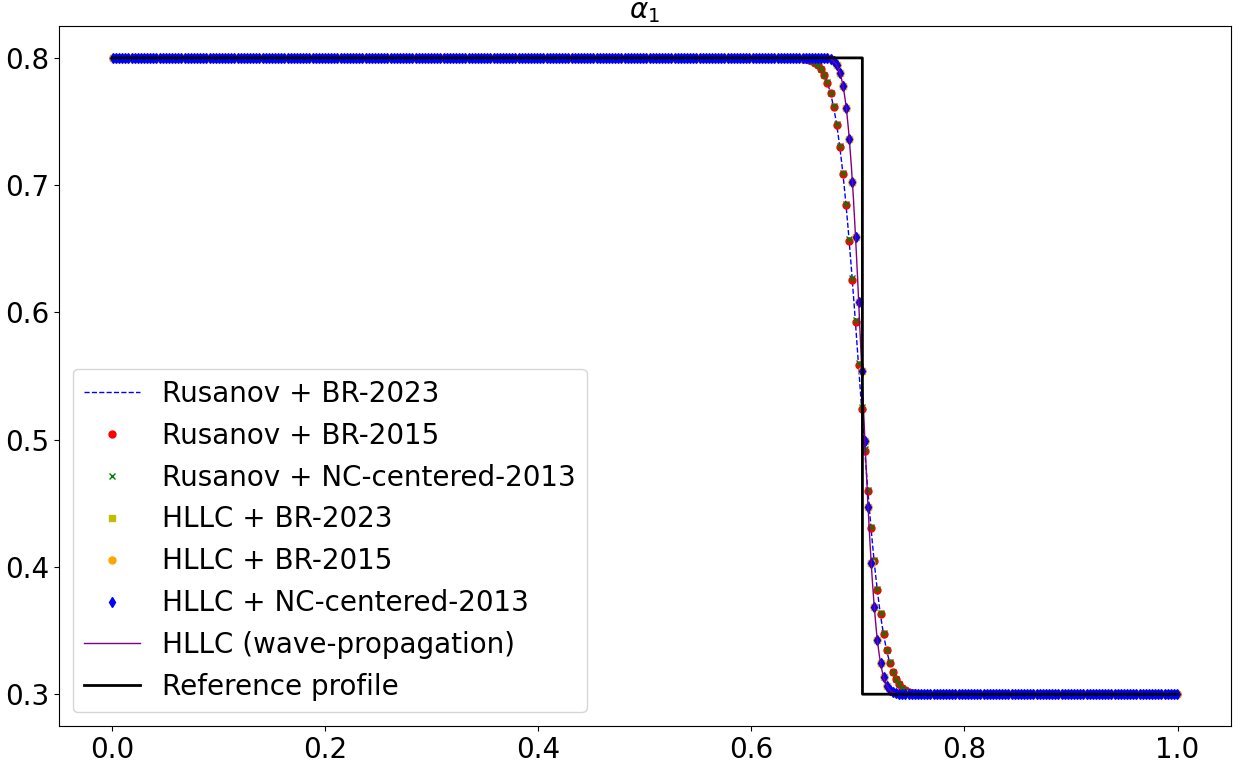}
	\end{subfigure}
	\begin{subfigure}{0.475\textwidth}
		\centering
		\includegraphics[width = 0.95\textwidth]{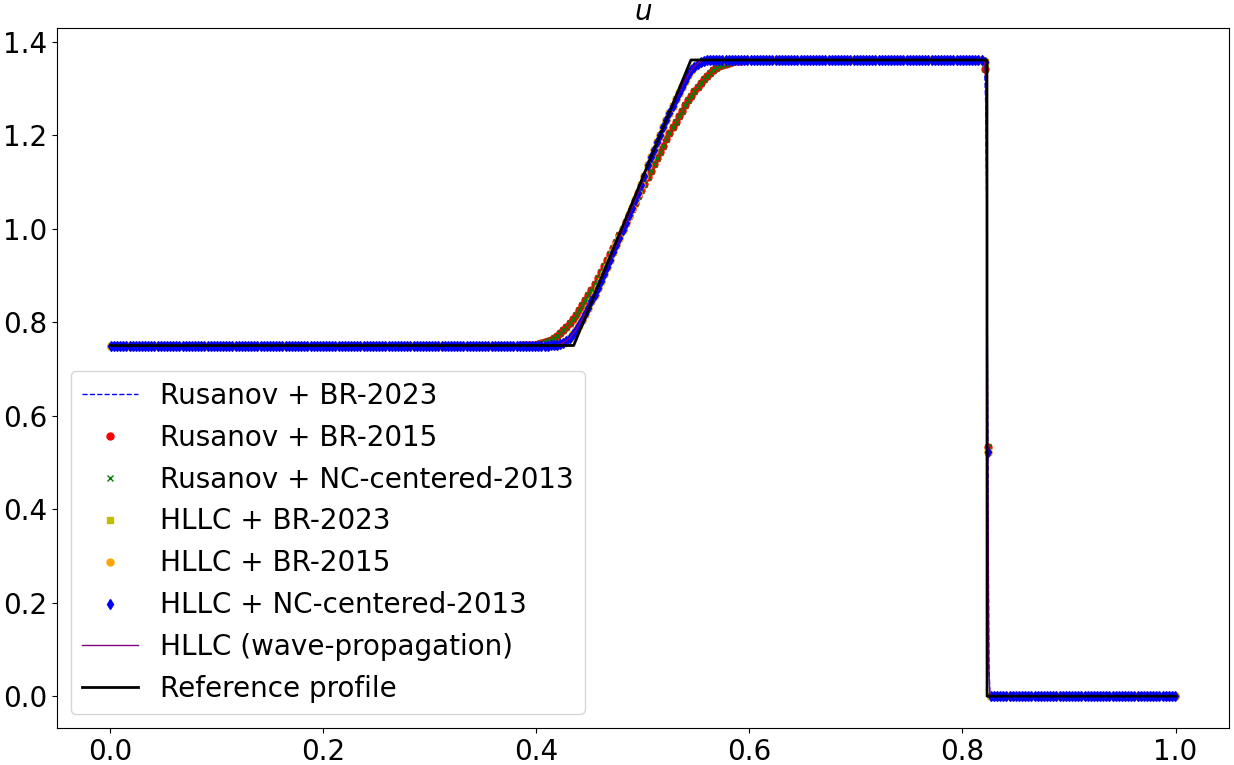}
	\end{subfigure}
	\begin{subfigure}{0.475\textwidth}
		\centering
		\includegraphics[width = 0.95\textwidth]{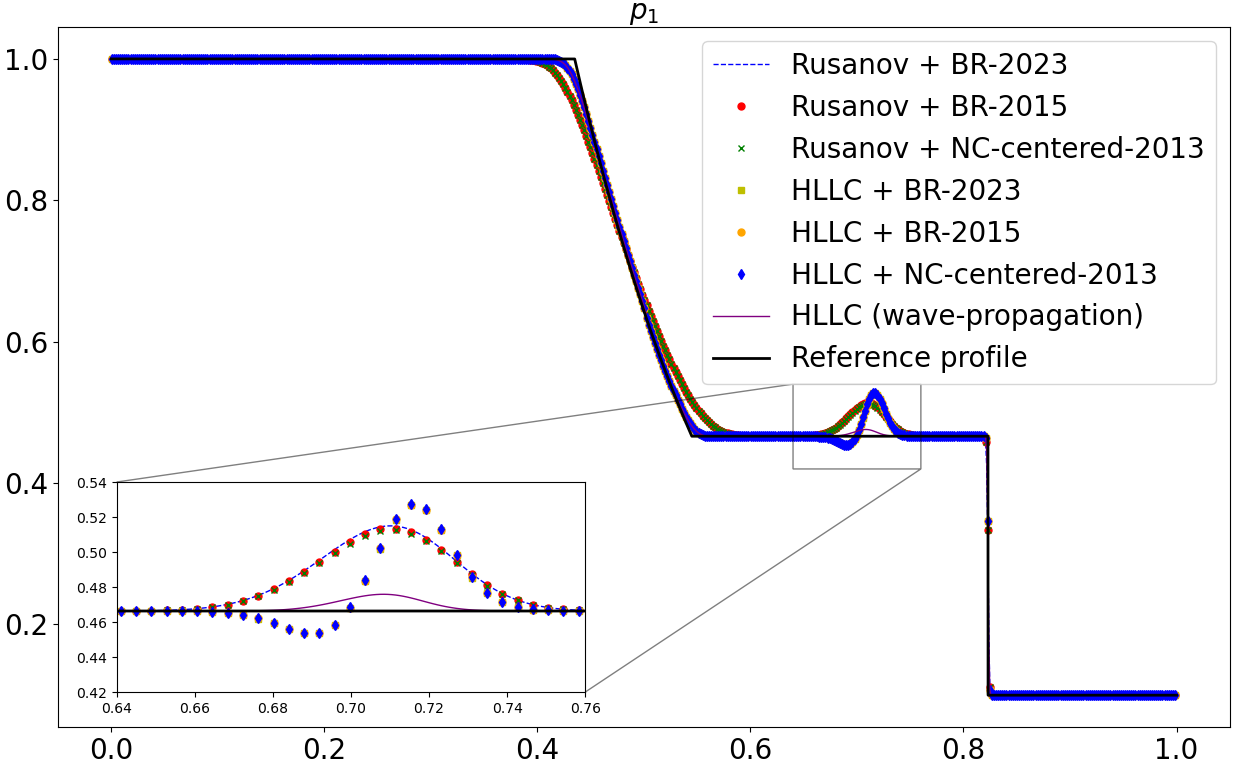}
	\end{subfigure}
	\begin{subfigure}{0.475\textwidth}
		\centering
		\includegraphics[width = 0.95\textwidth]{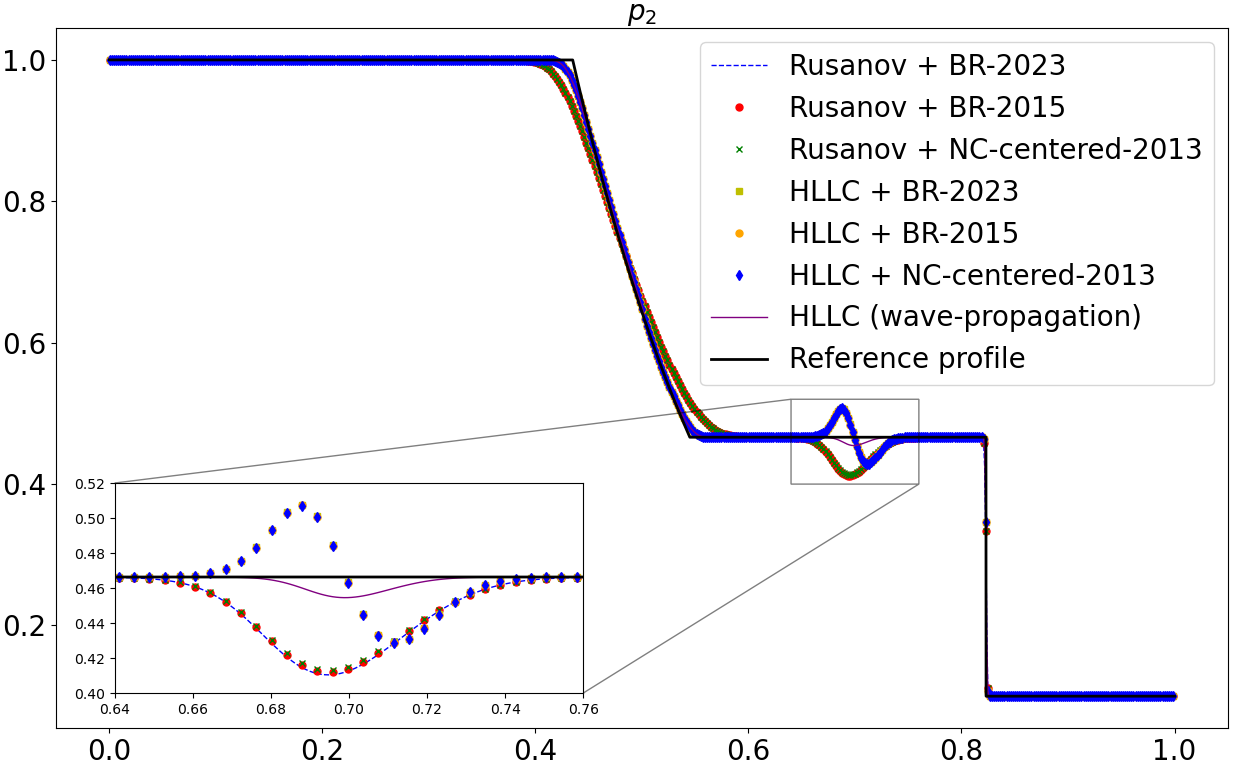}
	\end{subfigure}
	\caption{Sonic rarefaction test case with the coarse mesh, results at $t = T_{f}$. Top-left: volume fraction of phase $1$. Top-right: velocity. Bottom-left: pressure of phase $1$. Bottom-right: Pressure of phase $2$. Continuous purple lines: HLLC-type wave-propagation scheme. Dashed blue lines: Rusanov flux in combination with \textit{BR-2023} \eqref{eq:BR_Orlando_non_cons}. Red dots: Rusanov flux in combination with \textit{BR-2015} \eqref{eq:BR_Tumolo_non_cons}. Green crosses: Rusanov flux in combination with \Crouzet\ \eqref{eq:Crouzet_non_cons}. Yellow squares: HLLC flux in combination with \textit{BR-2023} \eqref{eq:BR_Orlando_non_cons}. Orange dots: HLLC flux in combination with \textit{BR-2015} \eqref{eq:BR_Tumolo_non_cons}. Blue diamonds: HLLC flux in combination with \Crouzet\ \eqref{eq:Crouzet_non_cons}. Continuous black lines: analytical solution of the Euler equations with passive transport of the volume fraction.}
	\label{fig:sonic_rarefaction_coarse}
\end{figure}

First, we consider a configuration which is inspired by the test case 3 in \cite{tokareva:2010} for the full Baer-Nunziato model. The configuration employed in \cite{tokareva:2010} considers indeed kinematic equilibrium and therefore it can be employed in the framework of a 6-equation single-velocity model. The solution for both phases consists of a right shock wave, a right travelling contact discontinuity and a left sonic rarefaction wave \cite{tokareva:2010}.
First, we employ the coarse mesh. As expected, the HLLC scheme provides a sharper solution along the contact discontinuity and the rarefaction (Figure \ref{fig:sonic_rarefaction_coarse}). The approaches for non-conservative terms based on the \textit{BR-2023} \eqref{eq:BR_Orlando_non_cons}, the \textit{BR-2015} \eqref{eq:BR_Tumolo_non_cons}, and the \Crouzet\ \eqref{eq:Crouzet_non_cons} provide similar results. The HLLC scheme for the conservative part of the system and the volume fraction in combination with the discretization of the non-conservative terms of the energy equations shows more oscillations for the phasic pressures with respect to the corresponding HLLC-type wave-propagation scheme (Figure \ref{fig:sonic_rarefaction_coarse}). For the sake of completeness, since the same EOS for the two phases is employed, the analytical solution of the Euler equations with passive transport of the volume fraction is also included and the relative errors in $l^{1}$ norms have been computed (Table \ref{tab:relative_errors_sonic_coarse}), which quantitatively confirm the previous considerations. In particular, when comparing HLLC-type methods, the wave-propagation scheme exhibits for the phasic pressures errors that are roughly half those obtained using the HLLC scheme for the conservative part of the system and the volume fraction, combined with the discretization of the non-conservative terms in the energy equations (Table \ref{tab:relative_errors_sonic_coarse}).

\begin{table}[h!]
	\centering
	\footnotesize
	\begin{tabularx}{0.75\columnwidth}{lXXXX}
		\toprule
		\multirow{2}{*}{\textbf{Scheme}} & \multicolumn{3}{c}{$l^{1}$ norm rel. error} & \\
		\cmidrule(l){2-5}
		& $\alpha_{1}$ & $u$ & $p_{1}$ & $p_{2}$ \\
		\midrule
		Rusanov + \textit{BR-2023} & \num{1.06e-2} & \num{8.53e-3} & \num{1.05e-2} & \num{1.11e-2} \\
		\midrule
		Rusanov + \textit{BR-2015} & \num{1.06e-2} & \num{8.53e-3} & \num{1.04e-2} & \num{1.10e-2} \\
		\midrule
		Rusanov + \textit{NC-centered-2013} & \num{1.06e-2} & \num{8.53e-3} & \num{1.03e-2} & \num{1.09e-2} \\
		\midrule
		HLLC + \textit{BR-2023} & \num{6.43e-3} & \num{3.33e-3} & \num{4.51e-3} & \num{4.53e-3} \\
		\midrule
		HLLC + \textit{BR-2015} & \num{6.43e-3} & \num{3.33e-3} & \num{4.52e-3} & \num{4.53e-3} \\
		\midrule
		HLLC + \textit{NC-centered-2013} & \num{6.43e-3} & \num{3.33e-3} & \num{4.53e-3} & \num{4.54e-3} \\
		\midrule
		HLLC (wave-propagation) & \num{6.43e-3} & \num{3.33e-3} & \num{2.52e-3} & \num{2.62e-3} \\
		\bottomrule
	\end{tabularx}
	\caption{Sonic rarefaction test case with the coarse mesh: relative errors in the $l^{1}$ norm for volume fraction, velocity, and phasic pressures at $t = T_{f}$.}
	\label{tab:relative_errors_sonic_coarse}
\end{table}

\begin{figure}[h!]
	\centering
	\begin{subfigure}{0.475\textwidth}
		\centering
		\includegraphics[width = 0.95\textwidth]{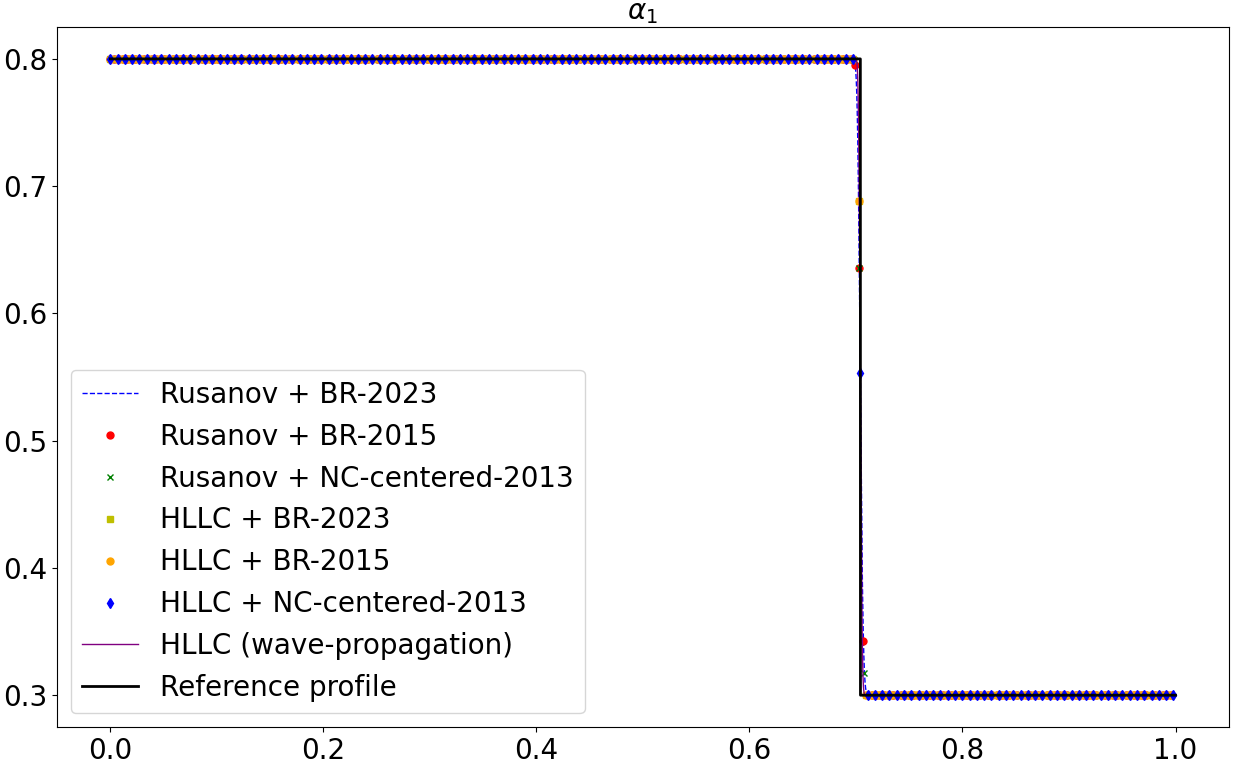}
	\end{subfigure}
	\begin{subfigure}{0.475\textwidth}
		\centering
		\includegraphics[width = 0.95\textwidth]{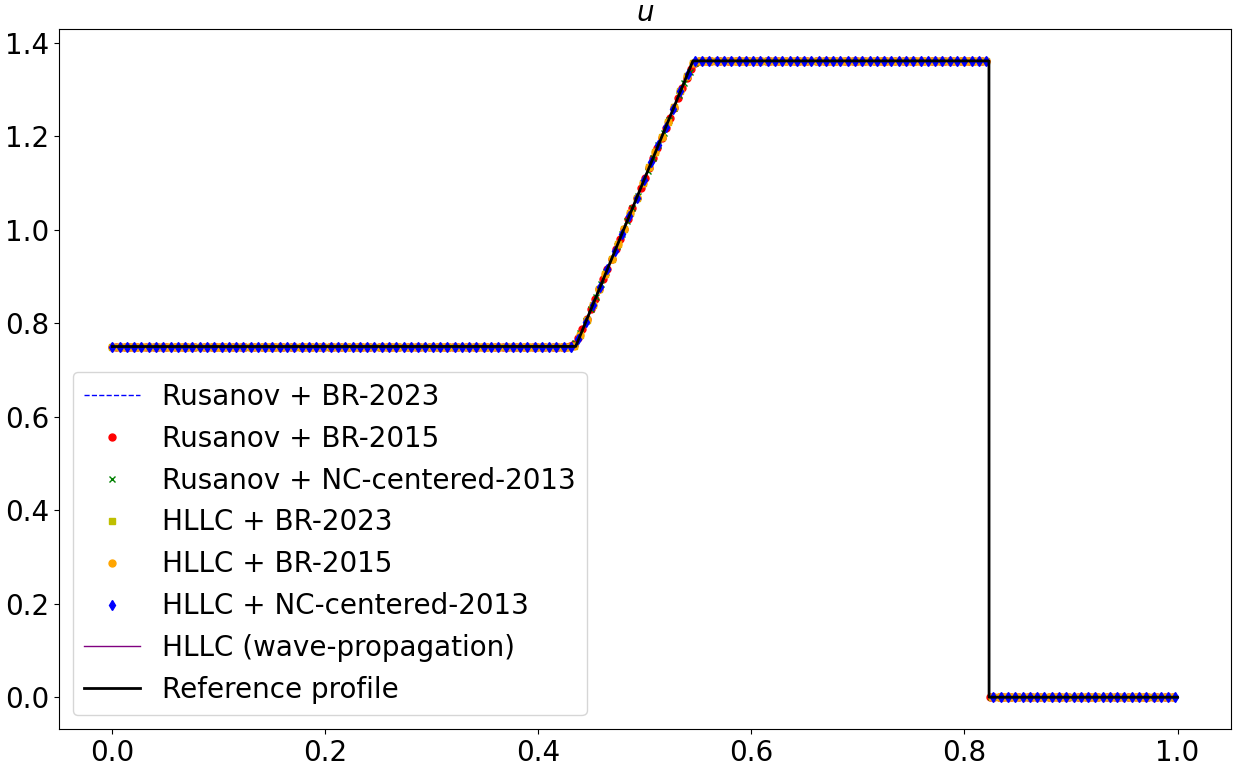}
	\end{subfigure}
	\begin{subfigure}{0.475\textwidth}
		\centering
		\includegraphics[width = 0.95\textwidth]{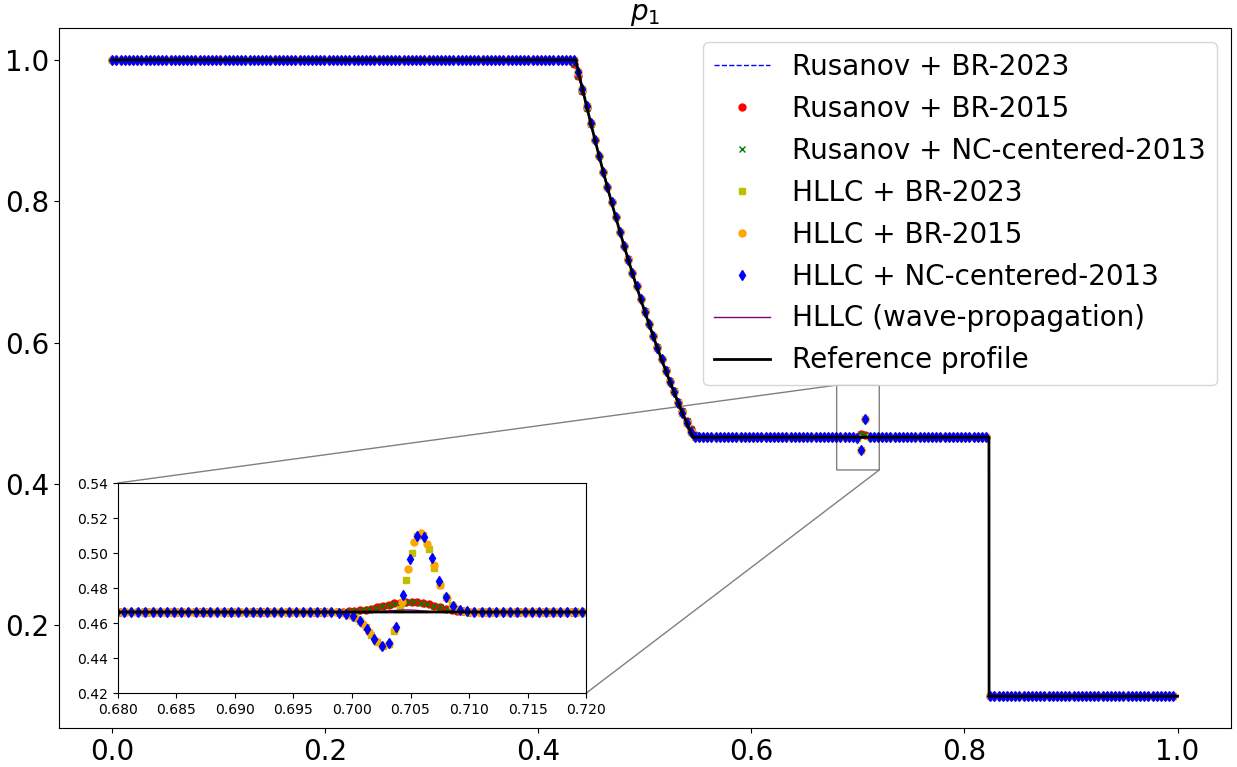}
	\end{subfigure}
	\begin{subfigure}{0.475\textwidth}
		\centering
		\includegraphics[width = 0.95\textwidth]{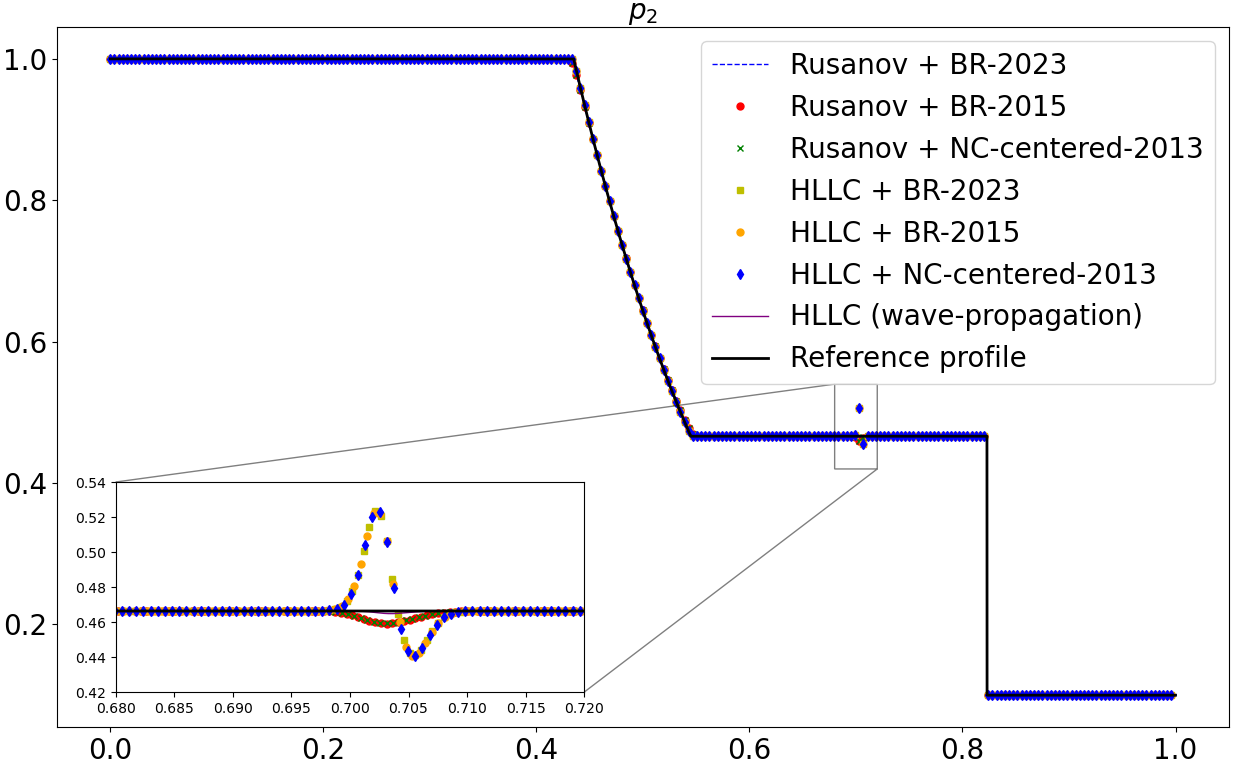}
	\end{subfigure}
	\caption{Sonic rarefaction test case with the fine mesh, results at $t = T_{f}$. Same description as Figure \ref{fig:sonic_rarefaction_coarse}.}
	\label{fig:sonic_rarefaction_fine}
\end{figure}

Next, we employ the fine mesh, so as to achieve mesh convergence. All the employed methods tend towards the expected solution of the BN system (Figure \ref{fig:sonic_rarefaction_fine}). However, we notice that, also at very high resolution, a bump is present in the phasic pressure fields in correspondence of the contact discontinuity, in particular using the HLLC scheme for the conservative part of the system and the volume fraction in combination with the three discretizations of the non-conservative terms of the energy equations. This consideration is quantitatively confirmed by the relative errors in the $l^{1}$ norm with respect to the reference profile, represented by the analytical solution of the Euler equations with passive transport of the volume fraction (Table \ref{tab:relative_errors_sonic_fine}). The wave-propagation scheme exhibits a relative error that is around 80\% lower than that of the other HLLC-type schemes. This is a purely numerical artefact as the exact solution of Riemann problems does not exhibit oscillations near contacts but only a discontinuity between the two neighbouring states. It is related to the fact that, along the contact discontinuity, only the mixture pressure $\pmix = \alpha_{1}p_{1} + \alpha_{2}p_{2}$ is preserved and not the phasic pressures (see the discussion in Section \ref{ssec:Riemann_solvers}). This behaviour does not appear solving the full Baer-Nunziato model \cite{tokareva:2010}. Hence, the 6-equation model yields different results compared to those obtained in \cite{tokareva:2010} by solving the BN model, even though the latter results still satisfy kinematic equilibrium.

\begin{table}[h!]
	\centering
	\footnotesize
	\begin{tabularx}{0.75\columnwidth}{lXXXX}
		\toprule
		\multirow{2}{*}{\textbf{Scheme}} & \multicolumn{3}{c}{$l^{1}$ norm rel. error} & \\
		\cmidrule(l){2-5}
		& $\alpha_{1}$ & $u$ & $p_{1}$ & $p_{2}$ \\
		\midrule
		Rusanov + \textit{BR-2023} & \num{1.29e-3} & \num{2.84e-4} & \num{3.21e-4} & \num{3.33e-4} \\
		\midrule
		Rusanov + \textit{BR-2015} & \num{1.29e-3} & \num{2.84e-4} & \num{3.20e-4} & \num{3.33e-4} \\
		\midrule
		Rusanov + \textit{NC-centered-2013} & \num{1.29e-3} & \num{2.84e-4} & \num{3.19e-4} & \num{3.33e-4} \\
		\midrule
		HLLC + \textit{BR-2023} & \num{8.04e-4} & \num{8.04e-5} & \num{3.13e-4} & \num{3.82e-4} \\
		\midrule
		HLLC + \textit{BR-2015} & \num{8.04e-4} & \num{8.04e-5} & \num{3.13e-4} & \num{3.82e-4} \\
		\midrule
		HLLC + \textit{NC-centered-2013} & \num{8.04e-4} & \num{8.04e-5} & \num{3.13e-4} & \num{3.82e-4} \\
		\midrule
		HLLC (wave-propagation) & \num{8.04e-4} & \num{8.04e-5} & \num{7.18e-5} & \num{7.34e-5} \\
		\bottomrule
	\end{tabularx}
	\caption{Sonic rarefaction test case with the fine mesh: relative errors in the $l^{1}$ norm for volume fraction, velocity, and phasic pressures at $t = T_{f}$.}
	\label{tab:relative_errors_sonic_fine}
\end{table}

\subsection{Low-density flow test case}
\label{ssec:low_density}

\begin{figure}[h!]
	\centering
	\begin{subfigure}{0.475\textwidth}
		\centering
		\includegraphics[width = 0.95\textwidth]{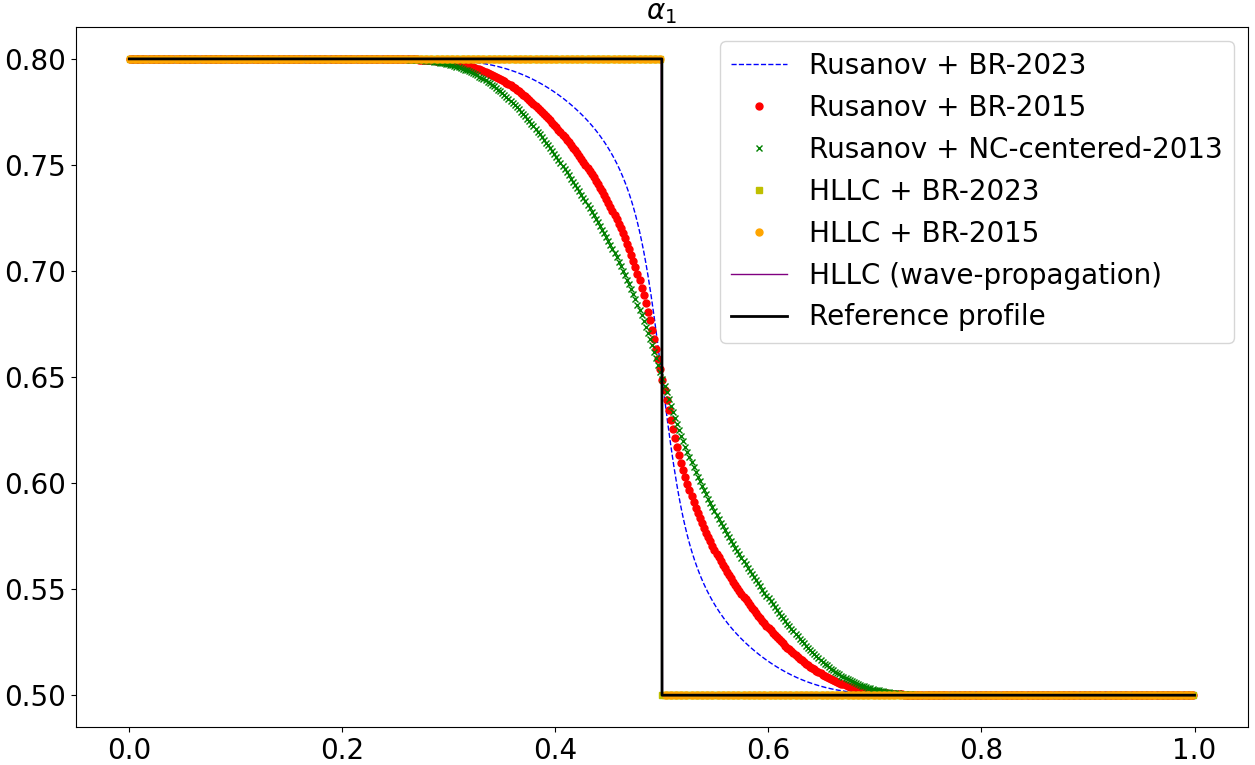}
	\end{subfigure}
	\begin{subfigure}{0.475\textwidth}
		\centering
		\includegraphics[width = 0.95\textwidth]{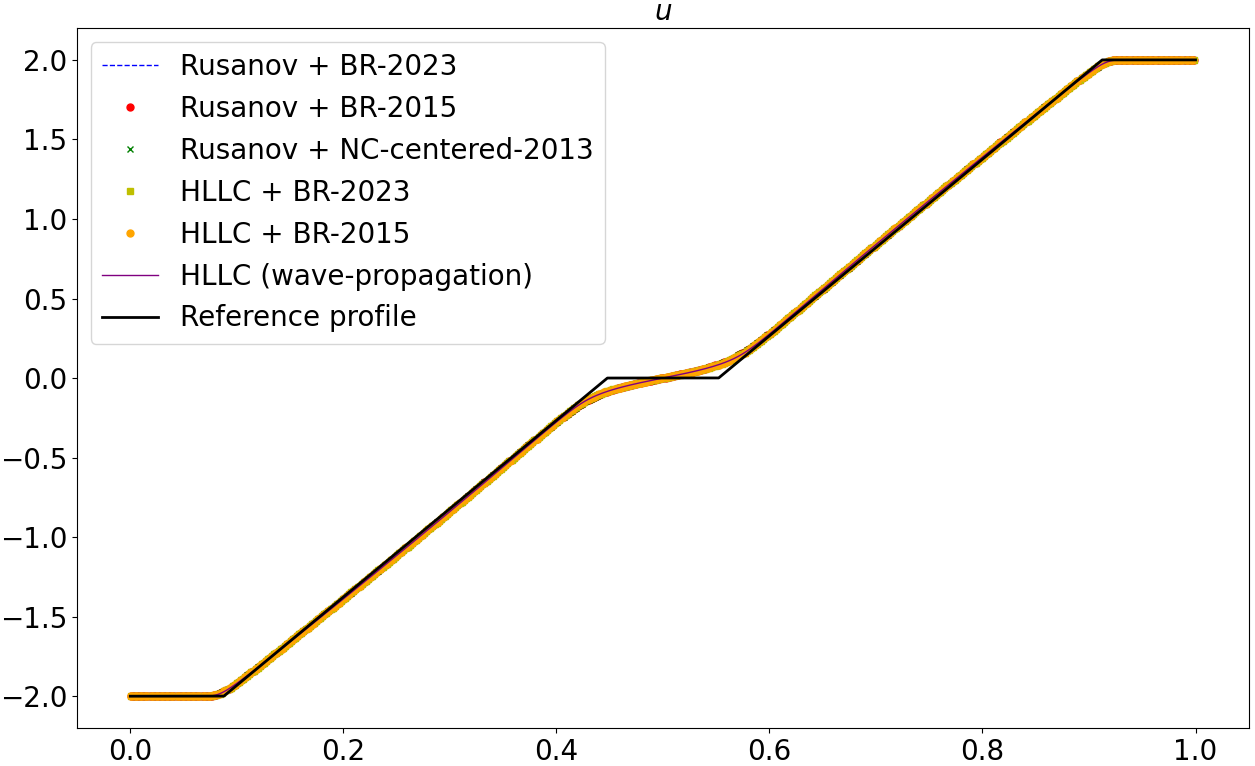}
	\end{subfigure}
	\begin{subfigure}{0.475\textwidth}
		\centering
		\includegraphics[width = 0.95\textwidth]{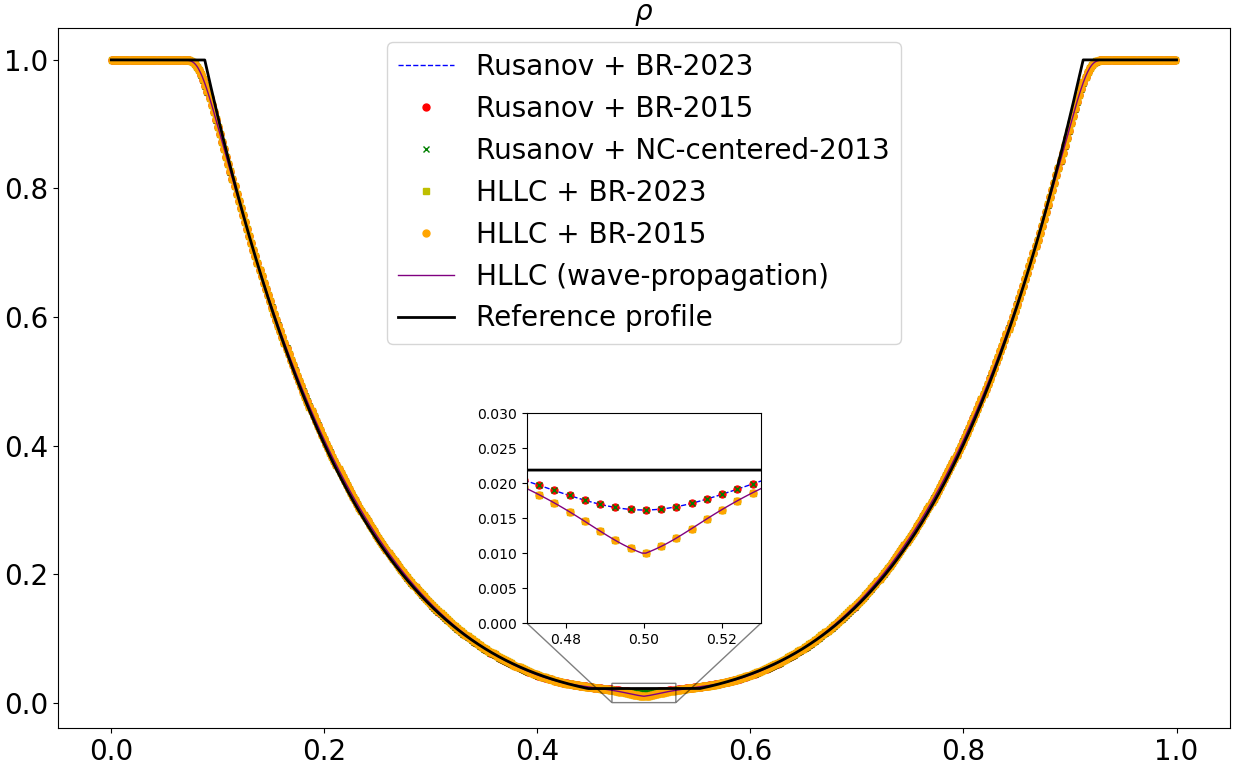}
	\end{subfigure}
	\begin{subfigure}{0.475\textwidth}
		\centering
		\includegraphics[width = 0.95\textwidth]{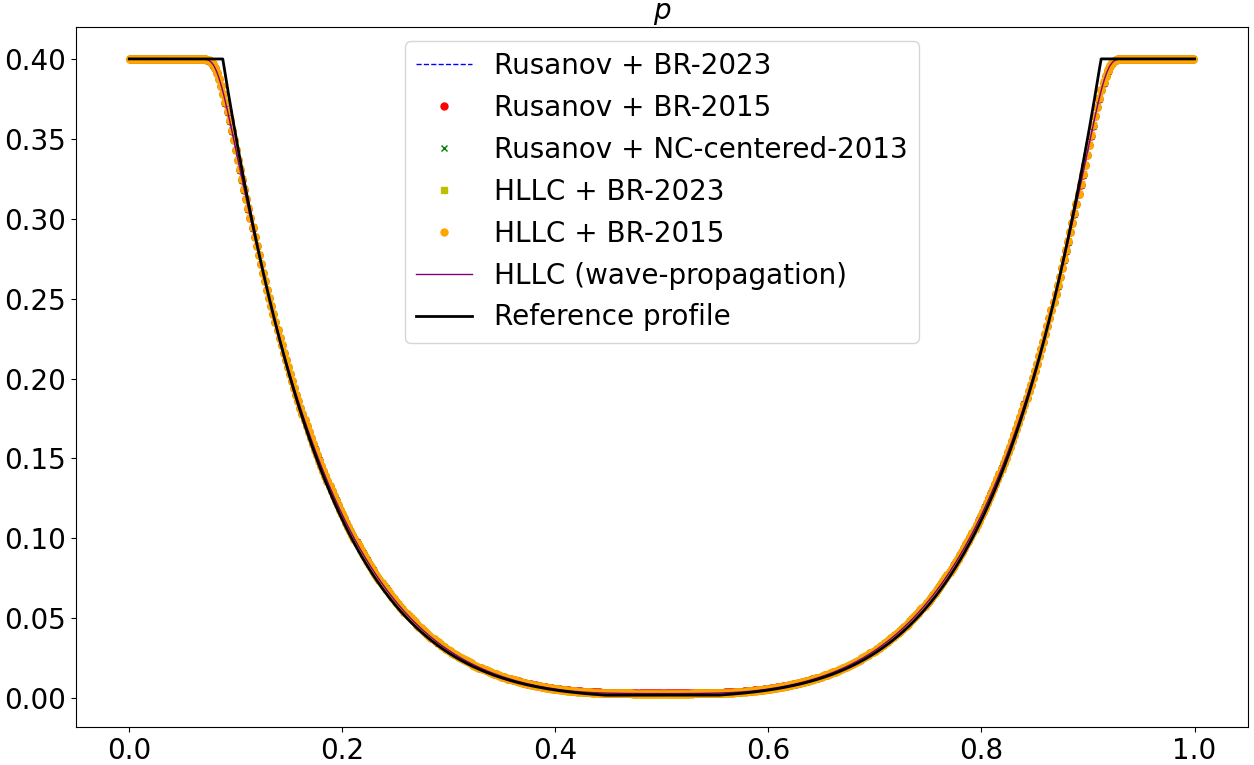}
	\end{subfigure}
	\caption{Low-density flow test case with the coarse mesh, results at $t = T_{f}$. Top-left: volume fraction of phase $1$. Top-right: velocity. Bottom-left: mixture density. Bottom-right: mixture pressure. Continuous purple lines: HLLC-type wave-propagation scheme. Dashed blue lines: Rusanov flux in combination with \textit{BR-2023} \eqref{eq:BR_Orlando_non_cons}. Red dots: Rusanov flux in combination with \textit{BR-2015} \eqref{eq:BR_Tumolo_non_cons}. Green crosses: Rusanov flux in combination with \Crouzet\ \eqref{eq:Crouzet_non_cons}. Yellow squares: HLLC flux in combination with \textit{BR-2023} \eqref{eq:BR_Orlando_non_cons}. Orange dots: HLLC flux in combination with \textit{BR-2015} \eqref{eq:BR_Tumolo_non_cons}. Continuous black lines: analytical solution of the Euler equations with passive transport of the volume fraction.}
	\label{fig:low_density_coarse}
\end{figure}

Next, we consider another configuration inspired by a test case reported in \cite{tokareva:2010}. More specifically, we employ test case 4, which is a two-phase extension of the so-called 123-problem \cite{re:2022, toro:2009}. Both phases consist of two symmetric rarefaction waves and a stationary contact wave. The region between the rarefaction waves is close to vacuum. Hence, this is a severe test problem with a low-density flow. First, we employ the coarse mesh. Similarly to what was done in Section \ref{ssec:sonic_rarefaction}, the analytical solution of the Euler equations with passive transport of the volume fraction is also included in Figure \ref{fig:low_density_coarse} and is considered as reference profile. The HLLC flux in combination with the \Crouzet\ approach for the non-conservative terms of the energy equations does not preserve the thermodynamic admissibility of the flow variables and corrupted pressure and density values are generated. In light of this result, it is not reported in Figure \ref{fig:low_density_coarse} and we will not further consider \Crouzet\ \eqref{eq:Crouzet_non_cons} in the following test cases. A positive solution is instead obtained employing the BR-type and the wave-propagation approach. A spike in the density arises across the contact discontinuity (Figure \ref{fig:low_density_coarse}). Moreover, when coupled with a Rusanov flux, the BR-type approach yields a sharper description of the contact discontinuity with respect to the \Crouzet\ \eqref{eq:Crouzet_non_cons} approach, in particular using the \textit{BR-2023} \eqref{eq:BR_Orlando_non_cons} (Figure \ref{fig:low_density_coarse}). This is quantitatively confirmed by the relative errors in the $l^{1}$ norm with respect to the reference profile (Table \ref{tab:relative_errors_low_density_coarse}).

\begin{table}[h!]
	\centering
	\footnotesize
	\begin{tabularx}{0.75\columnwidth}{lXXXX}
		\toprule
		\multirow{2}{*}{\textbf{Scheme}} & \multicolumn{3}{c}{$l^{1}$ norm rel. error} & \\
		\cmidrule(l){2-5}
		& $\alpha_{1}$ & $u$ & $\rho$ & $p$ \\
		\midrule
		Rusanov + \textit{BR-2023} & \num{1.84e-2} & \num{1.49e-2} & \num{9.736e-3} & \num{1.51e-2} \\
		\midrule
		Rusanov + \textit{BR-2015} & \num{2.69e-2} & \num{1.49e-2} & \num{9.73e-3} & \num{1.51e-2} \\
		\midrule
		Rusanov + \textit{NC-centered-2013} & \num{3.36e-2} & \num{1.49e-2} & \num{9.73e-3} & \num{1.51e-2} \\
		\midrule
		HLLC + \textit{BR-2023} & \num{2.26e-4} & \num{1.45e-2} & \num{1.00e-2} & \num{1.43e-2} \\
		\midrule
		HLLC + \textit{BR-2015} & \num{2.26e-4} & \num{1.45e-2} & \num{1.00e-2} & \num{1.43e-2} \\
		\midrule
		HLLC (wave-propagation) & \num{2.26e-4} & \num{1.45e-2} & \num{1.00e-2} & \num{1.43e-2} \\
		\bottomrule
	\end{tabularx}
	\caption{Low-density flow test case with the coarse mesh: relative errors in the $l^{1}$ norm for volume fraction, velocity, mixture density, and mixture pressure at $t = T_{f}$.}
	\label{tab:relative_errors_low_density_coarse}
\end{table}

\begin{figure}[h!]
	\centering
	\begin{subfigure}{0.475\textwidth}
		\centering
		\includegraphics[width = 0.95\textwidth]{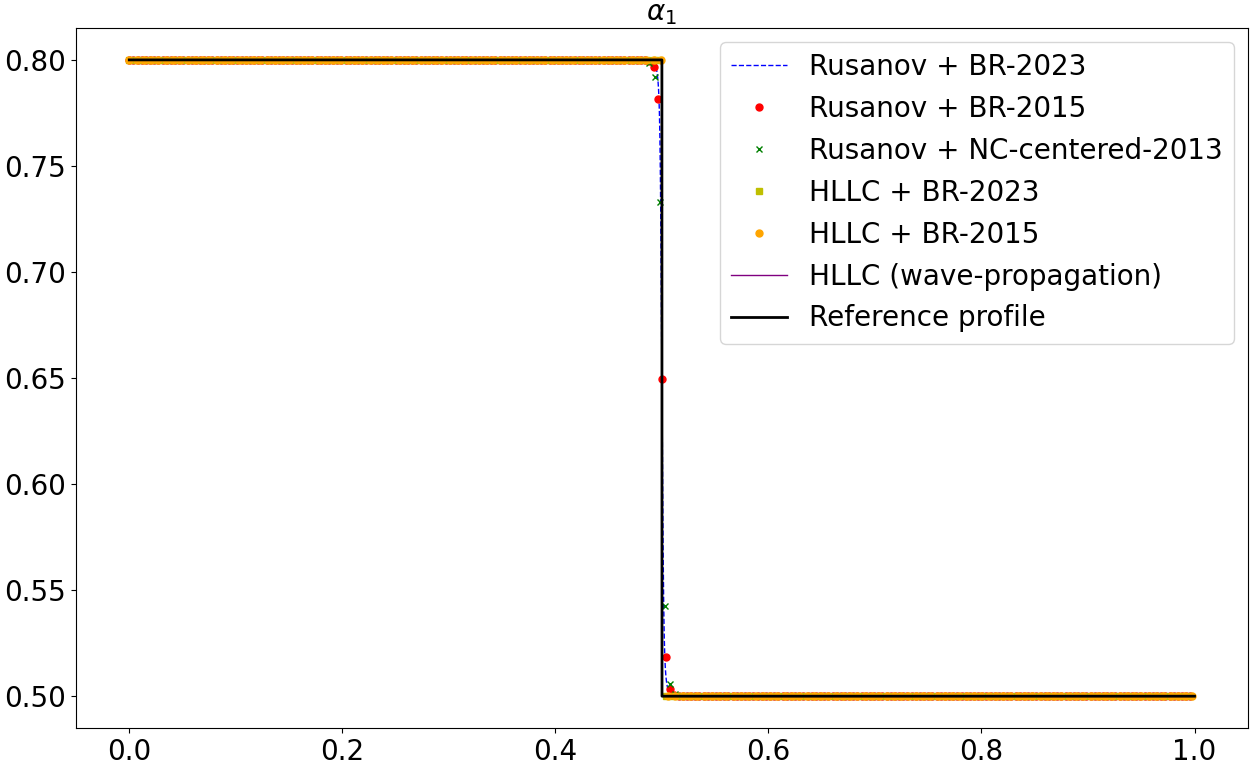}
	\end{subfigure}
	\begin{subfigure}{0.475\textwidth}
		\centering
		\includegraphics[width = 0.95\textwidth]{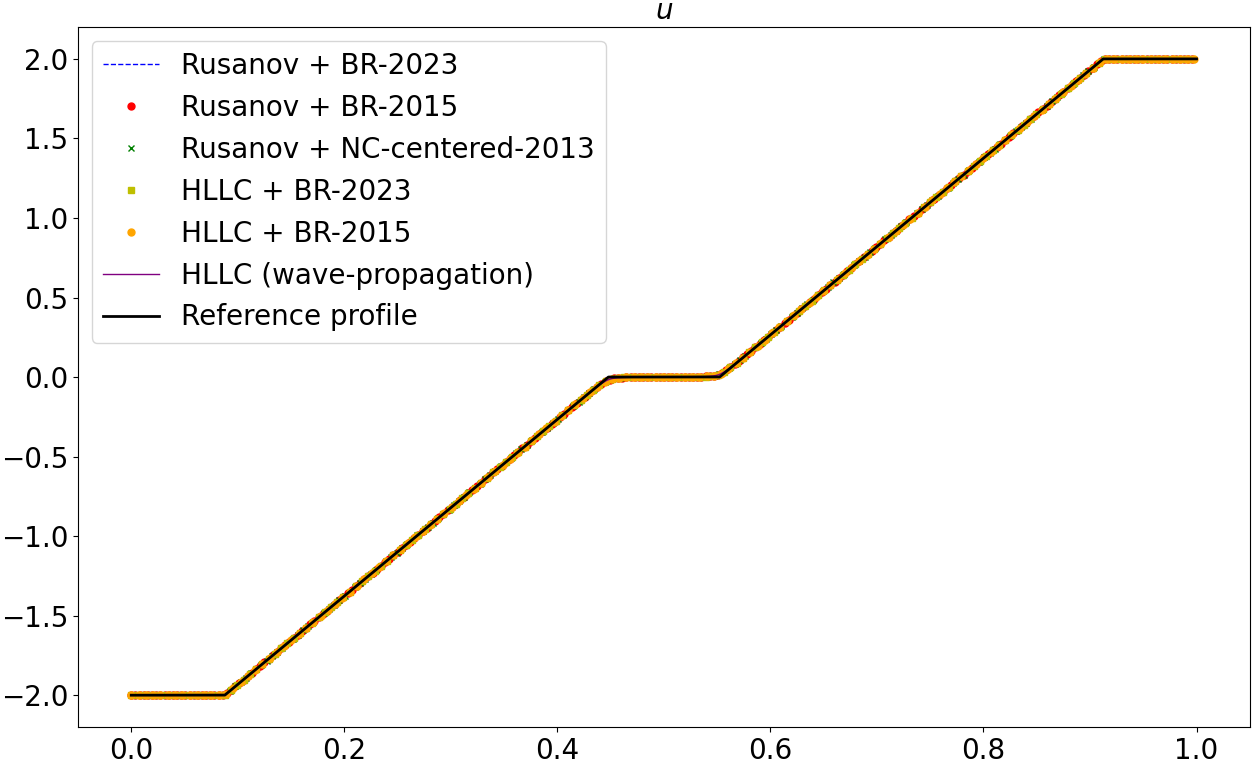}
	\end{subfigure}
	\begin{subfigure}{0.475\textwidth}
		\centering
		\includegraphics[width = 0.95\textwidth]{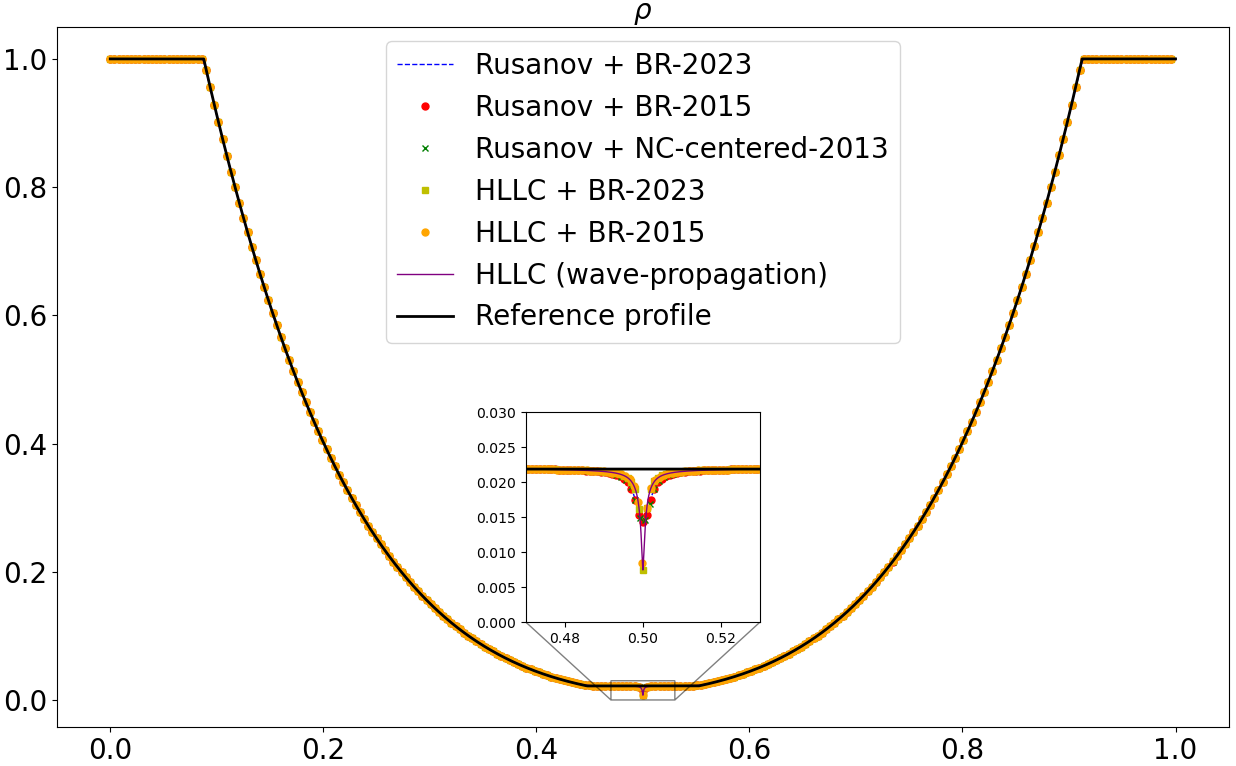}
	\end{subfigure}
	\begin{subfigure}{0.475\textwidth}
		\centering
		\includegraphics[width = 0.95\textwidth]{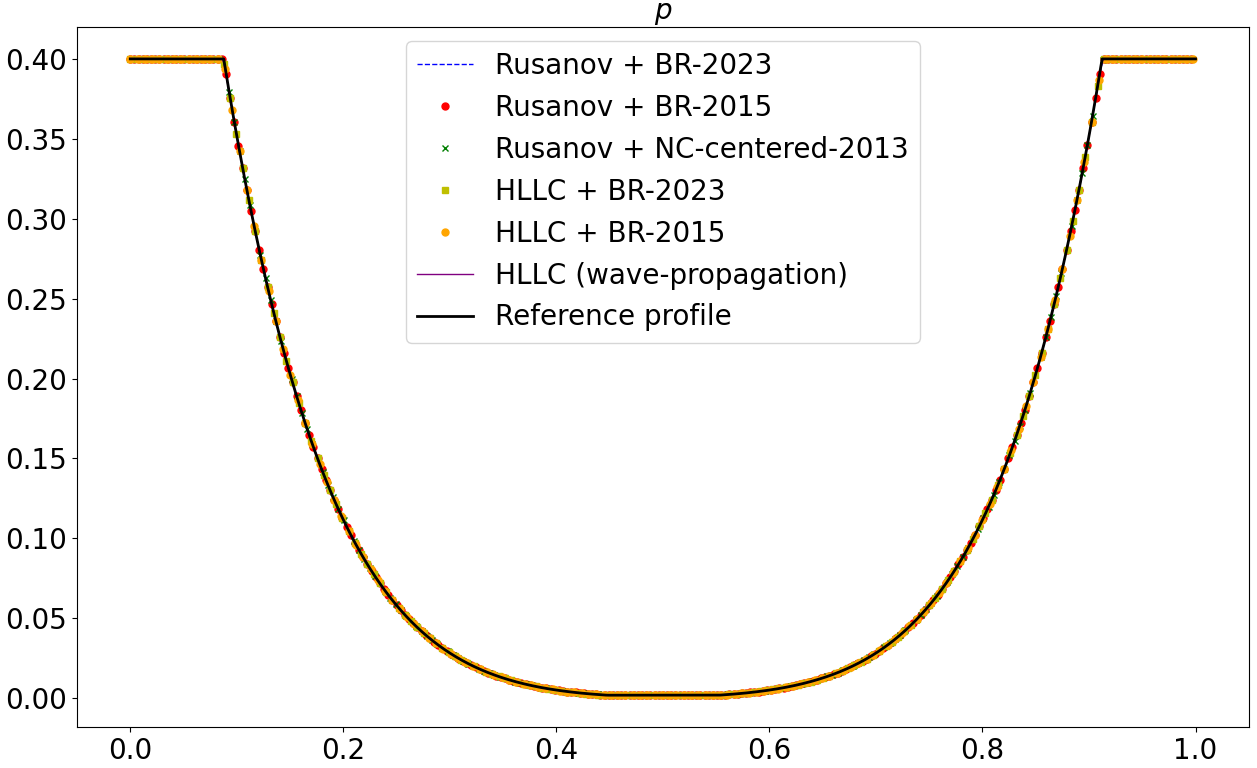}
	\end{subfigure}
	\caption{Low-density flow test case with the fine mesh, results at $t = T_{f}$. Same description as Figure \ref{fig:low_density_coarse}.}
	\label{fig:low_density_fine}
\end{figure}

Next, we employ the fine mesh. Analogous considerations to those reported in Section \ref{ssec:sonic_rarefaction} are valid. However, one can notice that the Rusanov flux tends to smear the contact discontinuity out even at very high resolution (Figure \ref{fig:low_density_fine}), as further confirmed by the relative errors in the $l^{1}$ norm with respect to the reference solution (Table \ref{tab:relative_errors_low_density_fine}). Moreover, the spike in the density across the contact discontinuity increases on the fine mesh, especially for HLLC-type schemes. This is a well known issue for this test case \cite{einfeldt:1991}.

\begin{table}[h!]
	\centering
	\footnotesize
	\begin{tabularx}{0.75\columnwidth}{lXXXX}
		\toprule
		\multirow{2}{*}{\textbf{Scheme}} & \multicolumn{3}{c}{$l^{1}$ norm rel. error} & \\
		\cmidrule(l){2-5}
		& $\alpha_{1}$ & $u$ & $\rho$ & $p$ \\
		\midrule
		Rusanov + \textit{BR-2023} & \num{6.60e-4} & \num{7.50e-4} & \num{4.82e-4} & \num{4.22e-4} \\
		\midrule
		Rusanov + \textit{BR-2015} & \num{8.95e-4} & \num{7.50e-4} & \num{4.82e-4} & \num{4.22e-4} \\
		\midrule
		Rusanov + \textit{NC-centered-2013} & \num{1.11e-3} & \num{7.50e-4} & \num{4.82e-4} & \num{4.22e-4} \\
		\midrule
		HLLC + \textit{BR-2023} & \num{4.62e-6} & \num{6.59e-4} & \num{4.50e-4} & \num{4.11e-4} \\
		\midrule
		HLLC + \textit{BR-2015} & \num{4.62e-6} & \num{6.59e-4} & \num{4.50e-4} & \num{4.11e-4} \\
		\midrule
		HLLC (wave-propagation) & \num{4.62e-6} & \num{6.59e-4} & \num{4.50e-4} & \num{4.11e-4} \\
		\bottomrule
	\end{tabularx}
	\caption{Low-density flow test case with the fine mesh: relative errors in the $l^{1}$ norm for volume fraction, velocity, mixture density, and mixture pressure at $t = T_{f}$.}
	\label{tab:relative_errors_low_density_fine}
\end{table}

\subsection{Water-air shock tube}
\label{ssec:water_air_shock_homogeneous}

\begin{figure}[h!]
	\centering
	\begin{subfigure}{0.475\textwidth}
		\centering
		\includegraphics[width = 0.95\textwidth]{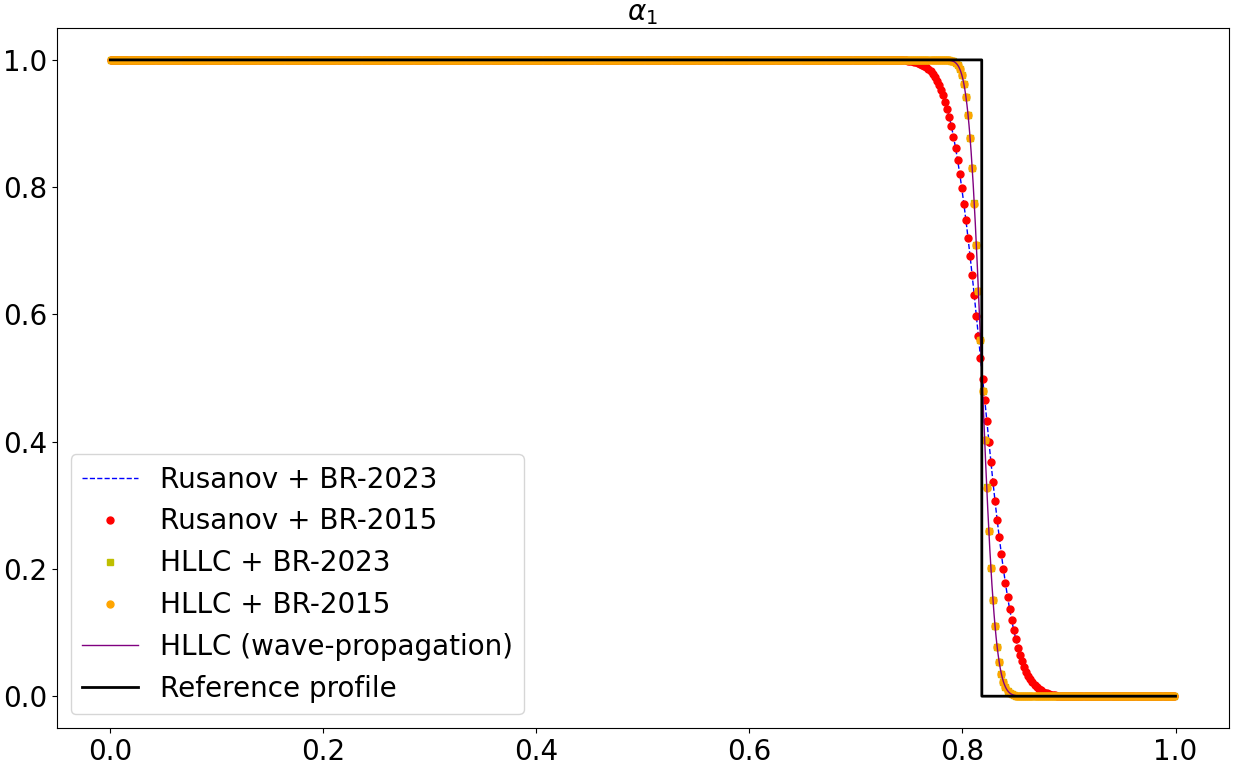}
	\end{subfigure}
	\begin{subfigure}{0.475\textwidth}
		\centering
		\includegraphics[width = 0.95\textwidth]{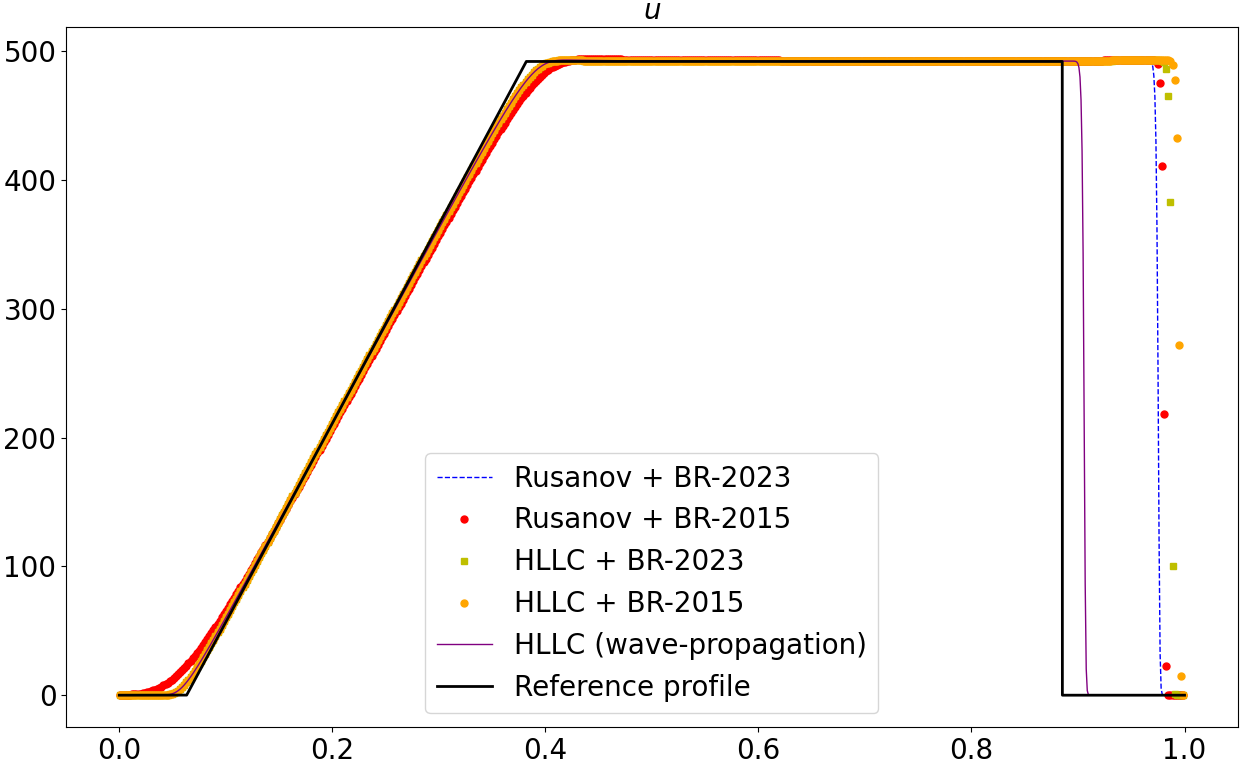}
	\end{subfigure}
	\begin{subfigure}{0.475\textwidth}
		\centering
		\includegraphics[width = 0.95\textwidth]{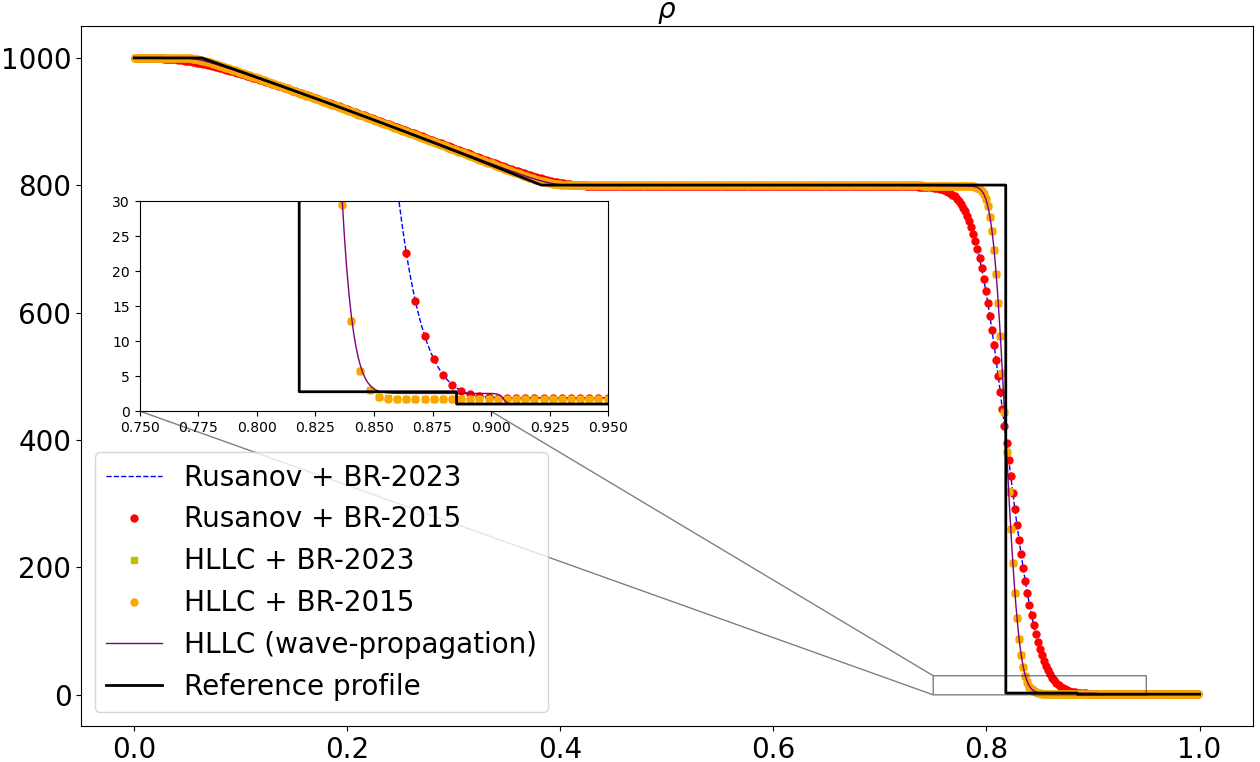}
	\end{subfigure}
	\begin{subfigure}{0.475\textwidth}
		\centering
		\includegraphics[width = 0.95\textwidth]{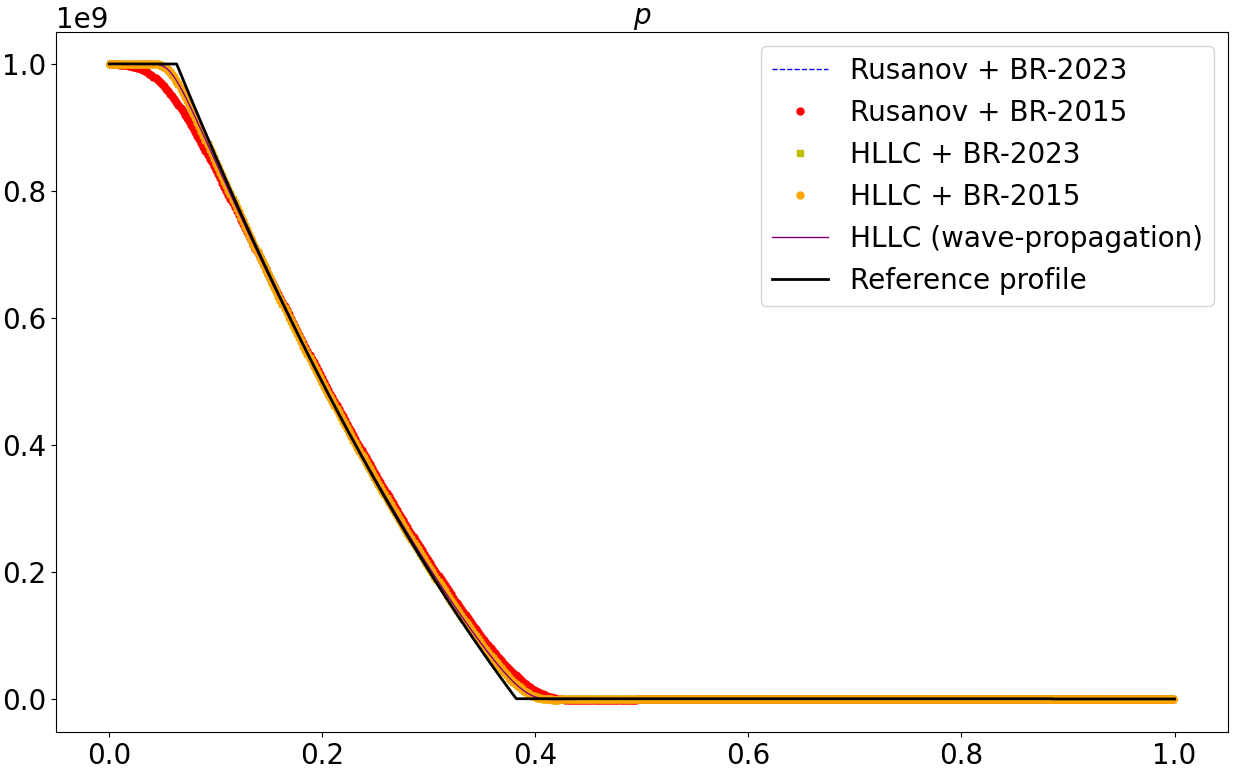}
	\end{subfigure}
	\caption{Water-air shock tube test case with the coarse mesh, results at $t = T_{f}$. Top-left: volume fraction of phase $1$. Top-right: velocity. Bottom-left: Mixture density. Bottom-right: Mixture pressure. Continuous purple lines: HLLC-type wave-propagation scheme. Dashed blue lines: Rusanov flux in combination with \textit{BR-2023} \eqref{eq:BR_Orlando_non_cons}. Red dots: Rusanov flux in combination with \textit{BR-2015} \eqref{eq:BR_Tumolo_non_cons}. Yellow squares: HLLC flux in combination with \textit{BR-2023} \eqref{eq:BR_Orlando_non_cons}. Orange dots: HLLC flux in combination with \textit{BR-2015} \eqref{eq:BR_Tumolo_non_cons}. Continuous black lines: analytical solution of the Riemann problem for the Euler--Euler equations with passive transport of the volume fraction. The results for the Rusanov schemes are established at $C = 0.29$.}
	\label{fig:water_air_shock_homogeneous_coarse}
\end{figure}

Next, we consider a water-air shock tube, analogous to those presented, e.g., in \cite{petitpas:2007, saurel:2009, schmidmayer:2023}. The shock tube is composed by two chambers each containing a nearly pure fluid (Tables \ref{tab:phase1_init}-\ref{tab:phase2_init}). First, we employ the coarse mesh. For what concerns the Rusanov flux, a physically meaningful solution is established only for Courant numbers $C < 0.3$. For larger values of the Courant number, the thermodynamic admissibility of phasic pressures is not preserved in general. Moreover, significant different results are obtained between the numerical methods, which are particularly evident from the velocity field (Figure \ref{fig:water_air_shock_homogeneous_coarse}). Since this test case is characterized by quasi-pure phases on each side of the material interface, we compare our results, for the sake of completeness, with the analytical solution of the Riemann problem for the Euler--Euler equations corresponding to a sharp-interface solution. (Figure \ref{fig:water_air_shock_homogeneous_coarse}). The analytical solution of the Riemann problem for the single-phase Euler equations governed by the SG-EOS \eqref{eq:SG_EOS} can be found, e.g., in \cite{ketcheson:2020}. One can easily notice that the numerical solution established with the HLLC-type wave-propagation scheme is closer to the analytical solution of the Euler--Euler equations compared to the other numerical schemes (Figure \ref{fig:water_air_shock_homogeneous_coarse}). The relative error in the $l^{1}$ norm with respect to the analytical solution of the Euler--Euler equations quantitatively confirms this observation (Table \ref{tab:relative_errors_water_air_shock_homogeneous_coarse}). However, we stress the fact that this solution should not be interpreted as the analytical solution for the 6-equation model, since the jump conditions are not uniquely defined. It is used here as a reference profile only because we are considering quasi-pure phases on each side of the material interface.

The results here obtained may be related to the fact that the fluctuations \eqref{eq:wave_propagation_waves} of the wave-propagation method are computed neglecting the non-conservative terms of the phasic total energies equations for the external 1-wave and 3-wave of the Riemann solution and the corresponding jump conditions reduce to the approximate jump relations of the 5-equation model presented in \cite{petitpas:2007} (see the discussion in Section \ref{ssec:Riemann_solvers}). The sensitivity of the shock location with respect to the numerical scheme of this model raises the question of the practical applicability --- of the model or the considered numerical schemes --- when no mechanical relaxation occurs, and by extension, when the mechanical relaxation occurs at finite rate only. We will show in Section \ref{sec:num_res_relax} that the differences between the different numerical methods are significantly reduced when the instantaneous mechanical relaxation is considered, thus showing a limited sensitivity of the 5-equation model's shock profiles to the numerical scheme for this configuration.

\begin{table}[h!]
	\centering
	\footnotesize
	\begin{tabularx}{0.75\columnwidth}{lXXXX}
		\toprule
		\multirow{2}{*}{\textbf{Scheme}} & \multicolumn{3}{c}{$l^{1}$ norm rel. error} & \\
		\cmidrule(l){2-5}
		& $\alpha_{1}$ & $u$ & $\rho$ & $p$ \\
		\midrule
		Rusanov + \textit{BR-2023} & \num{2.27e-2} & \num{1.48e-1} & \num{2.39e-2} & \num{3.33e-2} \\
		\midrule
		Rusanov + \textit{BR-2015} & \num{2.27e-2} & \num{1.56e-1} & \num{2.39e-2} & \num{3.33e-2} \\
		\midrule
		HLLC + \textit{BR-2023} & \num{9.60e-3} & \num{1.60e-1} & \num{1.02e-2} & \num{1.29e-2} \\
		\midrule
		HLLC + \textit{BR-2015} & \num{9.60e-3} & \num{1.70e-1} & \num{1.02e-2} & \num{1.29e-2} \\
		\midrule
		HLLC (wave-propagation) & \num{9.60e-3} & \num{3.57e-2} & \num{1.01e-2} & \num{1.25e-2} \\
		\bottomrule
	\end{tabularx}
	\caption{Water-air shock tube test case with the coarse mesh: relative errors in the $l^{1}$ norm for volume fraction, velocity, mixture density, and mixture pressure at $t = T_{f}$. The errors are computed with respect to the analytical solution of the Euler--Euler equations.}
	\label{tab:relative_errors_water_air_shock_homogeneous_coarse}
\end{table}

\begin{figure}[h!]
	\centering
	\begin{subfigure}{0.475\textwidth}
		\centering
		\includegraphics[width = 0.95\textwidth]{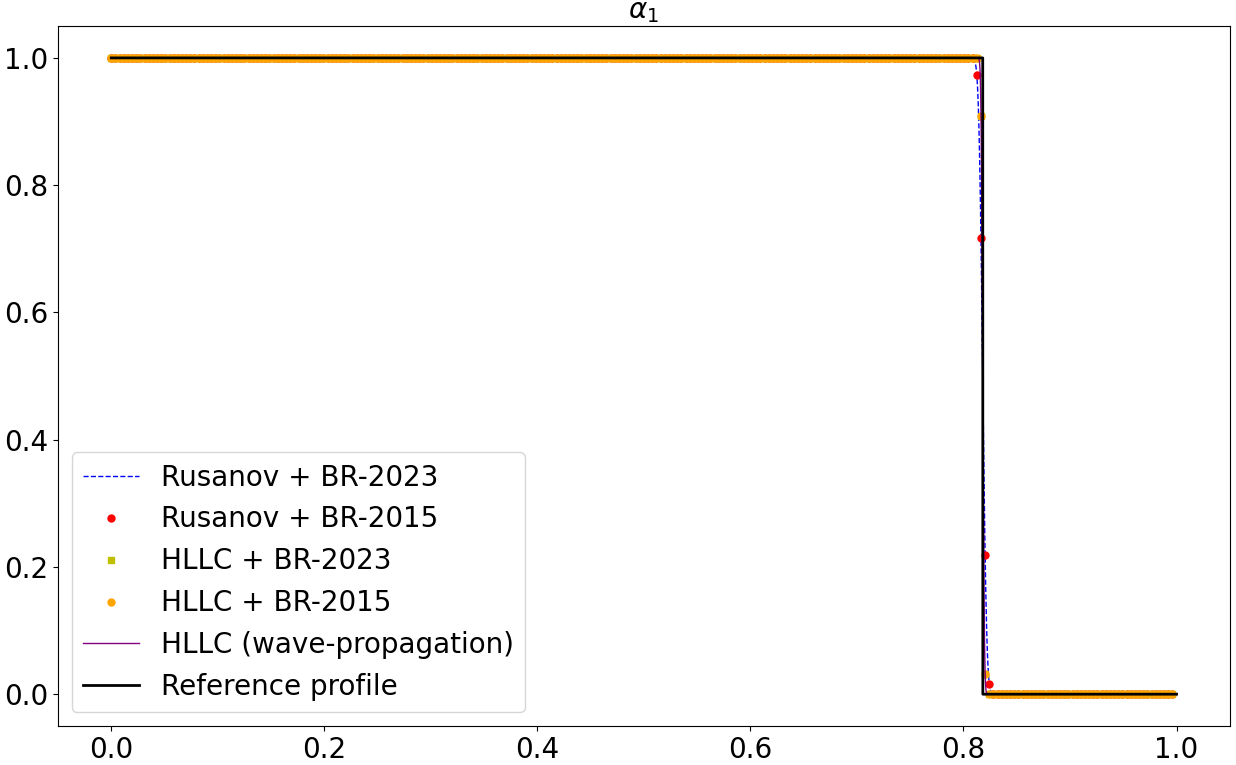}
	\end{subfigure}
	\begin{subfigure}{0.475\textwidth}
		\centering
		\includegraphics[width = 0.95\textwidth]{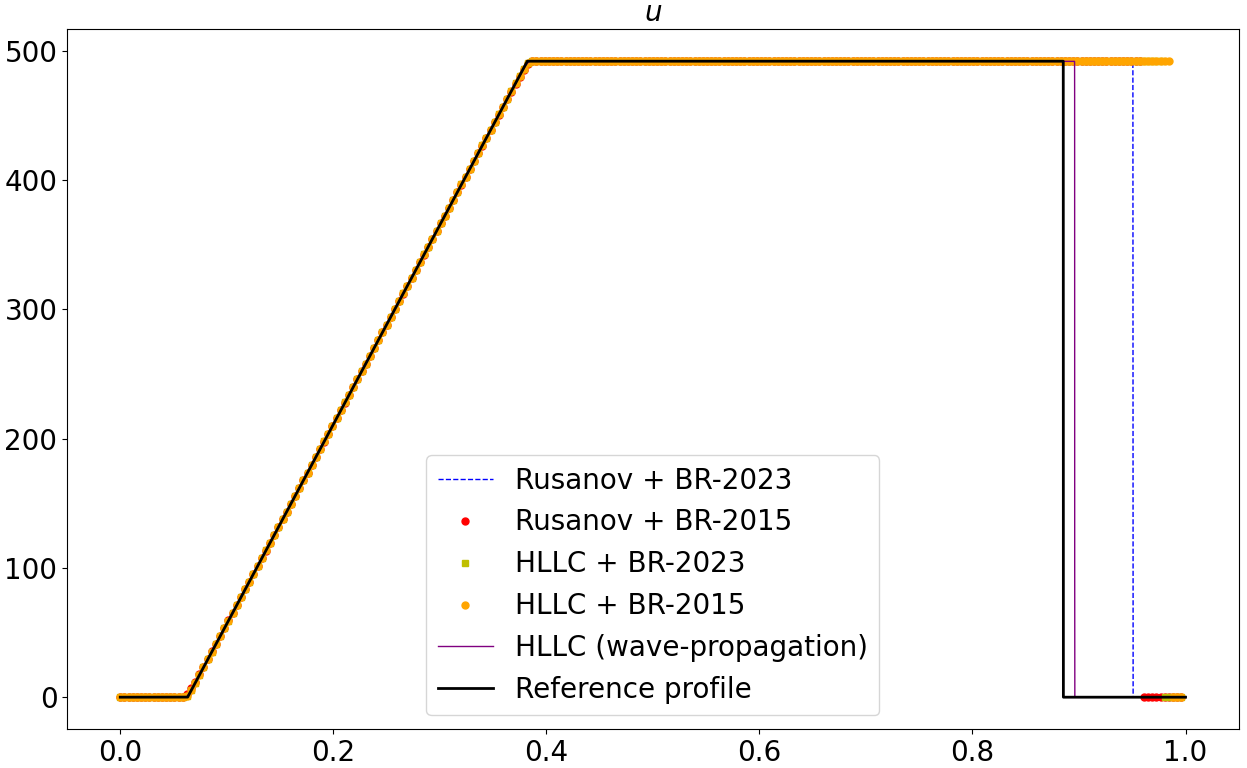}
	\end{subfigure}
	\begin{subfigure}{0.475\textwidth}
		\centering
		\includegraphics[width = 0.95\textwidth]{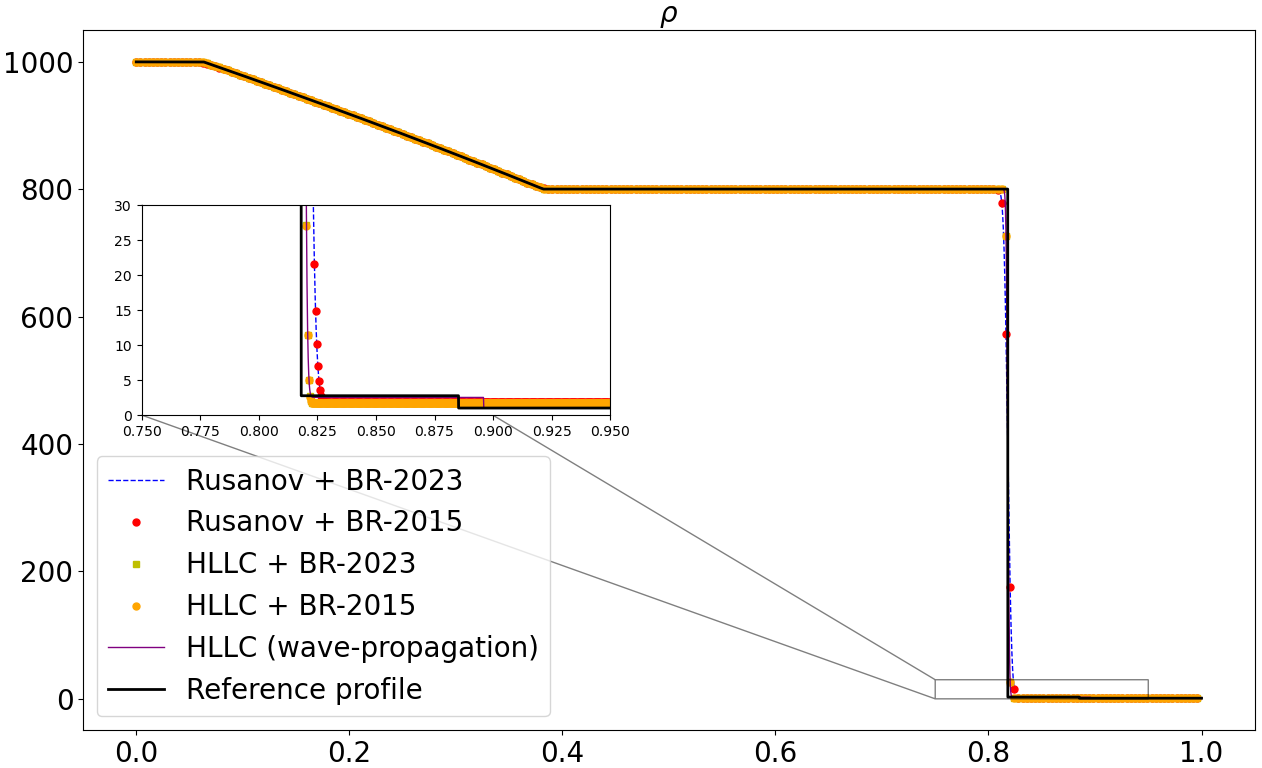}
	\end{subfigure}
	\begin{subfigure}{0.475\textwidth}
		\centering
		\includegraphics[width = 0.95\textwidth]{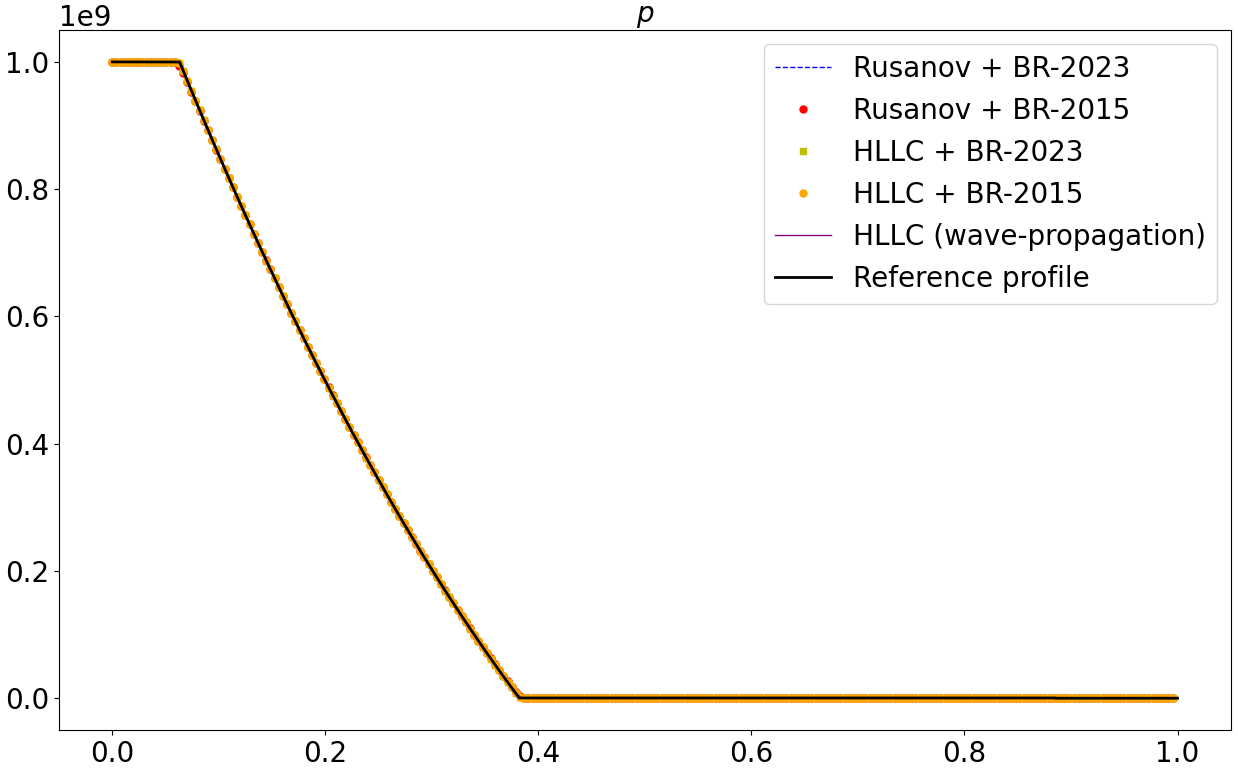}
	\end{subfigure}
	\caption{Water-air shock tube test case with the fine mesh, results at $t = T_{f}$. Same description as Figure \ref{fig:water_air_shock_homogeneous_coarse}.}
	\label{fig:water_air_shock_homogeneous_fine}
\end{figure}

Next, we employ the fine mesh so as to achieve mesh convergence. Analogous considerations to those reported for the coarse mesh are valid (Figure \ref{fig:water_air_shock_homogeneous_fine} and Table \ref{tab:relative_errors_water_air_shock_homogeneous_fine}). We also still notice visible differences in the results established with the \textit{BR-2023} \eqref{eq:BR_Orlando_non_cons} and the \textit{BR-2015} \eqref{eq:BR_Tumolo_non_cons} for the treatment of the non-conservative terms. This behaviour is likely dependent on the absence of a well-defined set of Rankine--Hugoniot conditions for this model and therefore different numerical schemes can converge to different solutions.

\begin{table}[h!]
	\centering
	\footnotesize
	\begin{tabularx}{0.75\columnwidth}{lXXXX}
		\toprule
		\multirow{2}{*}{\textbf{Scheme}} & \multicolumn{3}{c}{$l^{1}$ norm rel. error} & \\
		\cmidrule(l){2-5}
		& $\alpha_{1}$ & $u$ & $\rho$ & $p$ \\
		\midrule
		Rusanov + \textit{BR-2023} & \num{2.82e-3} & \num{9.92e-2} & \num{2.85e-3} & \num{1.43e-3} \\
		\midrule
		Rusanov + \textit{BR-2015} & \num{2.82e-3} & \num{1.11e-1} & \num{2.86e-3} & \num{1.49e-3} \\
		\midrule
		HLLC + \textit{BR-2023} & \num{1.19e-3} & \num{1.42e-1} & \num{1.35e-3} & \num{1.02e-3} \\
		\midrule
		HLLC + \textit{BR-2015} & \num{1.20e-3} & \num{1.51e-1} & \num{1.36e-3} & \num{1.08e-3} \\
		\midrule
		HLLC (wave-propagation) & \num{1.19e-3} & \num{1.63e-2} & \num{1.18e-3} & \num{3.90e-4} \\
		\bottomrule
	\end{tabularx}
	\caption{Water-air shock tube test case with the fine mesh: relative errors in the $l^{1}$ norm for volume fraction, velocity, mixture density, and mixture pressure at $t = T_{f}$. The errors are computed with respect to the analytical solution of the Euler--Euler equations.}
	\label{tab:relative_errors_water_air_shock_homogeneous_fine}
\end{table}

\begin{figure}[h!]
	\centering
	\begin{subfigure}{0.475\textwidth}
		\centering
		\includegraphics[width = 0.95\textwidth]{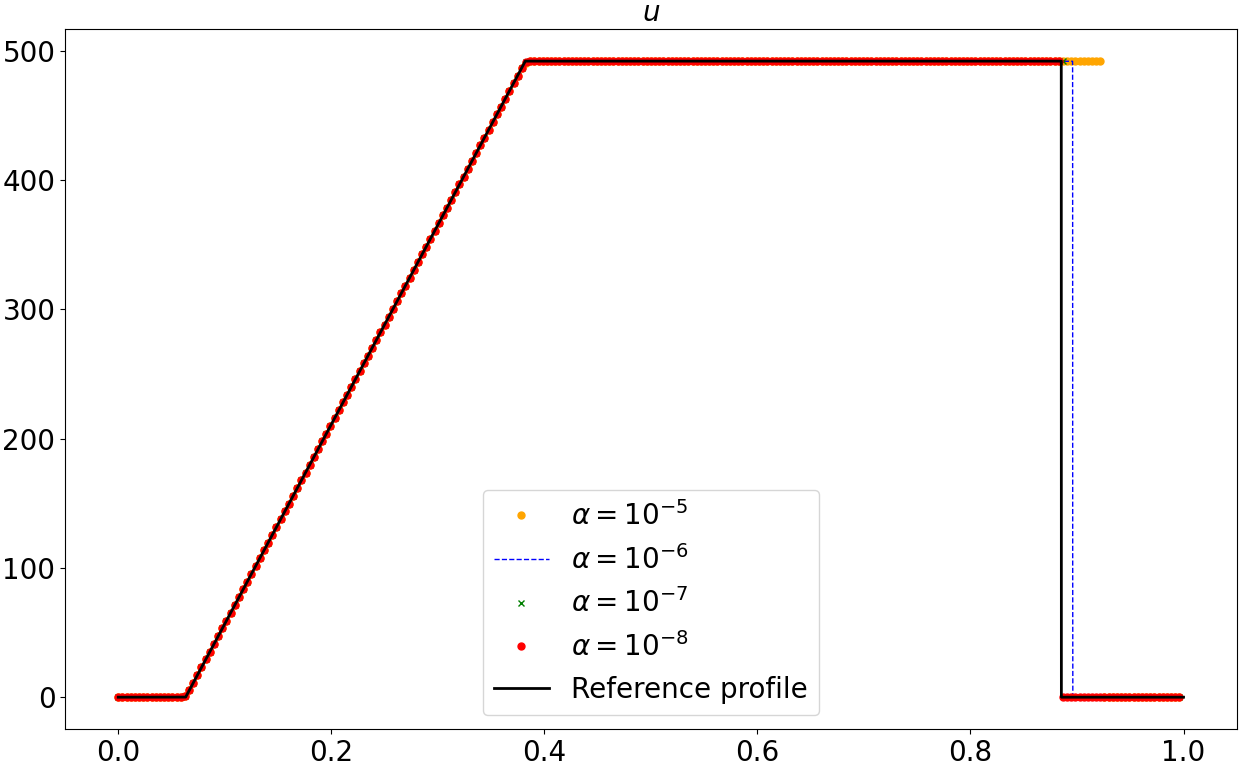}
	\end{subfigure}
	\begin{subfigure}{0.475\textwidth}
		\centering
		\includegraphics[width = 0.95\textwidth]{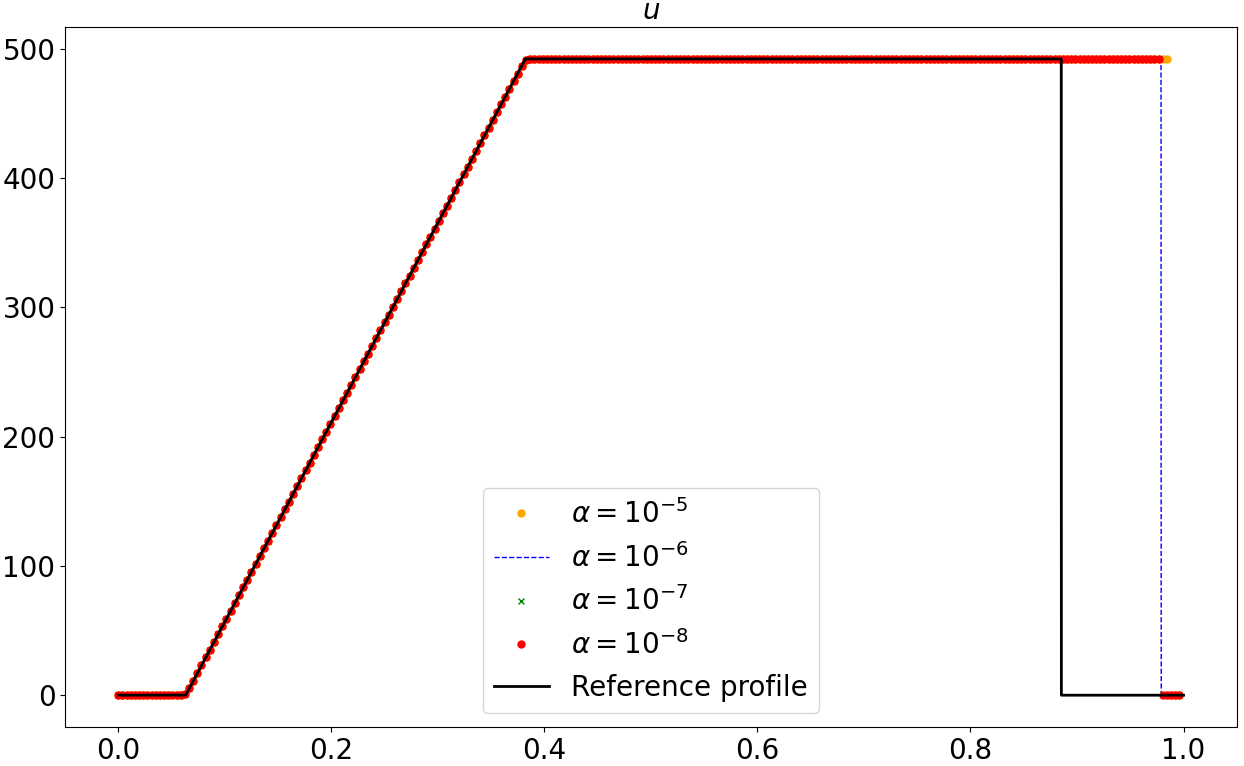}
	\end{subfigure}
	\caption{Water-air shock tube test case with the fine mesh as function of the ``residual'' volume fraction $\alpha_{1,R}$, results at $t = T_{f}$. Left: HLLC-type wave-propagation scheme. Right: HLLC flux in combination with \textit{BR-2023} \eqref{eq:BR_Orlando_non_cons}. Orange dots: $\alpha_{1,R}= 10^{-5}$. Dashed blue lines: $\alpha_{1,R}= 10^{-6}$. Green crosses: $\alpha_{1,R}= 10^{-7}$. Red dots: $\alpha_{1,R}= 10^{-8}$. Continuous black lines: analytical solution of the Riemann problem for the Euler--Euler equations with passive transport of the volume fraction.}
	\label{fig:water_air_shock_homogeneous_fine_vol_fraction}
\end{figure}

In order to further analyze the behaviour of the numerical schemes when considering nearly pure fluids, we perform a sensitivity analysis with respect to the ``residual '' volume fraction. More specifically, we consider $\alpha_{1,R} = \cpth{10^{-5}, 10^{-6}, 10^{-7}, 10^{-8}}$. We present the results obtained using the HLLC-type wave-propagation scheme and the HLLC scheme combined with the BR-type approach \textit{BR-2023} \eqref{eq:BR_Orlando_non_cons}; analogous conclusions hold for the other schemes and are therefore omitted for brevity. It can be readily observed that the wave-propagation scheme ``converges'' to the solution of the Euler--Euler equations as $\alpha_{1,R} \to 0$, whereas the results obtained by employing a discretization for the conservative part of the system and a separate, simple discretization for the non-conservative terms exhibit only a weak dependence on the ``residual'' volume fraction and do not approach the analytical solution of the Euler--Euler equations (Figure \ref{fig:water_air_shock_homogeneous_fine_vol_fraction} and Table \ref{tab:relative_errors_water_air_shock_homogeneous_fine_vol_fraction}). These results are likely a further manifestation of the absence of a uniquely defined set of jump conditions, which the wave-propagation scheme appears to be able to overcome likely because, as already mentioned, the non-conservative terms in the phasic total energy equations are neglected for the external waves. Hence, each phase follows its own Hugoniot curve, and thus the behavior of a pure phase is recovered naturally. A simple treatment of the non-conservative terms, as those considered in this work, does not seem to be sufficient when dealing with quasi-pure phases for the 6-equation model; therefore, more sophisticated non-conservative schemes should be considered or developed. This, however, goes beyond the scope of the present work.

\begin{table}[h!]
	\centering
	\footnotesize
	\begin{tabularx}{0.65\columnwidth}{lll}
		\toprule
		\multirow{2}{*}{\textbf{``Residual'' volume fraction}} & \multicolumn{1}{c}{$l^{1}$ norm rel. error} & \\
		\cmidrule(l){2-3}
		& HLLC (wave-propagation) & HLLC + \textit{BR-2023} \\
		\midrule
		$10^{-5}$ & \num{5.70e-2} & \num{1.54e-1} \\
		\midrule
		$10^{-6}$ & \num{1.63e-2} & \num{1.42e-1} \\
		\midrule
		$10^{-7}$ & \num{2.97e-3} & \num{1.40e-1} \\
		\midrule
		$10^{-8}$ & \num{7.00e-4} & \num{1.40e-1} \\
		\bottomrule
	\end{tabularx}
	\caption{Water-air shock tube test case with the fine mesh: relative errors in the $l^{1}$ norm for the velocity as function of the ``residual'' volume fraction $\alpha_{1,R}$ at $t = T_{f}$.}
	\label{tab:relative_errors_water_air_shock_homogeneous_fine_vol_fraction}
\end{table}

\subsection{Epoxy-spinel shock}
\label{ssec:epoxy_spinel_homogeneous}

\begin{figure}[h!]
	\centering
	\begin{subfigure}{0.475\textwidth}
		\centering
		\includegraphics[width = 0.95\textwidth]{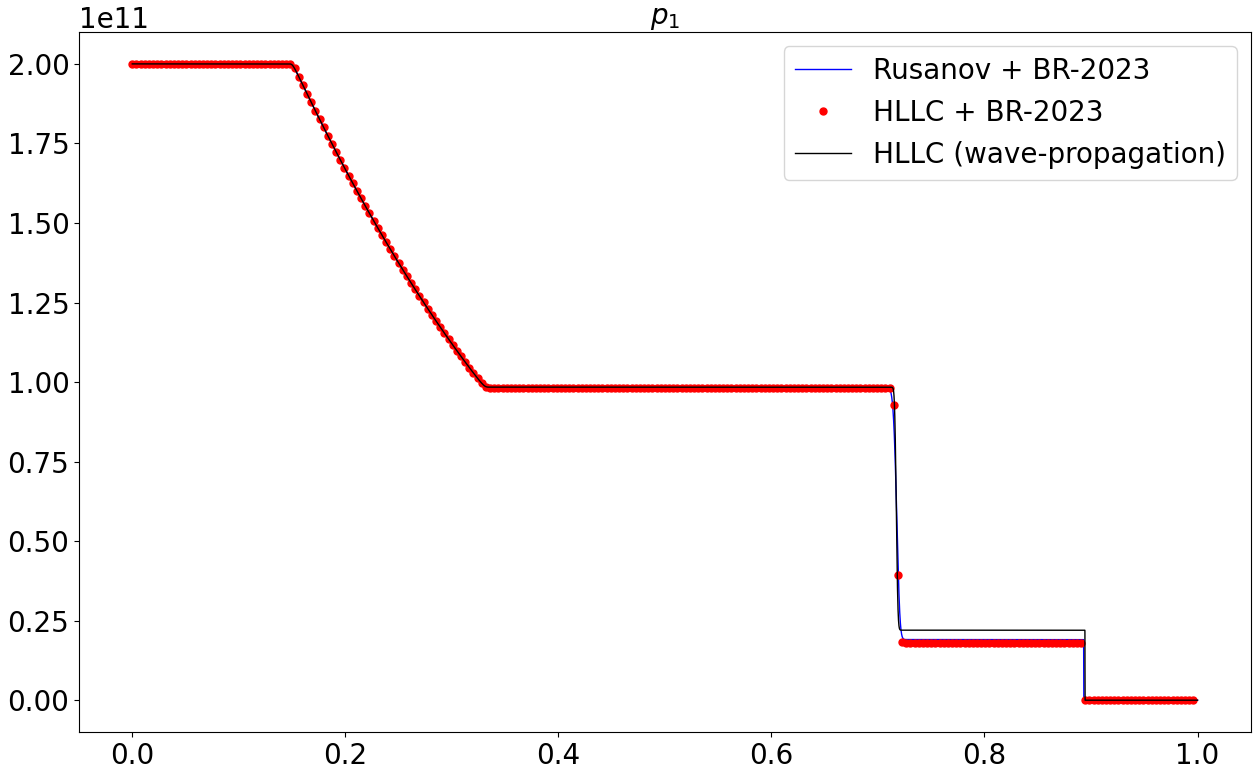}
	\end{subfigure}
	\begin{subfigure}{0.475\textwidth}
		\centering
		\includegraphics[width = 0.95\textwidth]{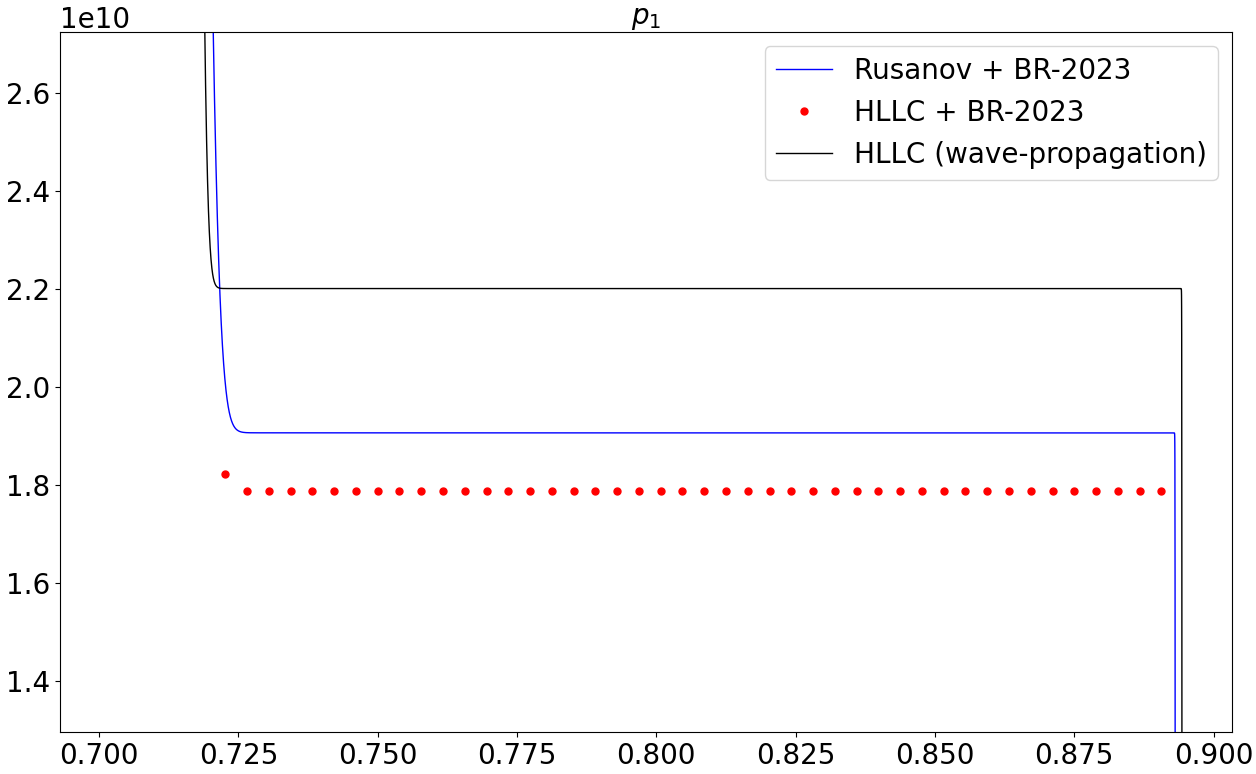}
	\end{subfigure}
	\begin{subfigure}{0.475\textwidth}
		\centering
		\includegraphics[width = 0.95\textwidth]{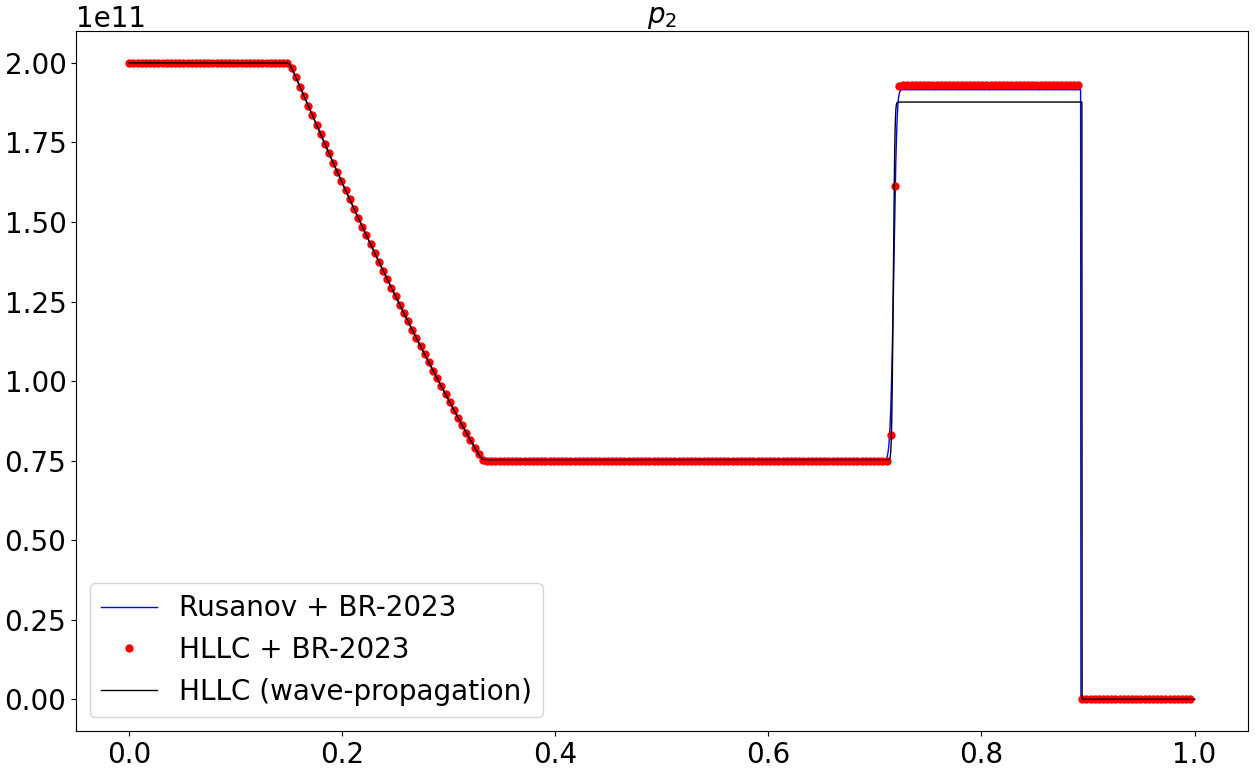}
	\end{subfigure}
	\begin{subfigure}{0.475\textwidth}
		\centering
		\includegraphics[width = 0.95\textwidth]{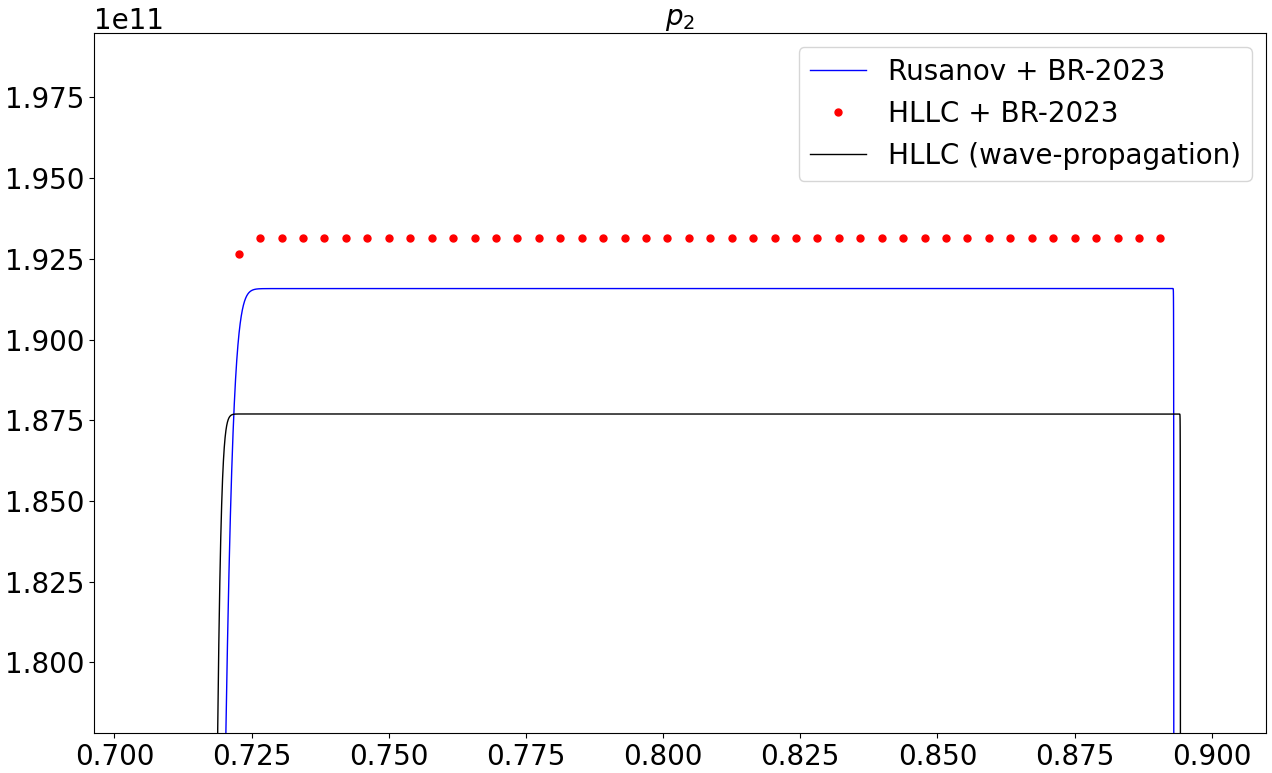}
	\end{subfigure}
	\caption{Epoxy-spinel shock test case with the fine mesh, results at $t = T_{f}$. Top-Left: pressure of phase $1$. Top-Right: zoom of the pressure of phase $1$. Bottom-Left: pressure of phase $2$. Bottom-Right: zoom of the pressure of phase $2$. Continuous black lines: HLLC-type wave-propagation method. Continuous blue lines: Rusanov scheme in combination with \textit{BR-2023} \eqref{eq:BR_Orlando_non_cons}. Red dots: HLLC scheme in combination with \textit{BR-2023} \eqref{eq:BR_Orlando_non_cons}.}
	\label{fig:epoxy_spinel_shock_homogeneous}
\end{figure}

In a final test, we consider the epoxy-spinel strong shock as presented in \cite{petitpas:2007}. We employ the fine mesh, so as to achieve mesh convergence. A reasonable agreement is established between the different numerical methods (Figure \ref{fig:epoxy_spinel_shock_homogeneous}). We point out that, in this test case, the non-conservative terms may play a more significant role compared to the previous test cases and the solution of the Riemann problem for the Euler--Euler equations does not provide a reference profile, since we are considering a mixture region instead of quasi-pure phases. However, one can easily notice visible differences in the phasic pressures in correspondence of the shock layer. We will further discuss this point in Section \ref{ssec:epoxy_spinel_relax}, where we show that these differences strongly influence the results obtained with the mechanical relaxation.

\subsection{Key findings on the numerical discretization of the 6-equation model}
\label{ssec:outcomes_6eqs}

We conclude this section by highlighting the main outcomes of our investigation on the numerical solution of the 6-equation model:
\begin{itemize}
	\item Spurious oscillations may arise around the contact discontinuity for the phasic pressures (Section \ref{ssec:sonic_rarefaction}). This occurs because, along the contact discontinuity, only the mixture pressure is preserved. This is a common issue for all numerical schemes, especially for HLLC-type schemes applied to the conservative part of the system and the volume fraction, combined with the discretization of the non-conservative terms in the energy equations. However, once mesh convergence is achieved, the HLLC wave-propagation scheme is the most effective at mitigating this issue.
	\item The HLLC scheme for the conservative part of the system and the volume fraction, combined with the \Crouzet\ approach \eqref{eq:Crouzet_non_cons} for the non-conservative terms of the energy equations, is highly sensitive to low-density flows (Section \ref{ssec:low_density}) and does not preserve in general thermodynamic admissibility. We refer to the upcoming Section \ref{sec:BR_path_conservative} for a possible explanation of this result.
	\item The 6-equation model exhibits a strong dependence on the discretization, even for mild shocks (Section \ref{ssec:water_air_shock_homogeneous}). This is likely related to the absence of a set of uniquely defined jump conditions. As we will see in Section \ref{sec:num_res_relax}, instantaneous relaxation generally yields consistent results, and an extreme pressure ratio is necessary to produce discretization-dependent solutions. Moreover, for the water-air shock test case in Section \ref{ssec:water_air_shock_homogeneous}, Rusanov-based schemes require a lower Courant number to achieve a physically meaningful solution.
\end{itemize}
Taking these considerations into account, among the considered methods, the HLLC wave-propagation scheme appears to be the most robust approach for mitigating the well-known issues related to the lack of uniquely defined jump conditions in the 6-equation model, especially when quasi-pure phases are involved.

\section{Reformulation of the BR-type approach in the framework of the path-conservative schemes}
\label{sec:BR_path_conservative}

In this section, we analyze the BR-type approach and show that it can be recast into the framework of path-conservative schemes. The path-conservative method \cite{pares:2006} has been originally developed for finite volume schemes \cite{castro:2006, pares:2006} and then extended to the Discontinuous Galerkin method \cite{dumbser:2009, dumbser:2010, gaburro:2024, fernandez:2022, rhebergen:2008}. We briefly recall some basic concepts. Consider again the non-conservative product $\mathbf{B}(w)\grad w$. Then, the path-conservative finite volume approximation on the element $K$ reads as follows
\begin{equation}
	\int_{K}\mathbf{B}(w)\grad w \mathrm{d}\Omega \approx \sum_{\Gamma in \partial K}\int_{\Gamma}\boldsymbol{\mathcal{D}}^{-}\rpth{w_{L}, w_{R}, \bm{n}} \mathrm{d}\Gamma,
\end{equation}
where $\bm{n}$ denotes the outward unit normal vector from the element. Following \cite{pares:2006}, in order to obtain a $\Psi$-conservative scheme, the function $\boldsymbol{\mathcal{D}}$ is such that
\begin{subequations}
	\begin{align}
		&\boldsymbol{\mathcal{D}}^{-}\rpth{w, w, \bm{n}} = 0 \quad \forall w, \bm{n}, \label{eq:PC_consistency} \\
		&\boldsymbol{\mathcal{D}}^{-}\rpth{w_{L}, w_{R}, \bm{n}} + \boldsymbol{\mathcal{D}}^{+}\rpth{w_{L}, w_{R}, \bm{n}} = \spth{\int_{0}^{1} \mathbf{B}\rpth{\Psi(w_{L}, w_{R}, \hat{s})}\pad{\Psi}{\hat{s}}d\hat{s}}\bm{n}, \label{eq:PC_conservative}
	\end{align}
\end{subequations}
where $\boldsymbol{\mathcal{D}}^{+}\rpth{w_{L}, w_{R}, \bm{n}} = \boldsymbol{\mathcal{D}}^{-}\rpth{w_{R}, w_{L}, -\bm{n}}$. The function $\Psi$ denotes a family of paths connecting $w_{L}$ and $w_{R}$ such that
\begin{equation*}
	\Psi(w_{L}, w_{R}, 0) = w_{L}, \qquad \Psi(w_{L}, w_{R}, 1) = w_{R}, \qquad \Psi(w, w, \hat{s}) = w.
\end{equation*}
We refer to \cite{pares:2006} for a detailed discussion of the regularity required by $\Psi$. For multidimensional systems, the path connecting the
two states $w_{L}$ and $w_{R}$ may in principle depend also on $\bm{n}$, but if the system is invariant under rotations, this dependency can be dropped \cite{castro:2009, dumbser:2009}. Condition \eqref{eq:PC_consistency} guarantees that if $\grad w = 0$, then no contribution arises. Condition \eqref{eq:PC_conservative} guarantees that if $\mathbf{B}$ represents the Jacobian of a flux, then a classical conservative scheme is retrieved \cite{pares:2006}. We stress the fact that the wave-propagation scheme discussed in Section \ref{ssec:disc_scheme} fits into this formalism. Denote now
\begin{equation}
	\boldsymbol{\mathcal{D}}^{-}\rpth{w_{L}, w_{R}, \bm{n}} = \rpth{\widehat{\mathbf{B}w}(w_{L}, w_{R}) - \widehat{\mathbf{B}}(w_{L}, w_{R})w_{L}}\bm{n},
\end{equation}
with $\widehat{\mathbf{B}w}(w_{L}, w_{R})$ and $\widehat{\mathbf{B}}(w_{L}, w_{R})$ as in \eqref{eq:BR}. One can immediately notice that condition \eqref{eq:PC_consistency} is fulfilled. Moreover,
\begin{equation}
	\boldsymbol{\mathcal{D}}^{-}\rpth{w_{L}, w_{R}, \bm{n}} + \boldsymbol{\mathcal{D}}^{+}\rpth{w_{L}, w_{R}, \bm{n}} = \spth{\widehat{\mathbf{B}}(w_{L}, w_{R})\rpth{w_{R} - w_{L}}}\bm{n}.
\end{equation}
If we consider a family of segments
\begin{equation}
	\Psi = w_{L} + \hat{s}\rpth{w_{R} - w_{L}},
\end{equation}
condition \eqref{eq:PC_conservative} reduces to
\begin{equation}
	\boldsymbol{\mathcal{D}}^{-}\rpth{w_{L}, w_{R}, \bm{n}} + \boldsymbol{\mathcal{D}}^{+}\rpth{w_{L}, w_{R}, \bm{n}} = \rpth{w_{R} - w_{L}}\spth{\int_{0}^{1}\mathbf{B}\rpth{\Psi(w_{L}, w_{R}, \hat{s})}\mathrm{d}\hat{s}}\bm{n},
\end{equation}
and therefore, considering $\widehat{\mathbf{B}}$ as an approximation of $\int_{0}^{1} \mathbf{B}\rpth{\Psi(w_{L}, w_{R}, \hat{s})}\mathrm{d}\hat{s}$, we conclude that also condition \eqref{eq:PC_conservative} is fulfilled. In particular, if we take
\begin{equation}
	\widehat{\mathbf{B}}(w_{L}, w_{R}) = \frac{1}{2}\rpth{\mathbf{B}(w_{L}) + \mathbf{B}(w_{R})},
\end{equation}
this is equivalent to consider a trapezoidal rule for $\int_{0}^{1} \mathbf{B}\rpth{\Psi(w_{L}, w_{R}, \hat{s})}\mathrm{d}\hat{s}$. The choice of a linear path is the standard one in the literature \cite{abgrall:2010, busto:2021, gaburro:2024, gatti:2024, fernandez:2022}, even though more complex paths can be considered \cite{delorenzo:2018, schwendeman:2006}, and also the choice of the trapezoidal rule to approximate the path integral is rather common \cite{busto:2021, gatti:2024}.

In conclusion, the BR-type approach can be seen as a particular approximate path-conservative scheme employing segments as paths. This holds independently of the specific definition of $\widehat{\mathbf{B}}(w_{L}, w_{R})$. The BR-type approach has already been employed for the discretization of non-conservative terms in hyperbolic systems \cite{orlando:2023, tumolo:2015} and the idea to consider integration by parts twice \cite{kopriva:2010} has been already adopted in several CFD applications to obtain well-balanced schemes in the so-called strong form \cite{abgrall:2017, arpaia:2022, arpaia:2026, vater:2019}. However, to the best of our knowledge, this is the first rigorous analysis that shows how a numerical strategy originally and typically employed for the discretization of diffusion operators \cite{bassi:1997a, ortleb:2020} is linked to a widely used framework for the discretization of hyperbolic non-conservative systems.

We conclude this section by performing a similar analysis for the \Crouzet\ approach \eqref{eq:Crouzet_non_cons}. We immediately notice that, as presented in \eqref{eq:Crouzet_non_cons}, it does not satisfy \eqref{eq:PC_consistency}. However, similarly to what we discussed for the wave-propagation scheme in Section \ref{ssec:disc_scheme}, we can rewrite \eqref{eq:discrete_scheme}, at least in the one-dimensional case, as
\begin{align}
	\mathbf{q}_{j}^{n+1} &= \mathbf{q}_{j}^{n} - \frac{\Delta t}{\Delta x}\bigg(\spth{\boldsymbol{\mathcal{H}}^{-}_{j + 1/2}\rpth{\mathbf{q}_{j}^{n}, \mathbf{q}_{j+1}^{n}} - \mathbf{F}(\mathbf{q}_{j}^{n}) - \boldsymbol{\mathcal{T}}^{-}(\mathbf{q}_{j}^{n}, \mathbf{q}_{j}^{n})} - \nonumber \\ &\phantom{= \mathbf{q}_{j}^{n} - \frac{\Delta t}{\Delta x}\bigg(}\spth{\boldsymbol{\mathcal{H}}^{+}_{j - 1/2}\rpth{\mathbf{q}_{j-1}^{n}, \mathbf{q}_{j}^{n}} - \mathbf{F}(\mathbf{q}_{j}^{n}) - \boldsymbol{\mathcal{T}}^{+}(\mathbf{q}_{j}^{n}, \mathbf{q}_{j}^{n})}\bigg).
\end{align}
Hence, returning to our non-generic conservative product, we can define
\begin{equation}
	\boldsymbol{\mathcal{D}}^{-}\rpth{w_{L}, w_{R}, \bm{n}} = \rpth{\mathbf{B}(w_{L})\ave{w} - \mathbf{B}(w_{L})w_{L}}\bm{n} = \frac{\mathbf{B}(w_{L})}{2}\rpth{w_{R} - w_{L}}\bm{n},
\end{equation}
so that,
\begin{equation}
	\boldsymbol{\mathcal{D}}^{-}\rpth{w_{L}, w_{R}, \bm{n}} + \boldsymbol{\mathcal{D}}^{+}\rpth{w_{L}, w_{R}, \bm{n}} = \rpth{w_{R} - w_{L}}\frac{\mathbf{B}(w_{L}) + \mathbf{B}(w_{R})}{2}\bm{n}.
\end{equation}
Hence, also for this specific scheme, we obtain that it corresponds, at least in the one-dimensional case, to considering a linear path with the trapezoidal rule to approximate $\int_{0}^{1} \mathbf{B}\rpth{\Psi(w_{L}, w_{R}, \hat{s})}\mathrm{d}\hat{s}$. The low robustness of this discretization for the low-density test case in Section \ref{ssec:low_density} does not therefore depend on the fact that it does not fit within the path-conservative framework, but is likely related to a lower numerical diffusion, in particular for the non-conservative terms of the energy equations. We aim to further investigate the treatment of non-conservative terms in a future dedicated work.

\section{Instantaneous mechanical relaxation}
\label{sec:relax}

In this section, we discuss the relaxation of the phasic pressures obtained from the discretization of the hyperbolic operator towards an equilibrium value. Following the operator splitting approach depicted at the beginning of Section \ref{sec:num_hyp}, we obtain the following system of ODEs
\begin{subequations}
	\begin{align}
		\odif{\alpha_{1}}{t} &= \mu\rpth{p_{1} - p_{2}}, \label{eq:alpha_relax} \\
		\odif{\alpha_{1}\rho_{1}}{t} &= 0, \label{eq:m1_relax} \\
		\odif{\alpha_{2}\rho_{2}}{t} &= 0, \label{eq:m2_relax} \\
		\odif{\rho\vel}{t} &= 0, \label{eq:u_relax} \\
		\odif{\alpha_{1}\rho_{1}E_{1}}{t} &= -\mu p_{I}\rpth{p_{1} - p_{2}}, \label{eq:E1_relax} \\
		\odif{\alpha_{2}\rho_{2}E_{2}}{t} &= \mu p_{I}\rpth{p_{1} - p_{2}}. \label{eq:E2_relax}
	\end{align}
\end{subequations}
Combining \eqref{eq:alpha_relax}, \eqref{eq:E1_relax}, and \eqref{eq:E2_relax}, one obtains \cite{pelanti:2014}
\begin{subequations}
	\begin{align}
		\odif{\alpha_{1}\rho_{1}E_{1}}{t} &= -p_{I}\odif{\alpha_{1}}{t}, \\
		\odif{\alpha_{2}\rho_{2}E_{2}}{t} &= p_{I}\odif{\alpha_{1}}{t},
	\end{align}
\end{subequations}
which, owing to \eqref{eq:m1_relax}, \eqref{eq:m2_relax}, and \eqref{eq:u_relax}, and, with a slight abuse of notation, can be rewritten as,
\begin{subequations}
	\begin{align}
		\alpha_{1}\rho_{1}\odif{e_{1}}{\alpha_{1}} &= -p_{I}, \label{eq:e1_relax_simpl} \\
		\alpha_{2}\rho_{2}\odif{e_{2}}{\alpha_{2}} &= p_{I}. \label{eq:e2_relax_simpl}
	\end{align}
\end{subequations}
Here, we choose the following expression for the interfacial pressure
\begin{equation}\label{eq:choice_pI}
	p_{I} = \frac{Z_{2}p_{1} + Z_{1}p_{2}}{Z_{1} + Z_{2}},
\end{equation}
where $Z_{k} = \rho_{k}c_{k}$ corresponds to the phasic acoustic impedance. Relation \eqref{eq:choice_pI} arises from the solution of the linearized Riemann problem for the Euler equations in the single-velocity limit. The first proposal to use the expression from the solution of the linearized Riemann problem for $\vel_{I}$ and $p_{I}$ appears in \cite{saurel:2003}, in the context of Baer--Nunziato-type models derived via the Discrete Equation Method \cite{abgrall:2003}.

In the following, for a given quantity $\chi$, we will denote by $\chi^{(0)}$ the value coming from the discretization of the hyperbolic operator, i.e. before the relaxation step, and by $\chi^{\star}$ the value at the mechanical equilibrium, i.e. after the relaxation process. Applying the midpoint rule for the integral approximation of the interfacial pressure in \eqref{eq:e1_relax_simpl}-\eqref{eq:e2_relax_simpl}, as proposed in \cite{pelanti:2014}, one obtains
\begin{subequations}
	\begin{align}
		e_{1}^{\star} - e_{1}^{(0)} &= -\frac{p_{I}^{(0)} + p_{I}^{\star}}{2}\frac{\alpha_{1}^{\star} - \alpha_{1}^{(0)}}{\alpha_{1}^{(0)}\rho_{1}^{(0)}}, \label{eq:e1_relax_integral} \\
		e_{2}^{\star} - e_{2}^{(0)} &= \frac{p_{I}^{(0)} + p_{I}^{\star}}{2}\frac{\alpha_{1}^{\star} - \alpha_{1}^{(0)}}{\alpha_{2}^{(0)}\rho_{2}^{(0)}}. \label{eq:e2_relax_integral}
	\end{align}
\end{subequations}
Recall that $\alpha_{k}^{(0)}\rho_{k}^{(0)} = \alpha_{k}^{(\star)}\rho_{k}^{(\star)}$ thanks to \eqref{eq:m1_relax}-\eqref{eq:m2_relax}. Next, imposing instantaneous mechanical equilibrium, i.e. $p_{1}^{\star} = p_{2}^{\star} = p^{\star}$, \eqref{eq:e1_relax_integral}-\eqref{eq:e2_relax_integral} reduce to a system of two equations in the two unknowns $p^{\star}$ and $\alpha_{1}^{\star}$. The mixture specific internal energy $e$ is defined as
\begin{equation}
	\rho e = \alpha_{1}\rho_{1}e_{1} + \alpha_{2}\rho_{2}e_{2}.
\end{equation}
The latter relation, assuming $p_{1} = p_{2} = p$ in the phasic pressure laws, determines implicitly the mixture pressure law 
$$p = p\rpth{e, \rho_{1}, \rho_{2}, \alpha_{1}}.$$ 
In the case of the SG-EOS, we obtain an explicit expression for the mixture pressure \cite{allaire:2002, pelanti:2014, saurel:2009}
\begin{equation}\label{eq:mixture_pressure_law_SG}
	p = p\rpth{e, \rho_{1}, \rho_{2}, \alpha_{1}} = \frac{\rho e - \rpth{\alpha_{1}\rho_{1}\eta_{1} + \alpha_{2}\rho_{2}\eta_{2}} - \rpth{\frac{\alpha_{1}\gamma_{1}\pi_{1}}{\gamma_{1} - 1} + \frac{\alpha_{2}\gamma_{2}\pi_{2}}{\gamma_{2} - 1}}}{\frac{\alpha_{1}}{\gamma_{1} - 1} + \frac{\alpha_{2}}{\gamma_{2} - 1}}.
\end{equation}
It is worth to remark that, thanks to the use of phasic total energies for the homogeneous system \eqref{eq:system_hyperbolic}, the numerical method guarantees that the relaxed pressure $p^{\star}$ verifies the mixture pressure law and satisfies by construction the conservation of the mixture energy, i.e. it is mixture-energy consistent \cite{pelanti:2014}. In the specific case of the SG-EOS, a quadratic equation of the form 
$$a\rpth{p^{\star}}^{2} + b p^{\star} + c = 0,$$
is obtained for the relaxed pressure $p^{\star}$, with
\begin{align}
	a =& 1 + \gamma_{2}\alpha_{1}^{(0)} + \gamma_{1}\alpha_{2}^{(0)}, \nonumber \\
	b =& \rpth{2\gamma_{1}\pi_{1} + \rpth{\gamma_{1} - 1}p_{I}^{(0)}}\alpha_{2}^{(0)} + \rpth{2\gamma_{2}\pi_{2} + \rpth{\gamma_{2} - 1}p_{I}^{(0)}}\alpha_{1}^{(0)} \nonumber \\
	&- \rpth{1 + \gamma_{2}}\alpha_{1}^{(0)}p_{1}^{(0)} - \rpth{1 + \gamma_{1}}\alpha_{2}^{(0)}p_{2}^{(0)}, \\
	c =& -\rpth{2\gamma_{2}\pi_{2} + \rpth{\gamma_{2} - 1}p_{I}^{(0)}}\alpha_{1}^{(0)}p_{1}^{(0)} - \rpth{2\gamma_{1}\pi_{1} + \rpth{\gamma_{1} - 1}p_{I}^{(0)}}\alpha_{2}^{(0)}p_{2}^{(0)}. \nonumber
\end{align}
Moreover, the value of the volume fraction at the equilibrium $\alpha_{1}^{\star}$ is \cite{pelanti:2014}
\begin{equation}
	\alpha_{1}^{\star} = \alpha_{1}^{(0)}\frac{\rpth{\gamma_{1} - 1}p^{\star} + 2p_{1}^{(0)} + 2\gamma_{1}\pi_{1} + \rpth{\gamma_{1} - 1}p_{I}^{(0)}}{\rpth{\gamma_{1} + 1}p^{\star} + 2\gamma_{1}\pi_{1} + \rpth{\gamma_{1} - 1}p_{I}^{(0)}}.
\end{equation}

As evident from this short presentation, the relaxation process depends on the choice of interfacial pressure $p_{I}$. Moreover, once a choice of the interfacial pressure is made, several approximation techniques for the relaxation process are available. Both the choice of the interfacial pressure and the numerical discretization of the relaxation operator may affect the post-relaxation state and impact the global robustness and accuracy of the numerical scheme. We refer to \cite{haegeman:2024} for a detailed discussion of the impact of different choices for the interfacial pressure and its numerical analysis. In the present work, we consider a single choice interfacial pressure given by \eqref{eq:choice_pI} and the corresponding numerical relaxation procedure described above. We remark that in the current method, without any form of implicitation, i.e. if $p_{I}^{\star}$ is not employed in the definition of $p_{I}$ for \eqref{eq:e1_relax_integral}-\eqref{eq:e2_relax_integral}, the positivity of the thermodynamic state may be lost for the proposed relaxation methods \cite{haegeman:2024}.

\section{Numerical results with mechanical relaxation}
\label{sec:num_res_relax}

In this section, we analyze the approaches for the instantaneous mechanical relaxation depicted in Section \ref{sec:relax} and their interaction with the numerical strategies for the hyperbolic operator.

\subsection{Water-air shock tube}
\label{ssec:water_air_shock_relax}

\begin{figure}[h!]
	\centering
	\begin{subfigure}{0.475\textwidth}
		\centering
		\includegraphics[width = 0.95\textwidth]{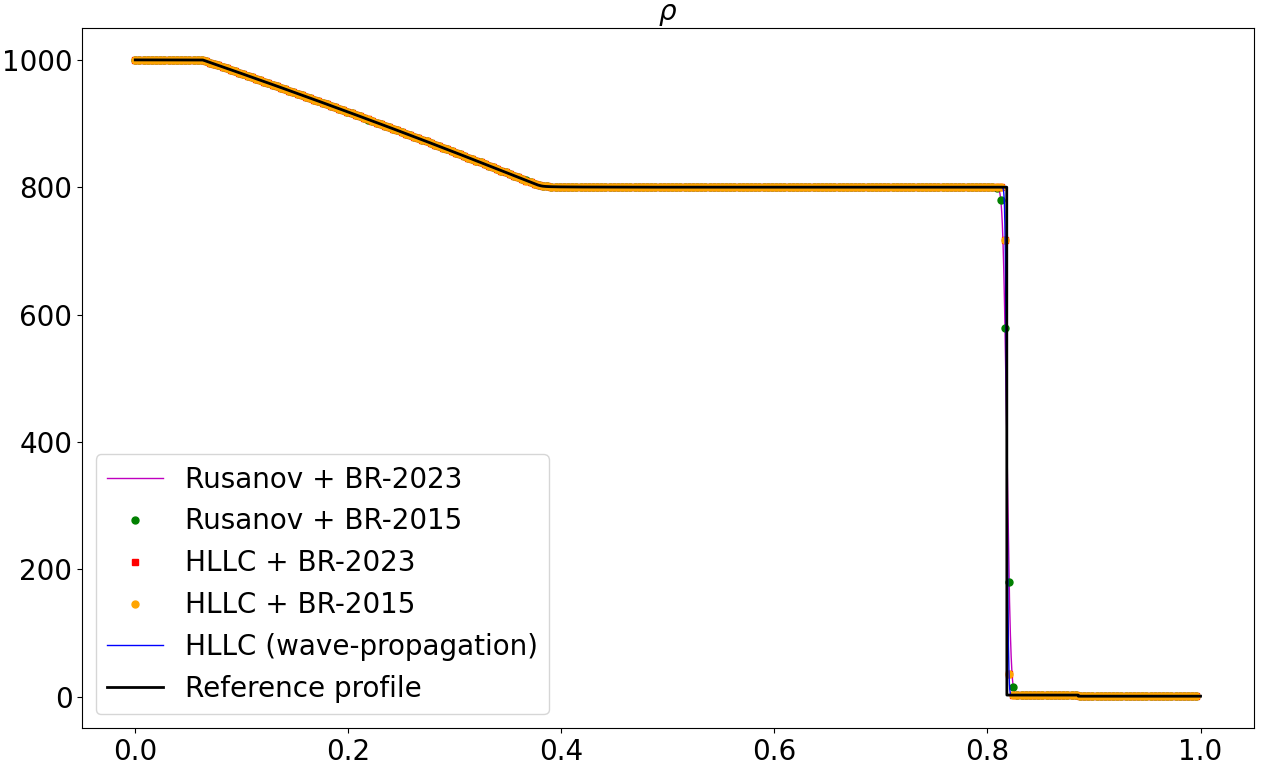}
	\end{subfigure}
	\begin{subfigure}{0.475\textwidth}
		\centering
		\includegraphics[width = 0.95\textwidth]{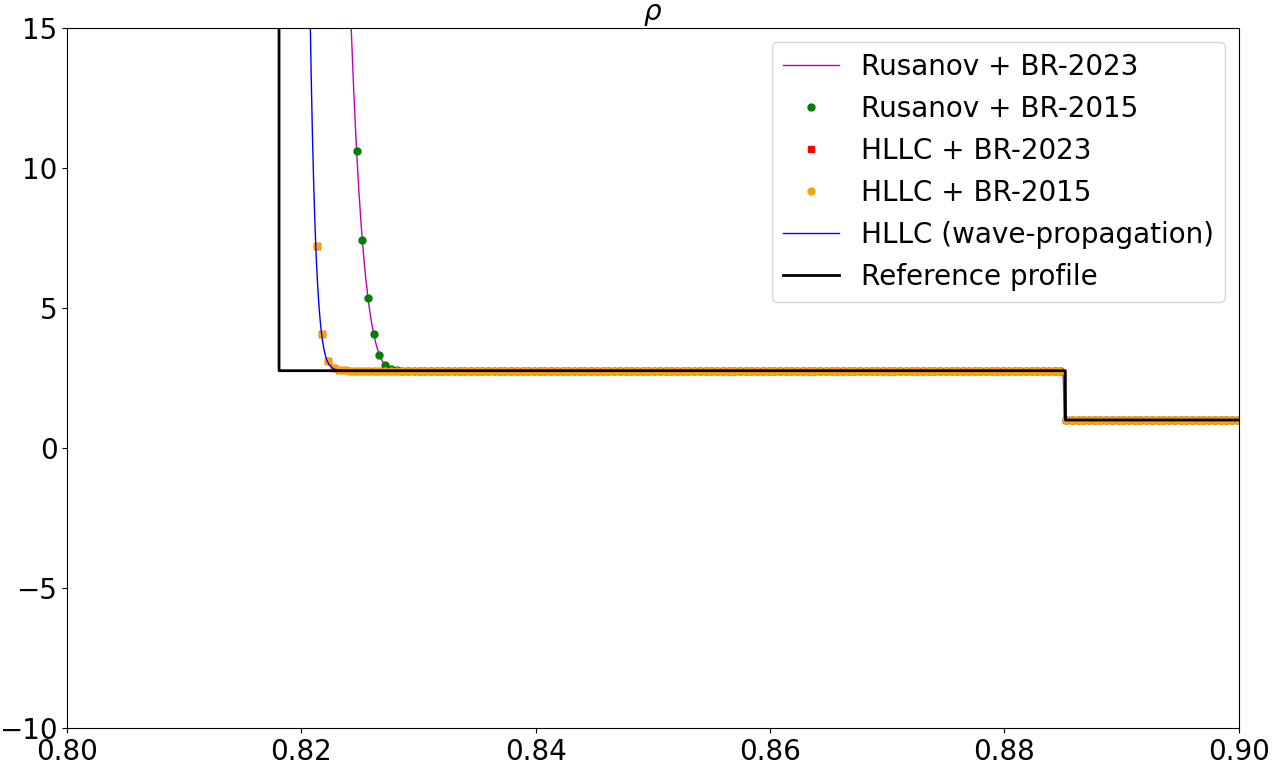}
	\end{subfigure}
	\begin{subfigure}{0.475\textwidth}
		\centering
		\includegraphics[width = 0.95\textwidth]{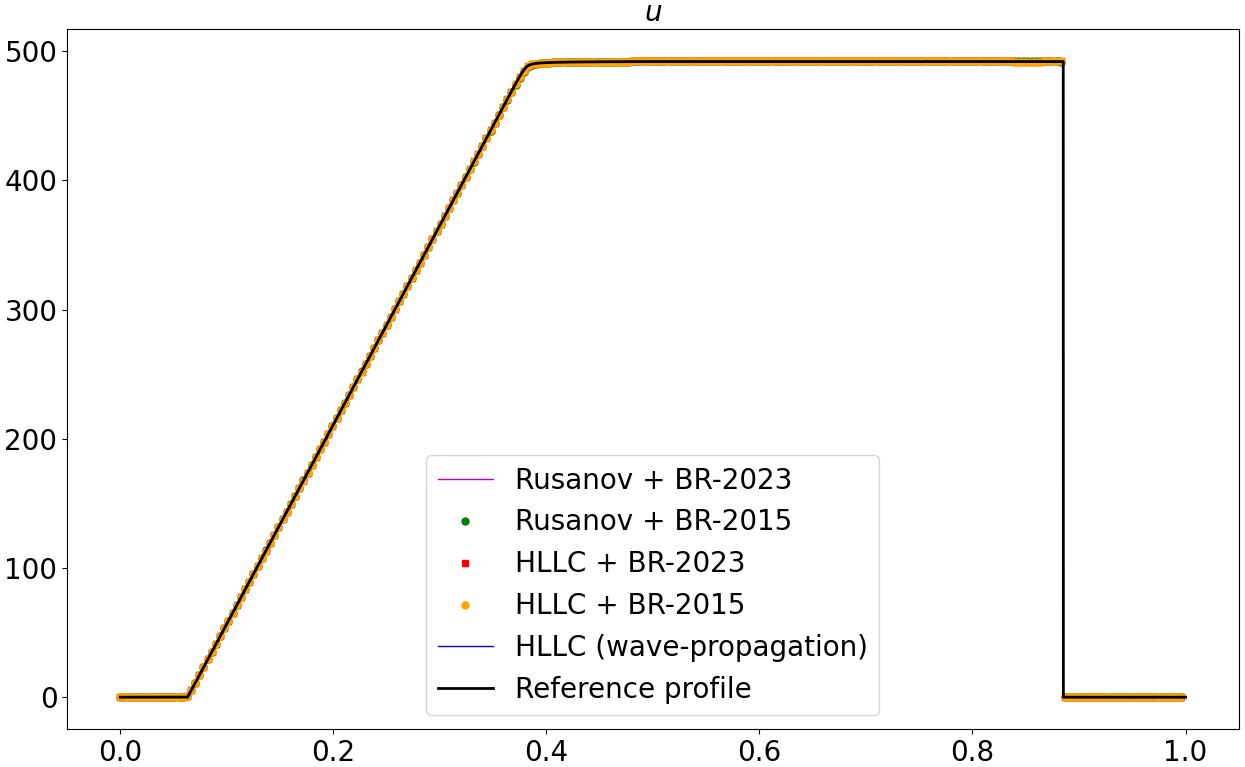}
	\end{subfigure}
	\begin{subfigure}{0.475\textwidth}
		\centering
		\includegraphics[width = 0.95\textwidth]{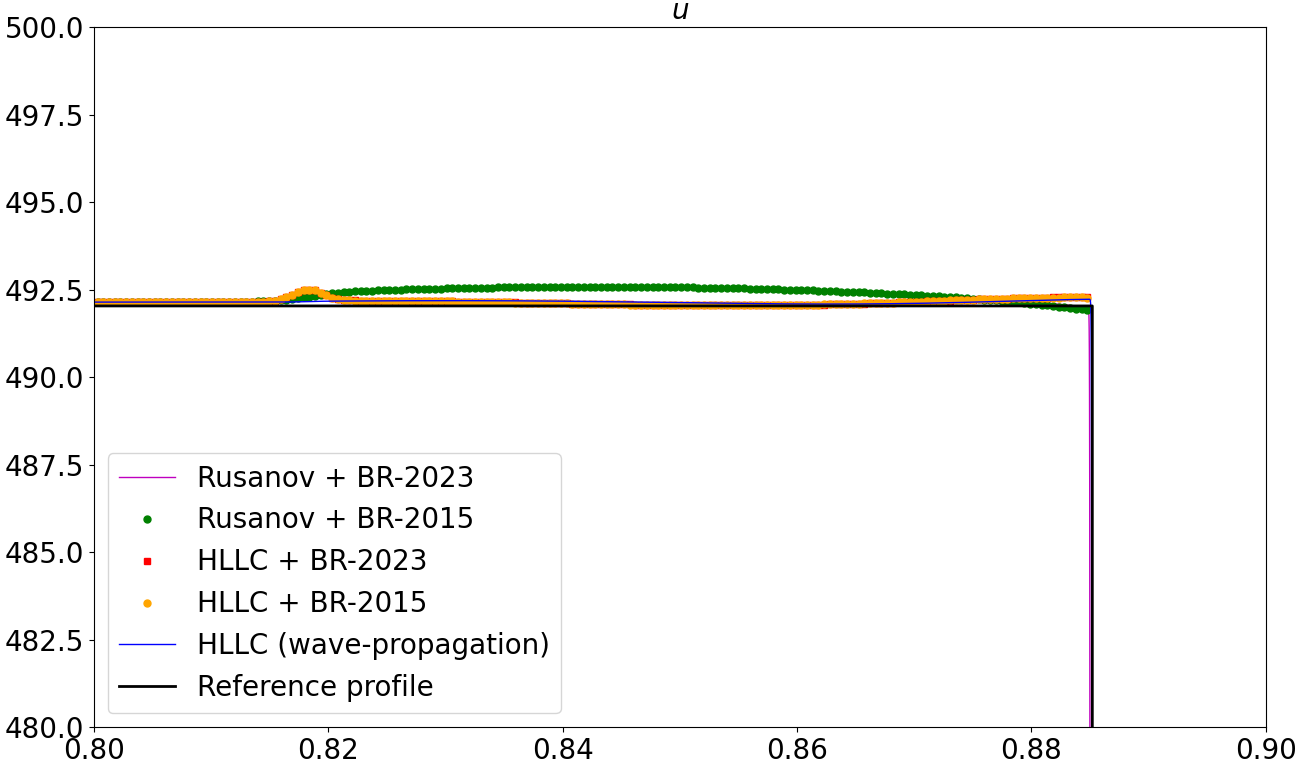}
	\end{subfigure}
	\caption{Water-air shock tube test case with instantaneous mechanical relaxation and fine mesh, results at $t = T_{f}$. Top-left: mixture density. Top-right: zoom of the mixture density. Bottom-left: velocity. Bottom-right: zoom of the velocity. Continuous black lines: solution obtained using the approximated jump conditions in \cite{petitpas:2007}. Continuous blue lines: HLLC-type wave-propagation scheme. Continuous magenta lines: Rusanov flux in combination with \textit{BR-2023} \eqref{eq:BR_Orlando_non_cons}. Green dots: Rusanov flux in combination with \textit{BR-2015} \eqref{eq:BR_Tumolo_non_cons}. Red squares: HLLC flux in combination with \textit{BR-2023} \eqref{eq:BR_Orlando_non_cons}. Orange dots: HLLC flux in combination with \textit{BR-2015} \eqref{eq:BR_Tumolo_non_cons}.}
	\label{fig:water_air_shock_relax_instantaneous}
\end{figure}

As a first test case of this section, we consider the instantaneous mechanical relaxation for the water-air shock tube described in Section \ref{ssec:water_air_shock_homogeneous}. We employ the fine mesh. The solution obtained considering the approximated jump relations for the 5-equation model of Kapila is reported for reference. A good agreement is established for all the schemes, in particular for what concerns the HLLC-type schemes (Figure \ref{fig:water_air_shock_relax_instantaneous}). Notice that, when the mechanical relaxation is considered, we can establish a physically meaningful solution with the Rusanov flux also at $C = 0.9$. One can easily notice that the HLLC-type wave-propagation scheme and the HLLC for the conservative portion of the system in combination with the treatment of the non-conservative terms converge towards the same numerical solution, even though the HLLC-type wave-propagation scheme exhibits slightly lower errors, when these are measured with respect to the solution obtained from the 5-equation model with approximated jump relations (Table \ref{tab:relative_errors_water_air_shock_relaxed_fine}).

\begin{table}[h!]
	\centering
	\footnotesize
	\begin{tabularx}{0.75\columnwidth}{lXXXX}
		\toprule
		\multirow{2}{*}{\textbf{Scheme}} & \multicolumn{3}{c}{$l^{1}$ norm rel. error} & \\
		\cmidrule(l){2-5}
		& $\alpha_{1}$ & $u$ & $\rho$ & $p$ \\
		\midrule
		Rusanov + \textit{BR-2023} & \num{2.80e-3} & \num{6.14e-4} & \num{2.66e-3} & \num{6.04e-4} \\
		\midrule
		Rusanov + \textit{BR-2015} & \num{2.80e-3} & \num{6.15e-4} & \num{2.66e-3} & \num{6.06e-4} \\
		\midrule
		HLLC + \textit{BR-2023} & \num{1.29e-3} & \num{4.55e-4} & \num{1.24e-3} & \num{3.56e-4} \\
		\midrule
		HLLC + \textit{BR-2015} & \num{1.29e-3} & \num{4.50e-4} & \num{1.24e-3} & \num{3.56e-4} \\
		\midrule
		HLLC (wave-propagation) & \num{1.20e-3} & \num{3.24e-4} & \num{1.15e-3} & \num{3.57e-4} \\
		\bottomrule
	\end{tabularx}
	\caption{Water-air shock tube test case with instantaneous mechanical relaxation and fine mesh: relative errors in the $l^{1}$ norm for volume fraction, velocity, mixture density, and mixture pressure at $t = T_{f}$.}
	\label{tab:relative_errors_water_air_shock_relaxed_fine}
\end{table}

\subsection{Water cavitation tube}
\label{ssec:cavitation}

\begin{figure}[h!]
	\centering
	\begin{subfigure}{0.475\textwidth}
		\centering
		\includegraphics[width = 0.95\textwidth]{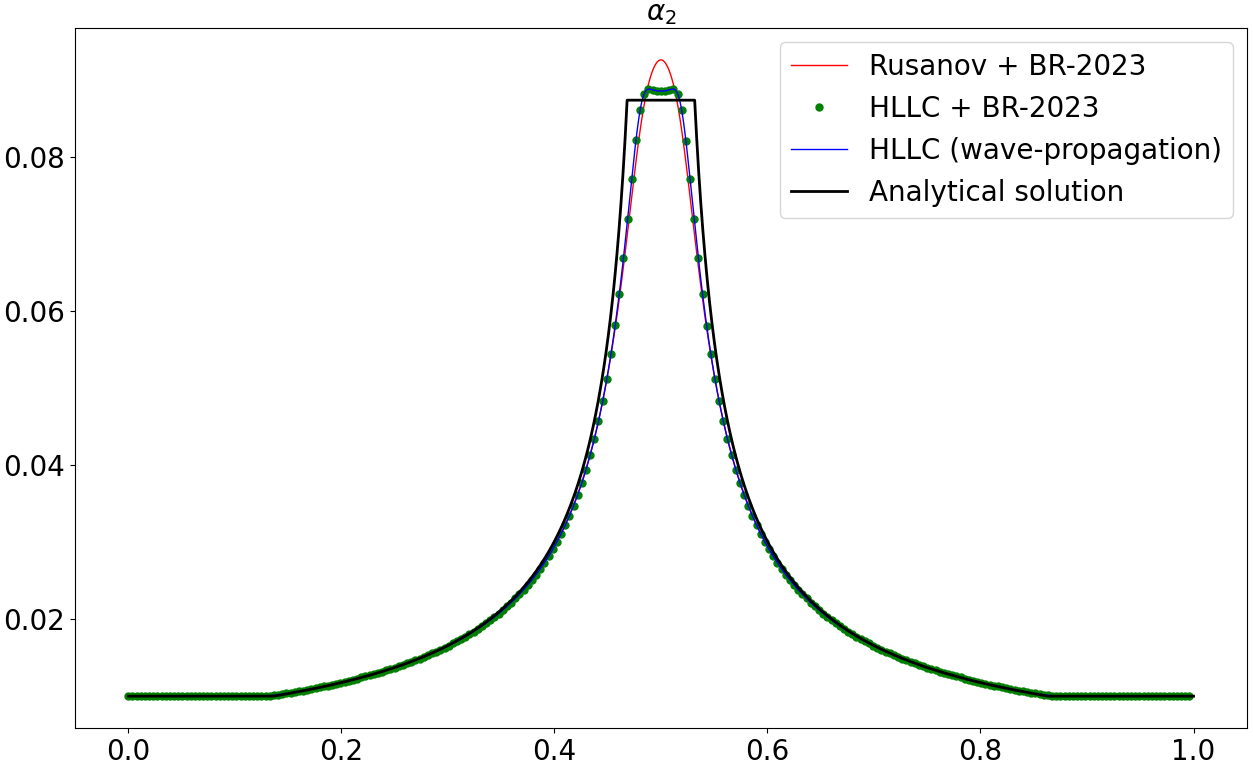}
	\end{subfigure}
	\begin{subfigure}{0.475\textwidth}
		\centering
		\includegraphics[width = 0.95\textwidth]{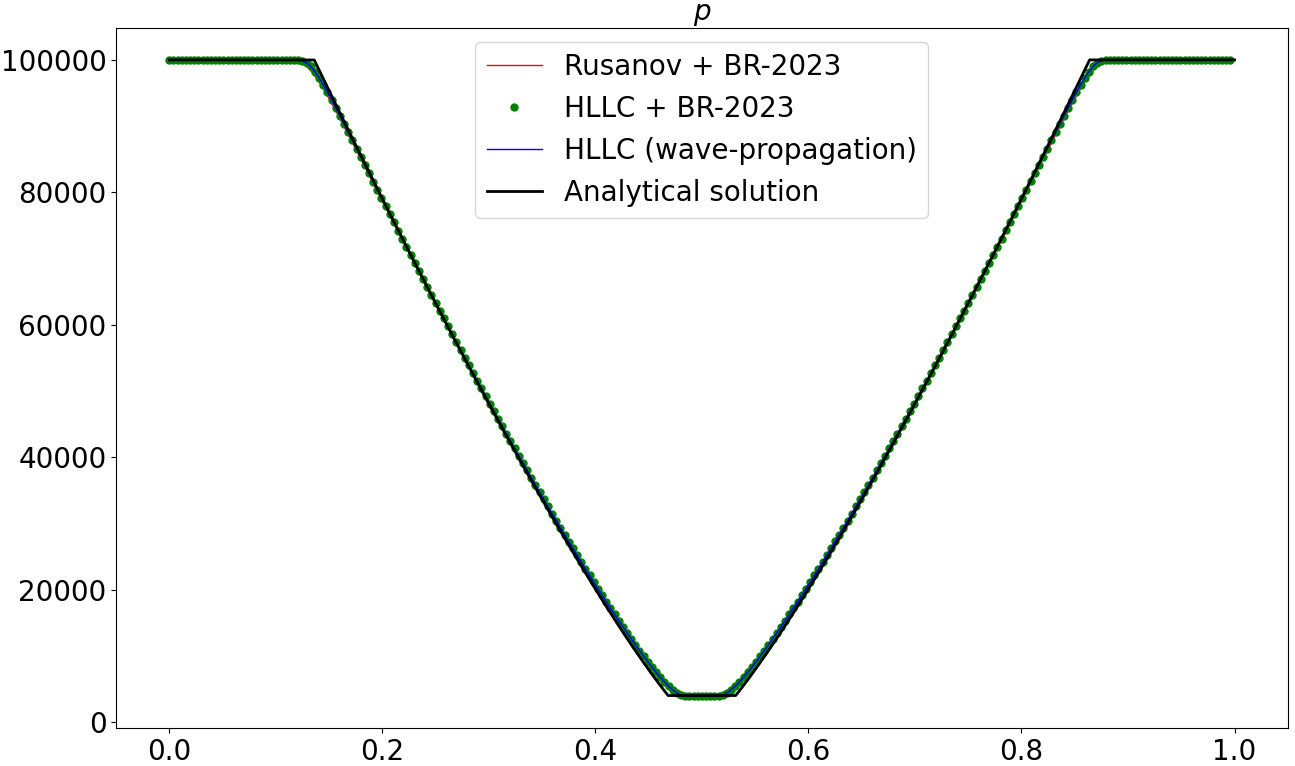}
	\end{subfigure}
	\begin{subfigure}{0.475\textwidth}
		\centering
		\includegraphics[width = 0.95\textwidth]{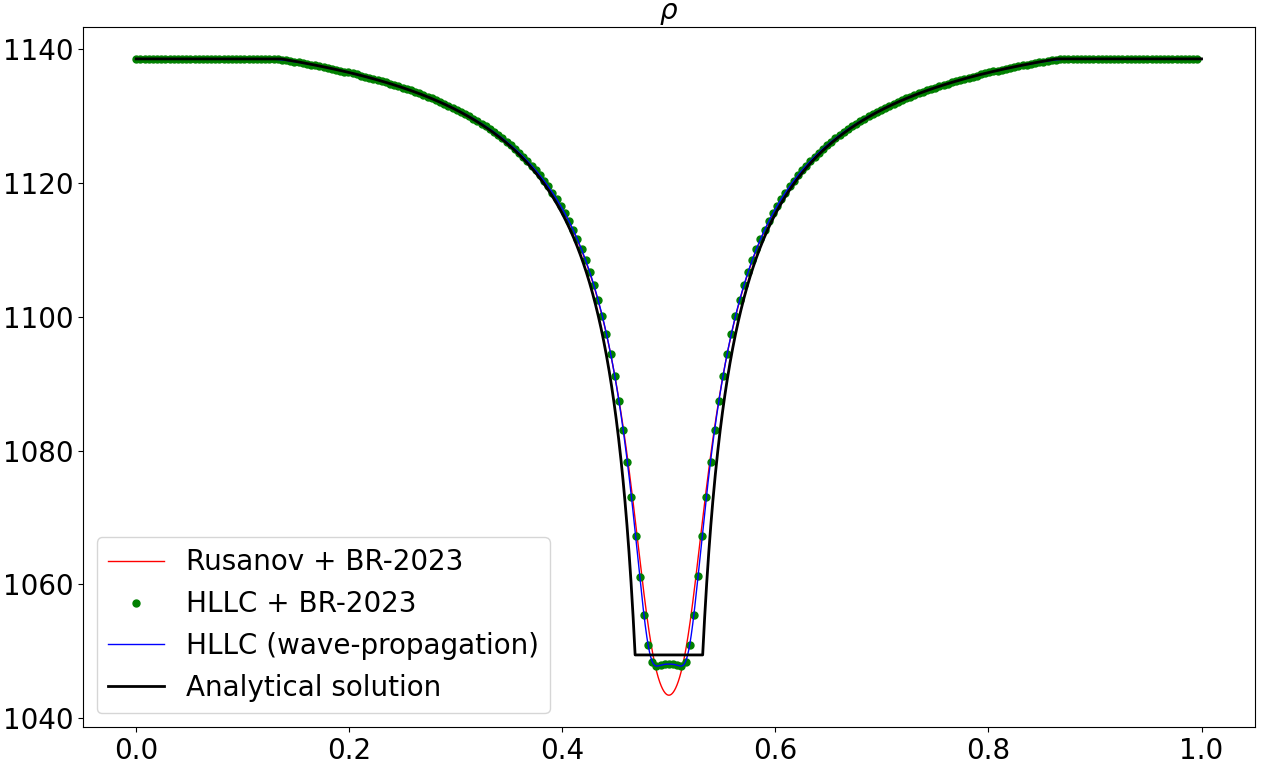}
	\end{subfigure}
	\begin{subfigure}{0.475\textwidth}
		\centering
		\includegraphics[width = 0.95\textwidth]{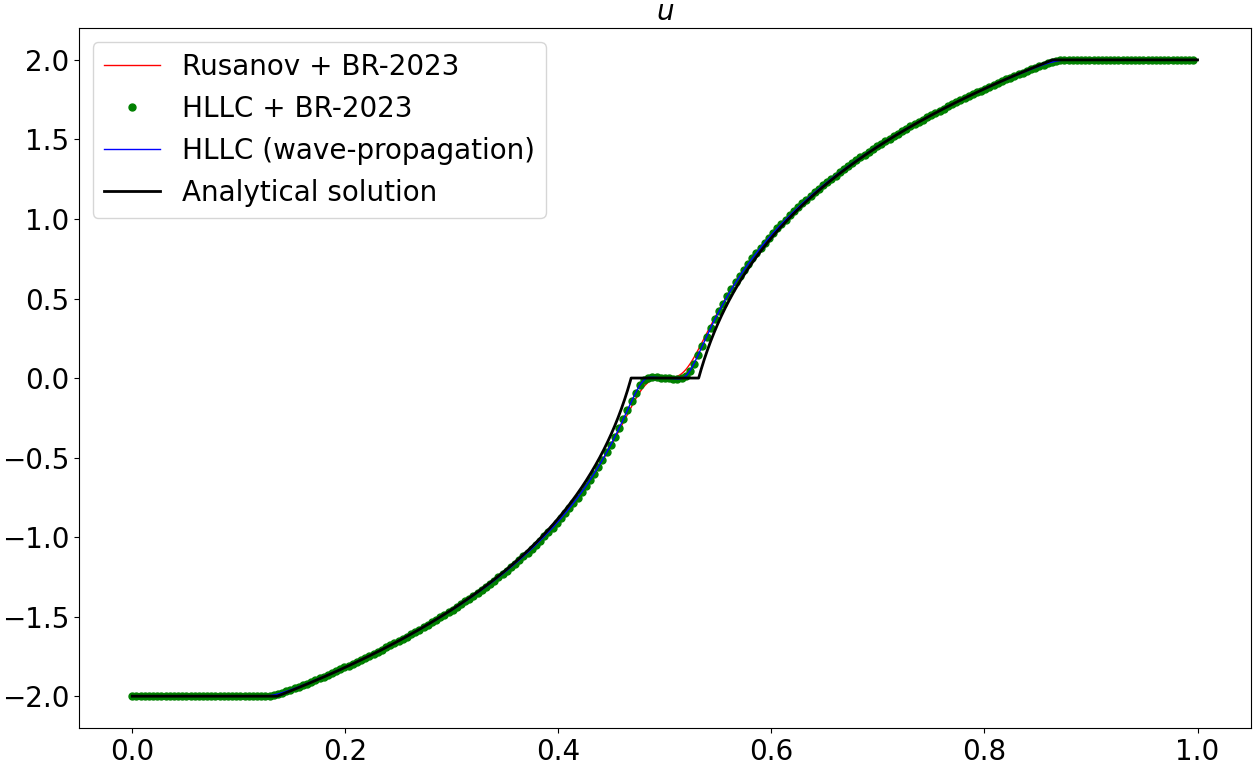}
	\end{subfigure}
	\caption{Double rarefaction test case with the fine mesh, results at $t = T_{f}$. Top-left: volume fraction of phase $2$. Top-right: mechanical equilibrium pressure. Bottom-left: mixture density. Bottom-right: velocity. Continuous black lines: analytical results. Continuous blue lines: HLLC-type wave-propagation scheme. Continuous red lines: Rusanov flux in combination with \textit{BR-2023} \eqref{eq:BR_Orlando_non_cons}. Green dots: HLLC flux in combination with \textit{BR-2023} \eqref{eq:BR_Orlando_non_cons}.}
	\label{fig:double_rarefaction}
\end{figure}

Next, we consider a water liquid-vapour cavitation (expansion) tube problem as presented in \cite{pelanti:2014}. The liquid, which we conventionally denote with phase $1$, contains a uniformly distributed small amount of vapour, i.e. phase $2$ (Table \ref{tab:phase2_init}). We focus on the instantaneous mechanical relaxation and we employ the fine mesh. The solution consists of two rarefactions propagating in opposite directions and producing a pressure decrease in the middle of the tube \cite{pelanti:2014}. An analytical solution is available for this test case. A good agreement is established between the different numerical strategies and with the analytical solution (Figure \ref{fig:double_rarefaction}). However, we notice visible differences in correspondence of the contact discontinuity (Figure \ref{fig:double_rarefaction}), in particular for what concerns the volume fraction and the mixture density. For the sake of completeness, since an analytical solution is available, we perform a convergence test. The expected convergence rates are achieved (Table \ref{tab:relative_errors_double_rarefaction}). The deviation of the observed order on the finest meshes simply indicates that mesh convergence has been reached. These results further confirm the robustness of the approaches employed in this work. Indeed, this test case is characterized by low values of the vapour volume fraction $\alpha_{2}$ and can be subject to the formation of the vacuum state \cite{toro:2009}.

\begin{table}[h!]
	\centering
	\footnotesize
	\begin{tabularx}{0.9\columnwidth}{lXXXXXX}
		\toprule
		\multirow{2}{*}{$N_{\text{cells}}$} & \multicolumn{2}{c}{Rusanov + \textit{BR-2023}} & \multicolumn{2}{c}{HLLC + \textit{BR-2023}} & \multicolumn{2}{c}{HLLC (wave-propagation)} \\
		\cmidrule(l){2-7}
		& $l^{1}$ rel. error $p$ & EOC & $l^{1}$ rel. error $p$ & EOC & $l^{1}$ rel. error $p$ & EOC \\
		\midrule
		$1024$ & \num{1.79e-1} & - & \num{1.28e-1} & - & \num{1.28e-1} & - \\
		\midrule
		$2048$ & \num{1.16e-1} & $0.6$ & \num{7.96e-2} & $0.7$ & \num{7.96e-2} & $0.7$ \\
		\midrule
		$4096$ & \num{7.09e-2} & $0.7$ & \num{4.66e-2} & $0.8$ & \num{4.66e-2} & $0.8$ \\
		\midrule
		$8192$ & \num{4.12e-2} & $0.8$ & \num{2.55e-2} & $0.9$ & \num{2.55e-2} & $0.9$ \\
		\midrule
		$16384$ & \num{2.27e-2} & $0.9$ & \num{1.35e-2} & $0.9$ & \num{1.35e-2} & $0.9$ \\
		\midrule
		$32768$ & \num{1.28e-2} & $0.8$ & \num{8.99e-3} & $0.6$ & \num{8.99e-3} & $0.6$ \\
		\midrule
		$65536$ & \num{9.32e-3} & $0.5$ & \num{8.23e-3} & $0.1$ & \num{8.23e-3} & $0.1$ \\
		\bottomrule
	\end{tabularx}
	\caption{Double rarefaction test case: experimental order of convergence (EOC) in the $l^{1}$ norm for the pressure at $t = T_{f}$.}
	\label{tab:relative_errors_double_rarefaction}
\end{table}

\subsection{Epoxy-spinel shock}
\label{ssec:epoxy_spinel_relax}

\begin{figure}[h!]
	\centering
	\begin{subfigure}{0.475\textwidth}
		\centering
		\includegraphics[width = 0.95\textwidth]{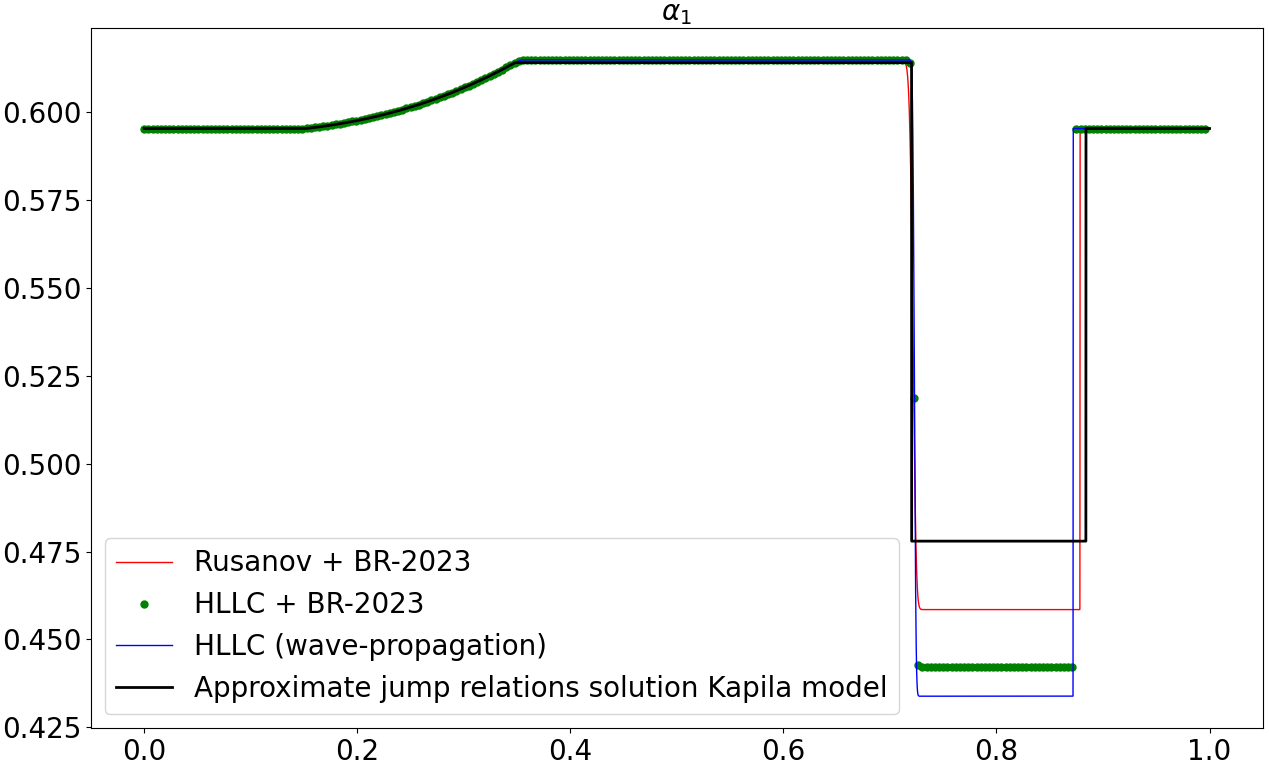}
	\end{subfigure}
	\begin{subfigure}{0.475\textwidth}
		\centering
		\includegraphics[width = 0.95\textwidth]{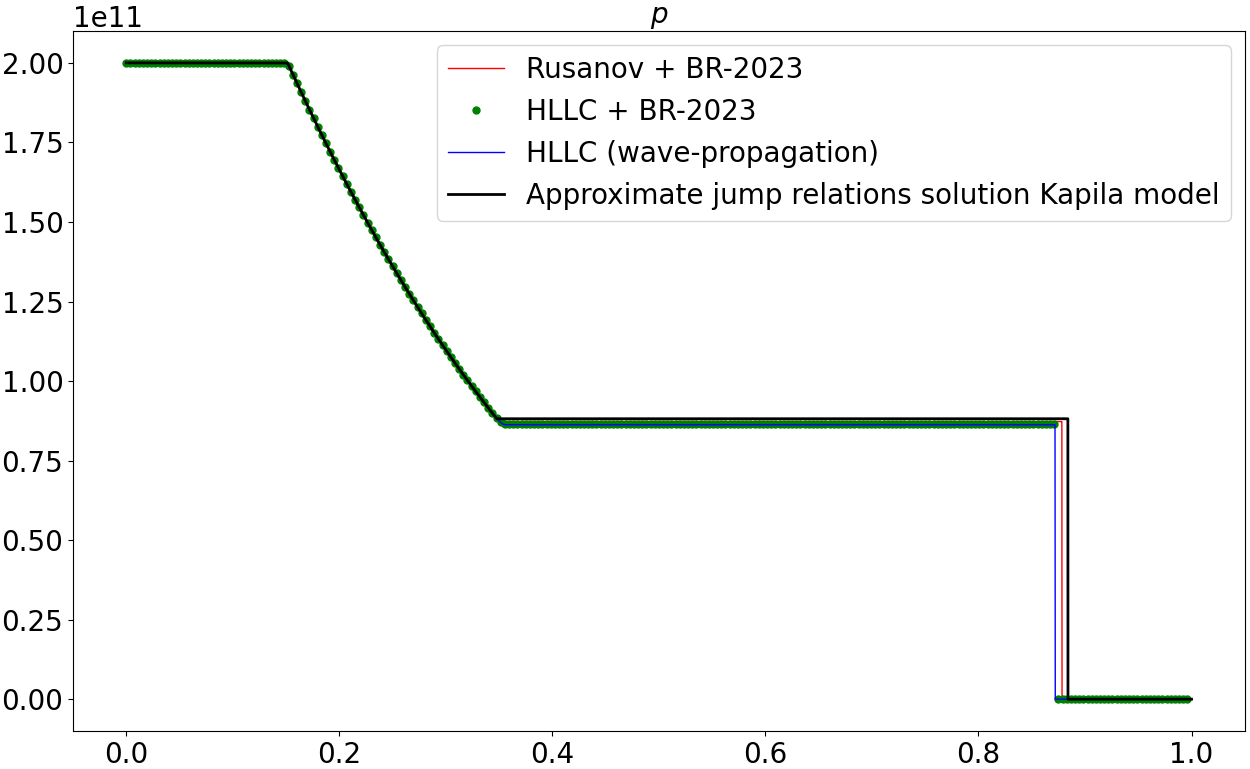}
	\end{subfigure}
	\begin{subfigure}{0.475\textwidth}
		\centering
		\includegraphics[width = 0.95\textwidth]{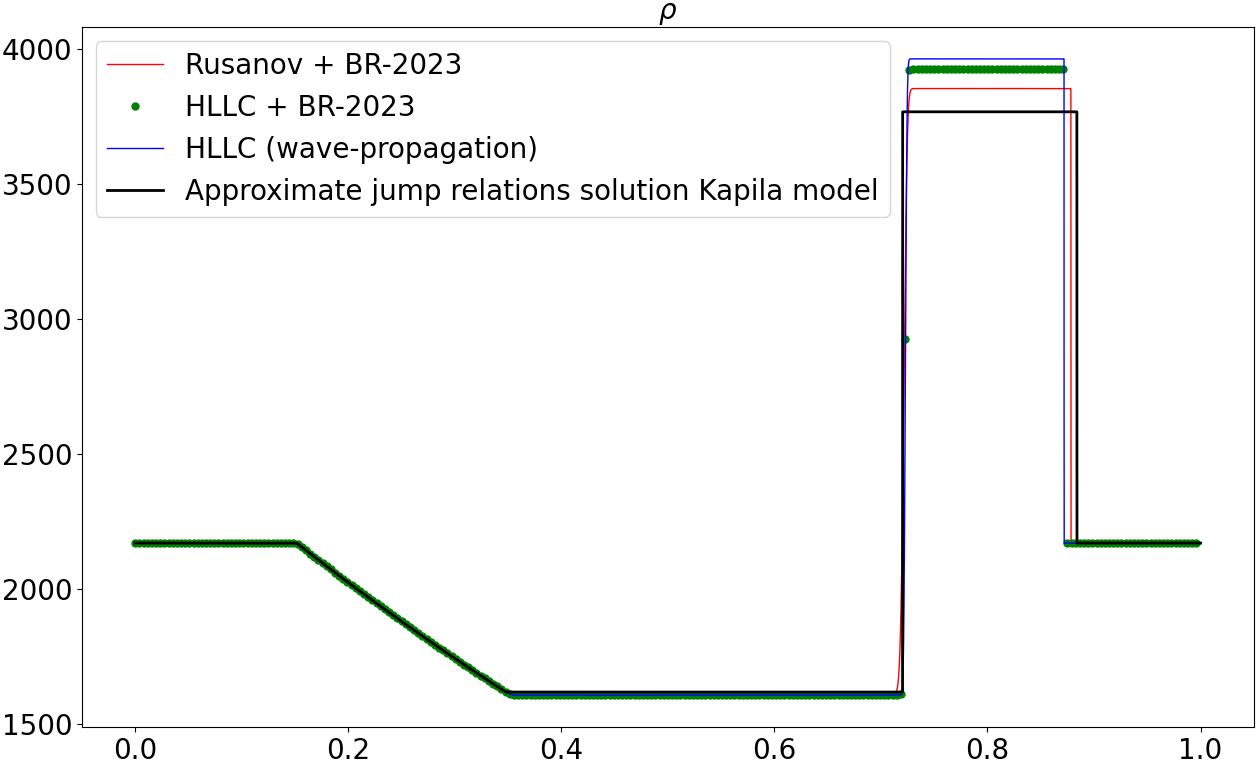}
	\end{subfigure}
	\begin{subfigure}{0.475\textwidth}
		\centering
		\includegraphics[width = 0.95\textwidth]{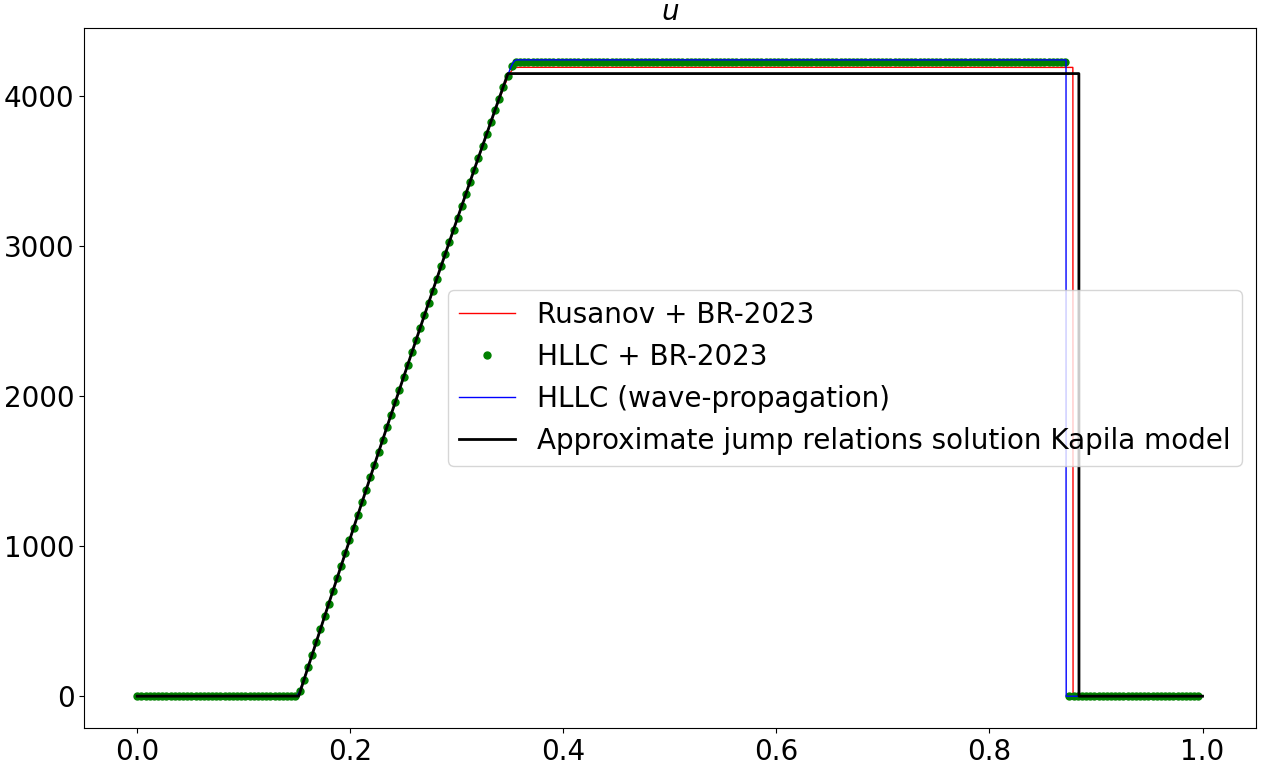}
	\end{subfigure}
	\caption{Epoxy-spinel shock test case with the fine mesh, results at $t = T_{f}$. Top-left: volume fraction of phase $1$. Top-right: mechanical equilibrium pressure. Bottom-left: mixture density. Bottom-right: velocity. Continuous black lines: solutions obtained using the approximated jump relations in \cite{petitpas:2007}. Continuous blue lines: HLLC-type wave-propagation scheme. Continuous red lines: Rusanov flux in combination with \textit{BR-2023} \eqref{eq:BR_Orlando_non_cons} Green dots: HLLC flux in combination with \textit{BR-2023} \eqref{eq:BR_Orlando_non_cons}.}
	\label{fig:epoxy_spinel_shock_relax}
\end{figure}

Next, we consider the epoxy-spinel strong shock already introduced in Section \ref{ssec:epoxy_spinel_homogeneous}. We employ the fine mesh, so as to achieve mesh convergence. One can easily verify that the solution does not converge to the solution obtained using the approximated jump conditions in \cite{petitpas:2007} (Figure \ref{fig:epoxy_spinel_shock_relax}). As discussed in \cite{petitpas:2007}, the approximated jump relations are valid only in the case of a weak shock. This discrepancy depends on the partition of internal energy in the shock layer (see the discussions in \cite{petitpas:2007, saurel:2009}). The convergence value of the intermediate state of the volume fraction for the Riemann problem is strongly affected by the choice of the numerical flux for the conservative part and visible differences arise also between the HLLC-type wave-propagation scheme and the combination with the BR-type approach (Figure \ref{fig:epoxy_spinel_shock_relax}). Note that this result is totally dependent on the interaction between the numerical method employed to solve the homogeneous problem \eqref{eq:system_hyperbolic} and the mechanical relaxation. A uniform volume fraction is indeed imposed as initial condition (see Table \ref{tab:phase1_init}) and therefore, without relaxation, the initial value is preserved. Moreover, we recall that in this test case we are dealing with a mixture region rather than quasi-pure phases and therefore non-conservative terms may have a greater impact.

\begin{figure}[h!]
	\centering
	\includegraphics[width = 0.7\textwidth]{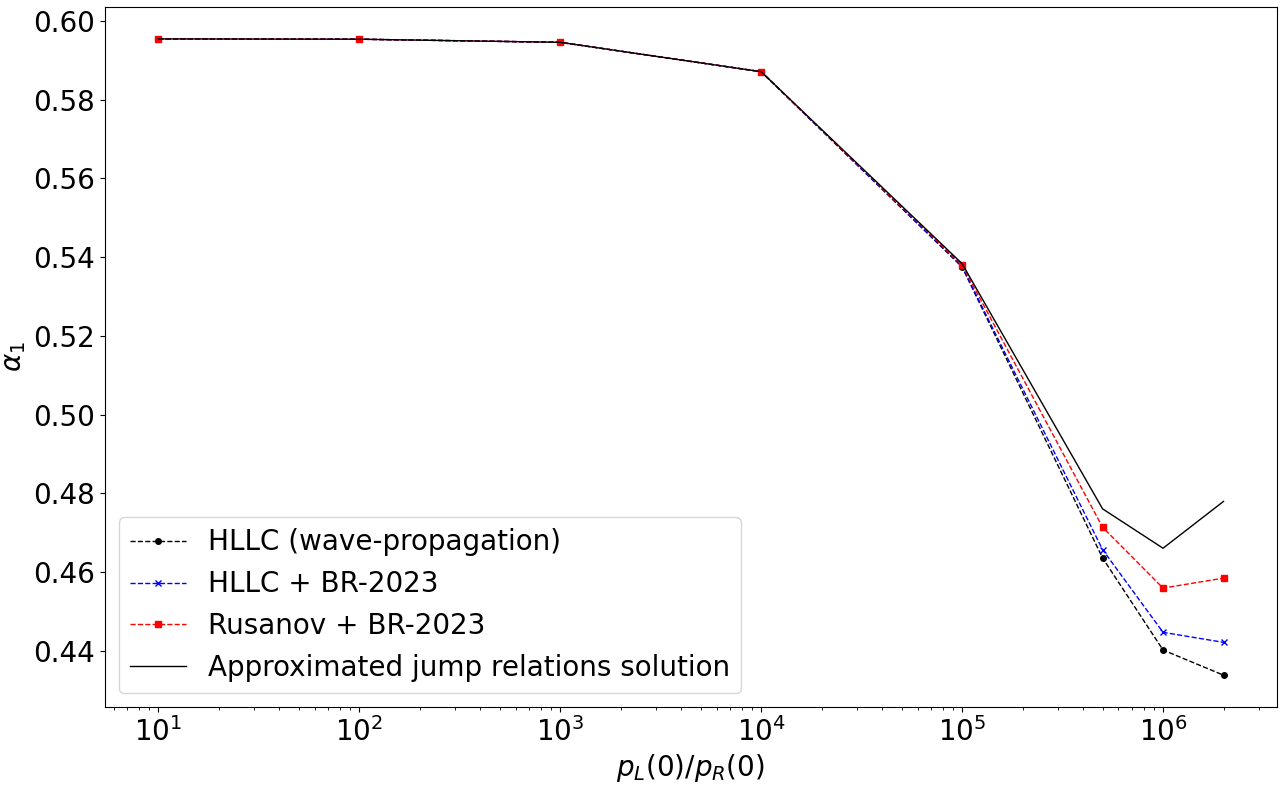}
	\caption{Epoxy-spinel shock test case with the fine mesh, intermediate state of the volume fraction ($y$-coordinate) as a function of the initial pressure ratio ($x$-coordinate), results at $t = T_{f}$. The continuous black line shows the value predicted by the approximated jump relations in \cite{petitpas:2007}. The black dots report the results obtained using the HLLC-type wave-propagation scheme, the blue crosses show the results obtained using the HLLC flux in combination with \textit{BR-2023} \eqref{eq:BR_Orlando_non_cons} for the non-conservative terms of the energy equations, whereas the red square represent the results obtained using the Rusanov flux in combination with \textit{BR-2023} \eqref{eq:BR_Orlando_non_cons} for non-conservative terms.}
	\label{fig:epoxy_spinel_shock_relax_alpha1}
\end{figure}

In order to further analyze this behaviour, we consider different initial pressure ratios, i.e. we set for the initial conditions
$$p_{1,L} = p_{2,L} = \cpth{10^{6}, 10^{7}, 10^{8}, 10^{9}, 10^{10}, 5 \times 10^{10}, 10^{11}}.$$
One can notice that visible difference arise starting from $p_{1,L}/p_{1,R} = p_{2,L}/p_{2,R} = 5 \times 10^{5}$ (Figure \ref{fig:epoxy_spinel_shock_relax_alpha1}). The value predicted by the approximated jump relations in \cite{petitpas:2007} is added for reference. This result reflects the absence of a complete set of well-defined Rankine--Hugoniot conditions for the 5-equation model proposed in \cite{kapila:2001}. The main result for this test case is that the differences observed in the volume fraction field are quite moderate when compared to the extreme pressure ratio needed to generate them. It thus shows that the theoretical shortcomings of the model, i.e. non-uniqueness of shock profiles, has only a limited impact in practical applications. The use of robust numerical strategies allows to obtain physically relevant numerical solutions.

\subsection{2D Riemann problem}
\label{ssec:2D_Riemann}

In a final test, we consider a two-dimensional Riemann problem inspired from the configuration C1 employed in \cite{dumbser:2013} and analyzed also in \cite{dumbser:2016, fraysse:2016}. In \cite{dumbser:2013}, the test is solved for the full Baer--Nunziato model, with relaxation parameters chosen so as to mimic the Kapila model. The computational domain is $\Omega = \rpth{-0.5, 0.5} \times \rpth{0.5, 0.5}$. The EOS parameters are $\gamma_{1} = 1.4, \gamma_{2} = 1.67, \pi_{1} = \pi_{2} = 0, \eta_{1} = \eta_{2} = 0$. The initial conditions are reported in Table \ref{tab:2D_Riemann_init}. The final time is $T_{f} = \SI{0.15}{\second}$. Following \cite{dumbser:2013}, we adopt reflective boundary conditions. In order to enhance the computational efficiency, the multiresolution capabilities of \texttt{samurai}\footnote{\url{https://github.com/hpc-maths/samurai}} \cite{bellotti:2022} are employed. More specifically, we consider a minimum resolution which would correspond to a uniform mesh with $N_{\text{cells}} = 256$ per direction and maximum resolution which would correspond to a uniform mesh with $N_{\text{cells}} = 2048$ per direction, i.e. $2048^{2} \approx 4 \cdot 10^{6}$ cells. The Courant number is set to $C = 0.45$.

\begin{table}[h!]
	\centering
	\footnotesize
	\begin{tabularx}{0.5\columnwidth}{lrrrrrrr}
		\toprule
		& $\alpha_{1}$ & $\rho_{1}$ & $p_{1}$ & $u$ & $v$ & $\rho_{2}$ & $p_{2}$ \\
		\midrule
		$\mathbf{Q}_{1}: \rpth{x > 0, y > 0}$ & $0.8$ & $2.0$ & $2.0$ & $0.0$ & $0.0$ & $1.5$ & $2.0$ \\
		$\mathbf{Q}_{2}: \rpth{x < 0, y > 0}$ & $0.4$ & $1.0$ & $1.0$ & $0.0$ & $0.0$ & $0.5$ & $1.0$ \\
		$\mathbf{Q}_{3}: \rpth{x < 0, y < 0}$ & $0.8$ & $2.0$ & $2.0$ & $0.0$ & $0.0$ & $1.5$ & $2.0$ \\
		$\mathbf{Q}_{4}: \rpth{x > 0, y < 0}$ & $0.4$ & $1.0$ & $1.0$ & $0.0$ & $0.0$ & $0.5$ & $1.0$ \\
		\bottomrule
	\end{tabularx}
	\caption{Initial conditions for the 2D Riemann problem.}
	\label{tab:2D_Riemann_init}
\end{table}

\begin{figure}[h!]
	\centering
	\begin{subfigure}{0.475\textwidth}
		\centering
		\includegraphics[width = 0.95\textwidth]{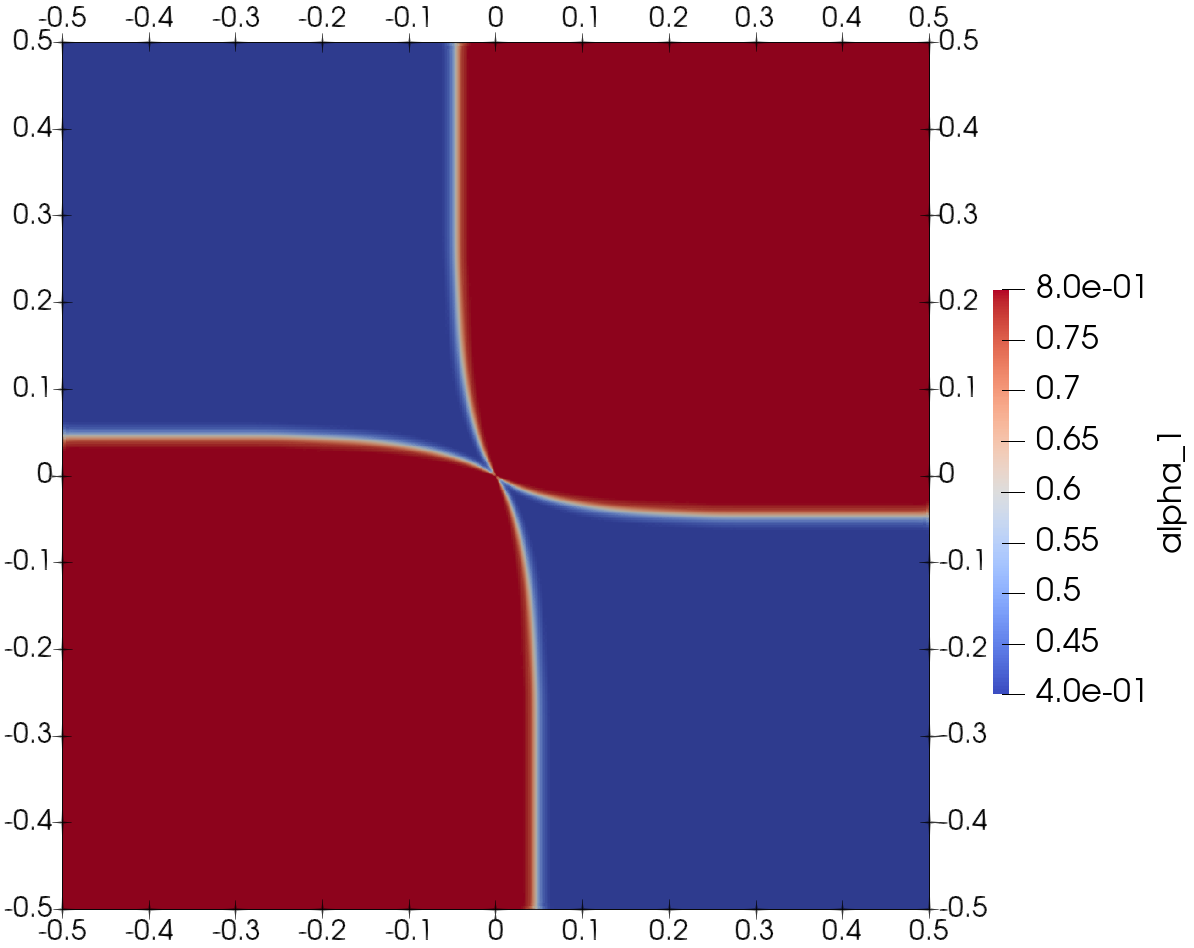}
	\end{subfigure}
	\begin{subfigure}{0.475\textwidth}
		\centering
		\includegraphics[width = 0.95\textwidth]{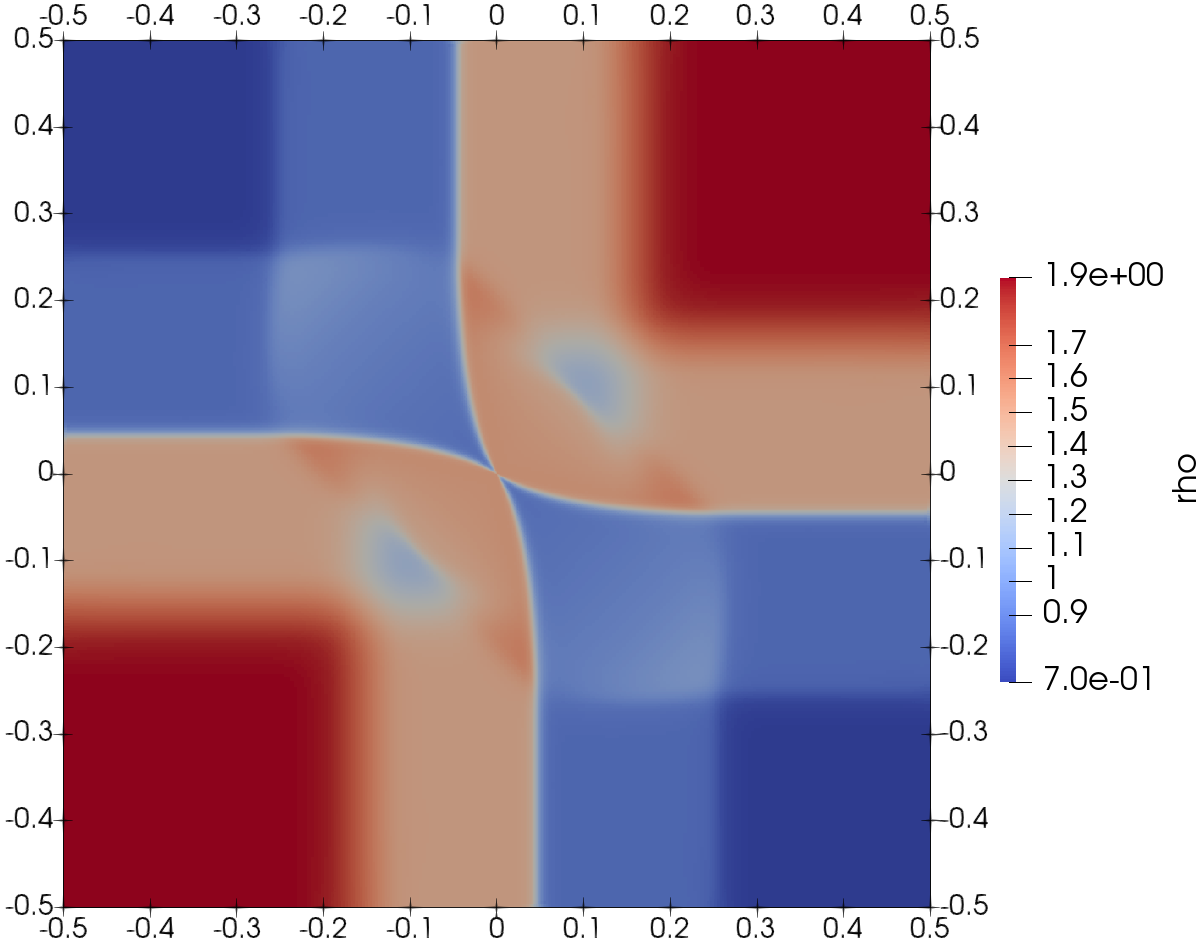}
	\end{subfigure}
	\begin{subfigure}{0.475\textwidth}
		\centering
		\includegraphics[width = 0.95\textwidth]{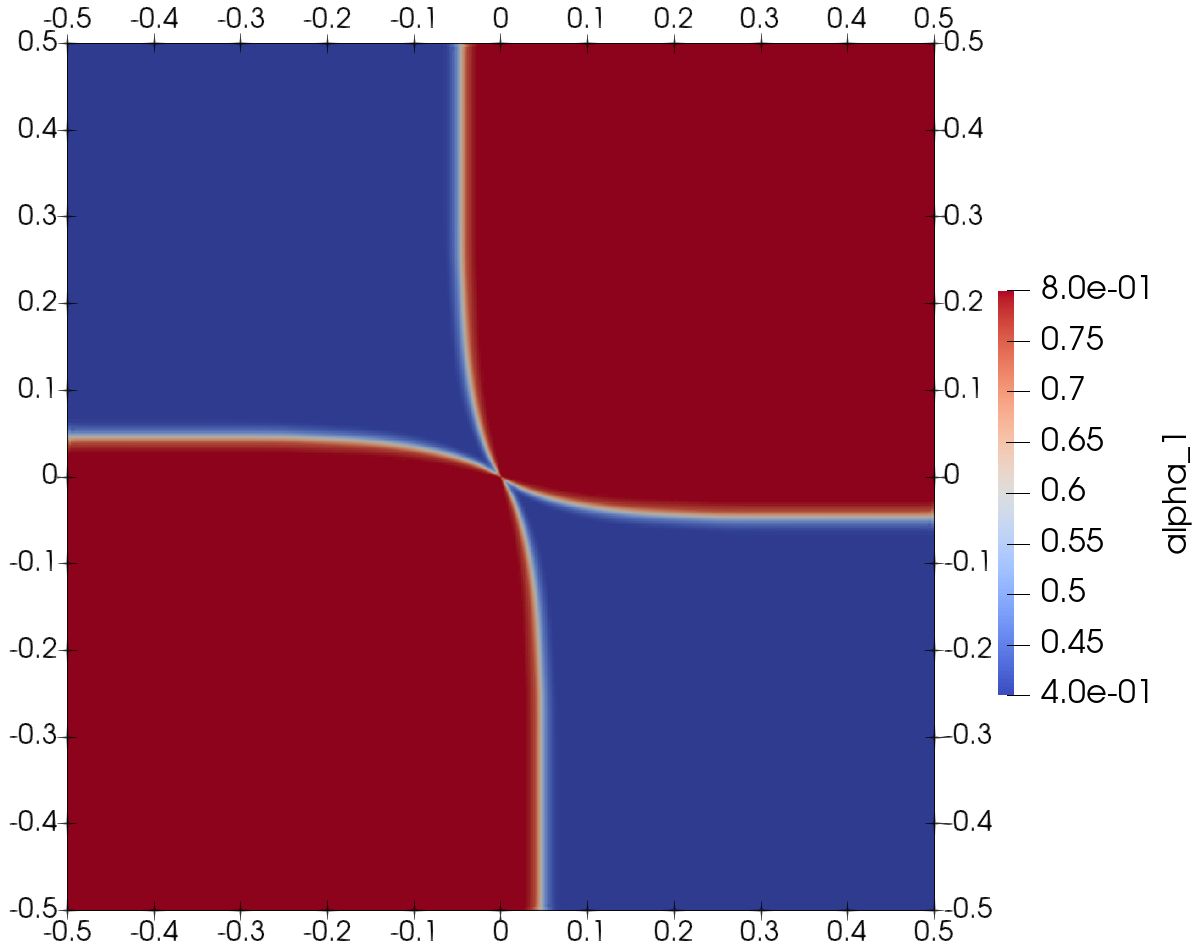}
	\end{subfigure}
	\begin{subfigure}{0.475\textwidth}
		\centering
		\includegraphics[width = 0.95\textwidth]{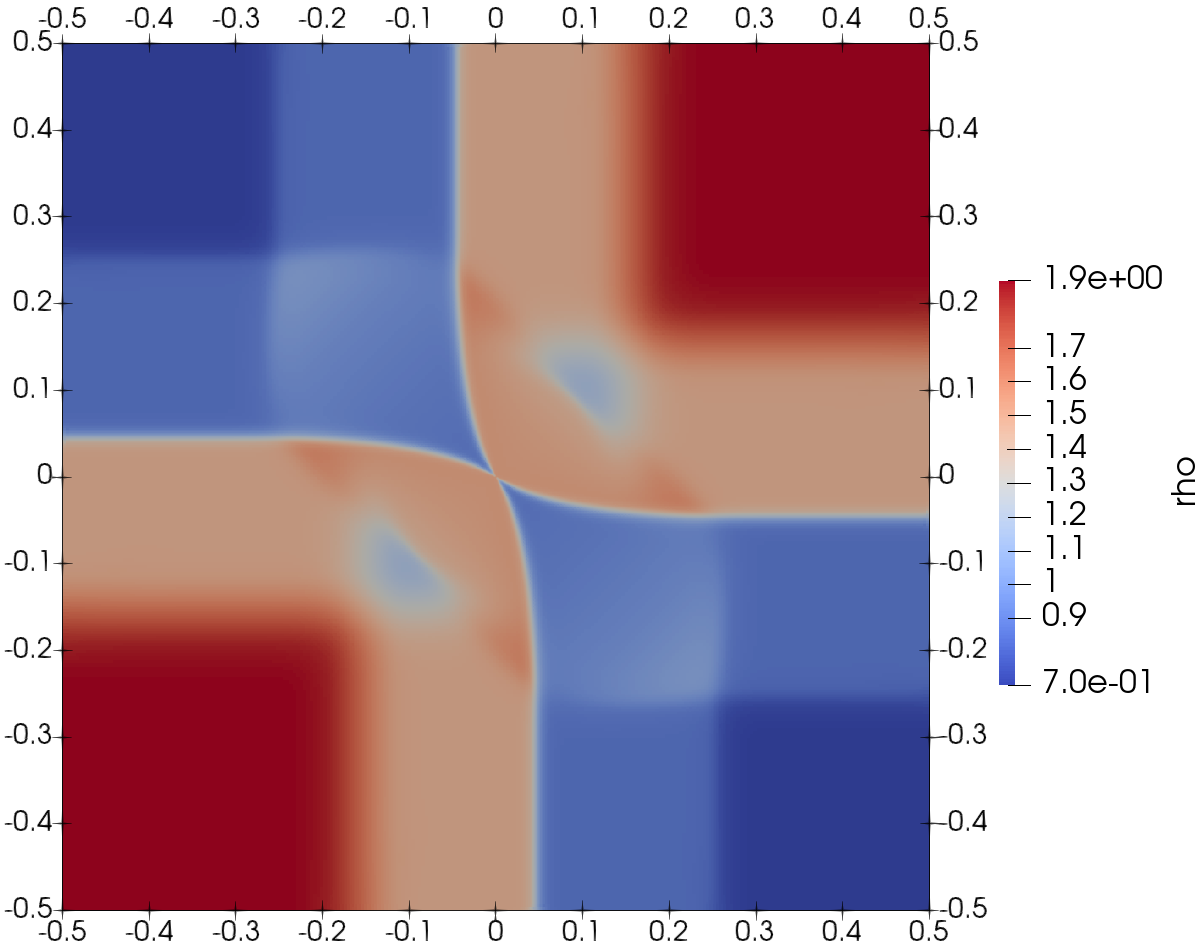}
	\end{subfigure}
	\begin{subfigure}{0.475\textwidth}
		\centering
		\includegraphics[width = 0.95\textwidth]{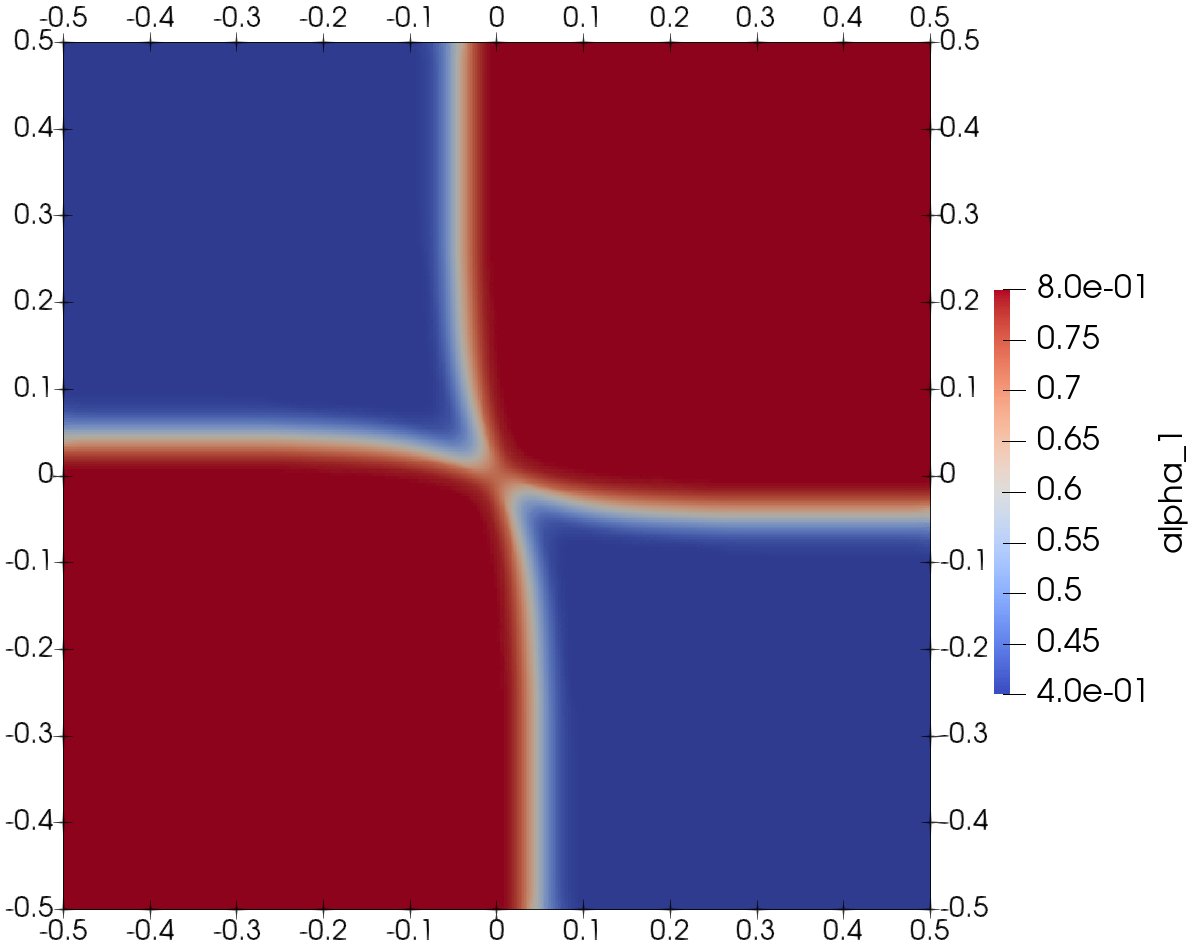}
	\end{subfigure}
	\begin{subfigure}{0.475\textwidth}
		\centering
		\includegraphics[width = 0.95\textwidth]{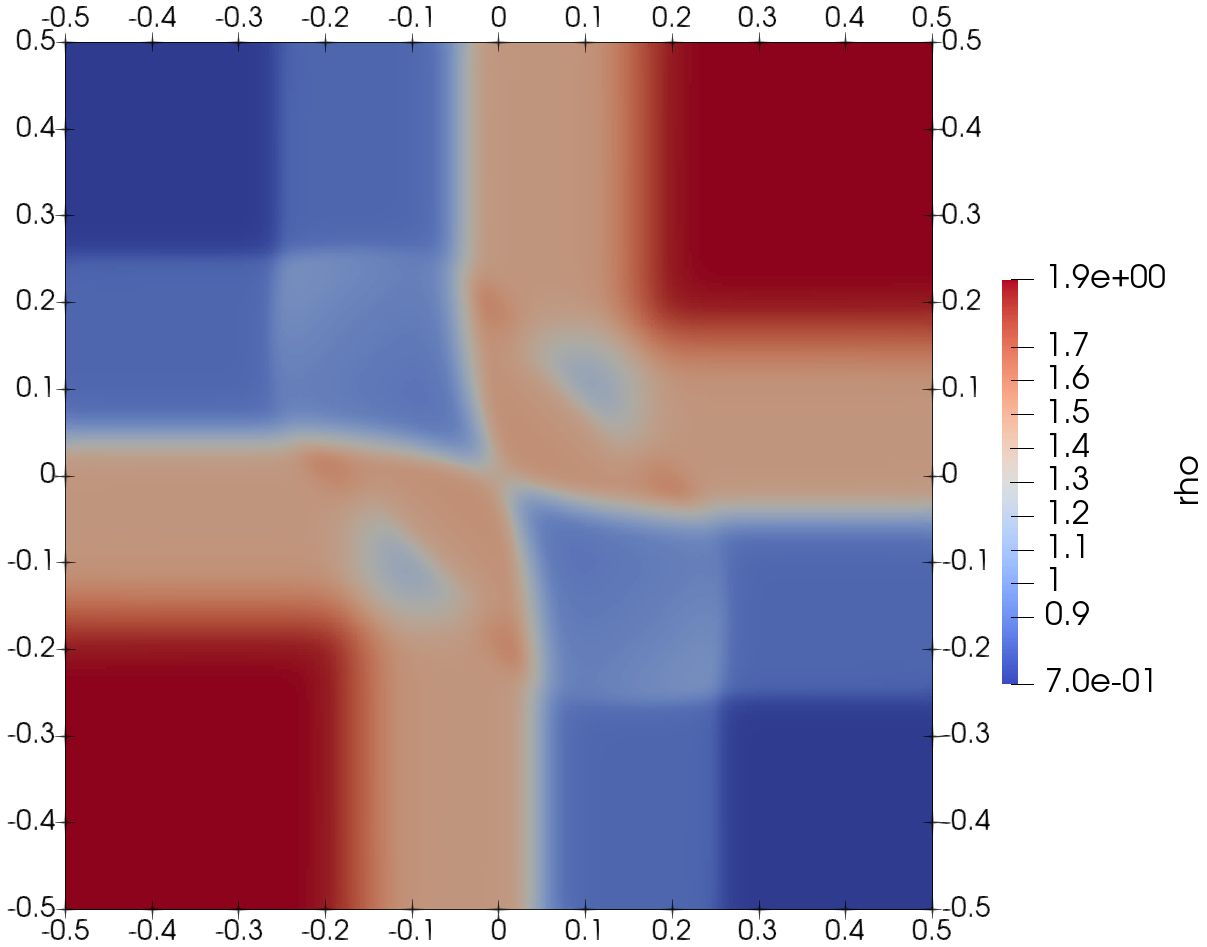}
	\end{subfigure}
	\caption{2D Riemann problem test case, results at $t = T_{f}$. Left column: volume fraction of phase $1$. Right column: mixture density. Top: HLLC-type wave-propagation. Middle: HLLC flux in combination with \textit{BR-2023} \eqref{eq:BR_Orlando_non_cons}. Bottom: Rusanov in combination with \textit{BR-2023} \eqref{eq:BR_Orlando_non_cons}.}
	\label{fig:2D_Riemann_Tf}
\end{figure}

\begin{figure}[h!]
	\centering
	\begin{subfigure}{0.475\textwidth}
		\centering
		\includegraphics[width = 0.95\textwidth]{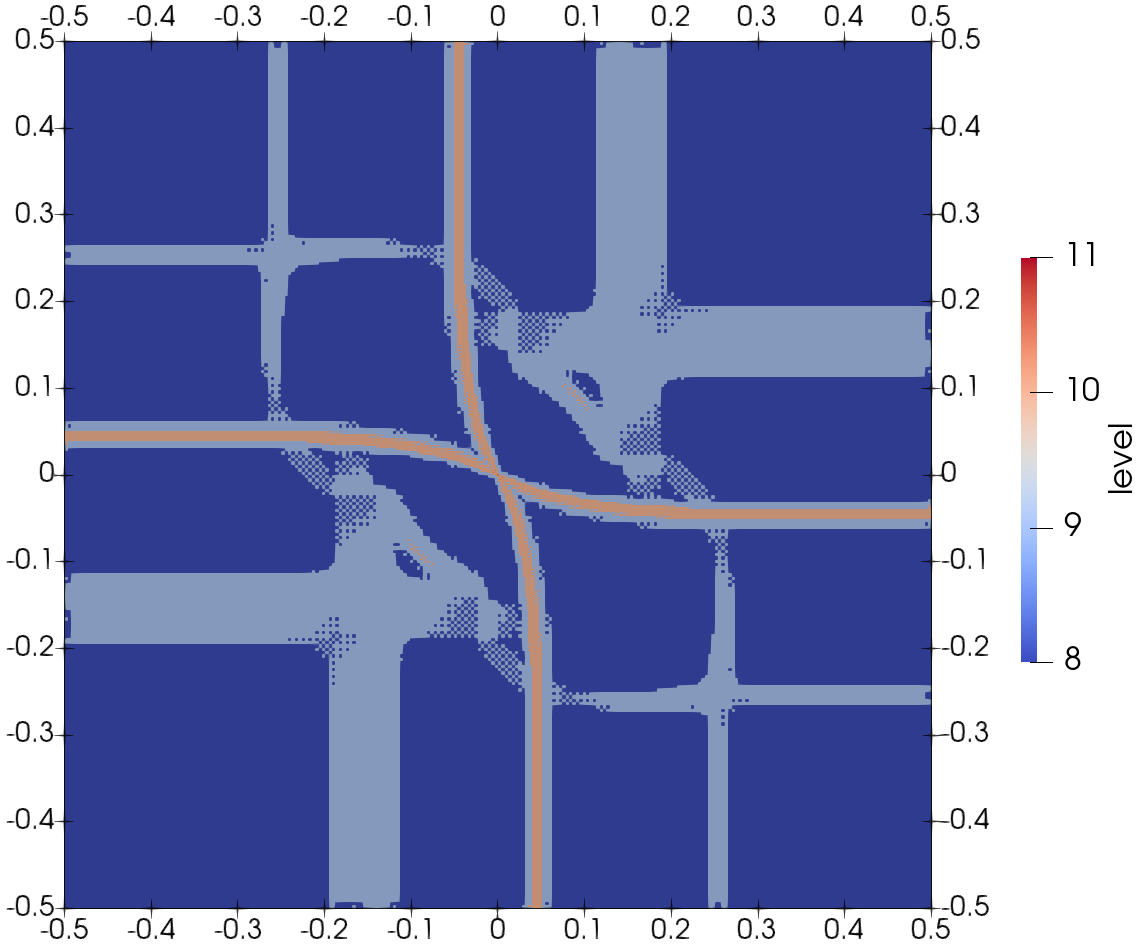}
	\end{subfigure}
	\begin{subfigure}{0.475\textwidth}
		\centering
		\includegraphics[width = 0.95\textwidth]{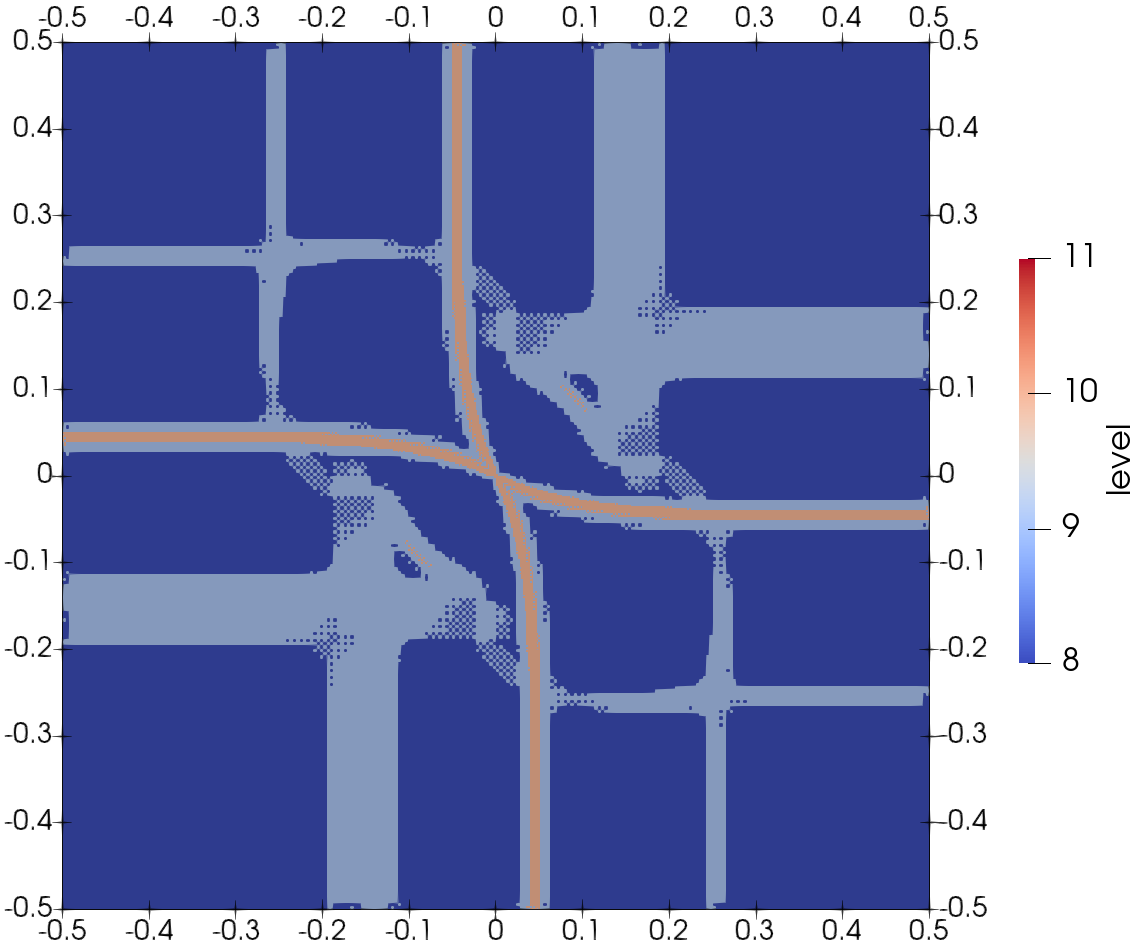}
	\end{subfigure}
	\begin{subfigure}{0.475\textwidth}
		\centering
		\includegraphics[width = 0.95\textwidth]{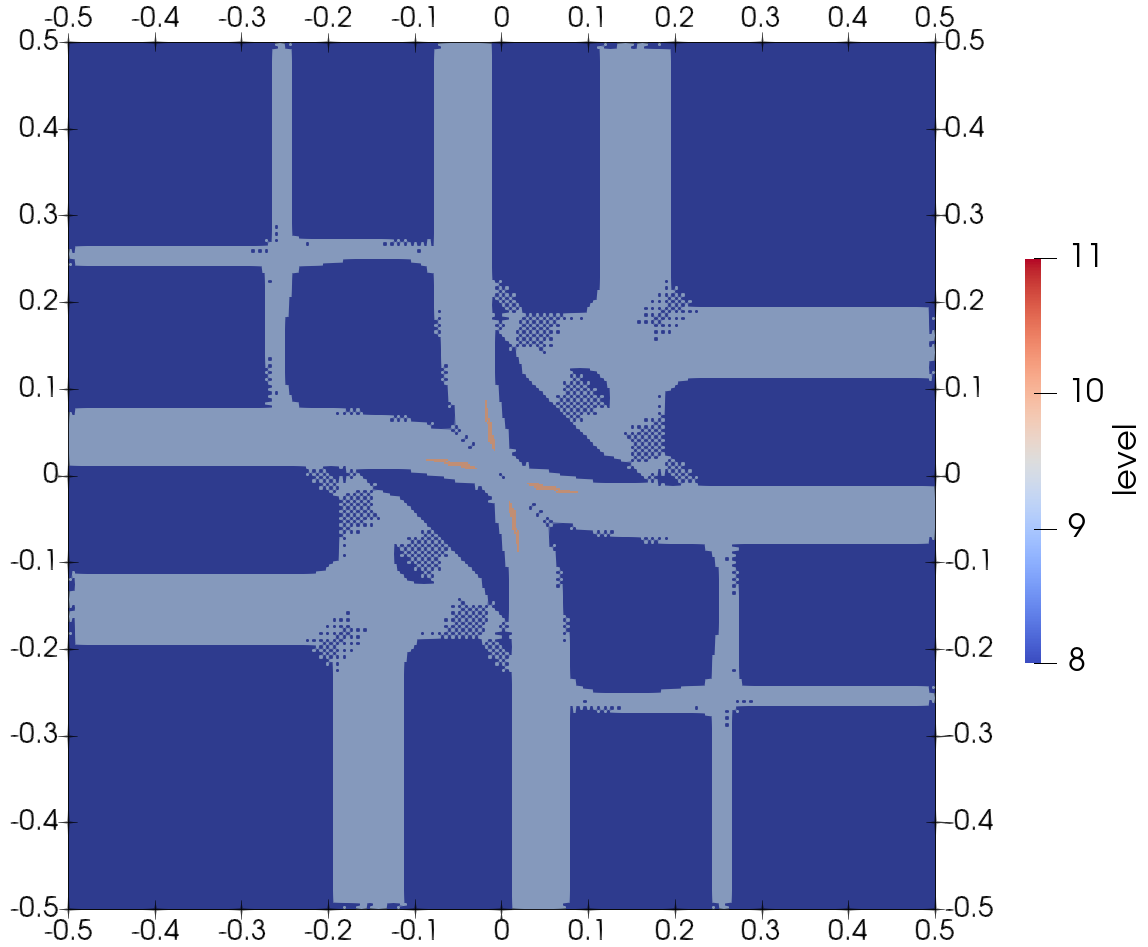}
	\end{subfigure}
	\caption{2D Riemann problem test case, computational mesh at different refinement levels at $t = T_{f}$. Top-left: HLLC-type wave-propagation scheme ($N_{\text{cells}} = 133 036$). Top-right: HLLC flux in combination with \textit{BR-2023} \eqref{eq:BR_Orlando_non_cons} ($N_{\text{cells}} = 134 530$). Bottom: Rusanov flux in combination with \textit{BR-2023} \eqref{eq:BR_Orlando_non_cons} ($N_{\text{cells}} = 126 556$). Level 8 corresponds to the minimum resolution, equivalent to a uniform mesh with $N_{\text{cells}} = 256$ per direction, whereas level 11 corresponds to the maximum resolution, equivalent to a uniform mesh with $N_{\text{cells}} = 2048$ per direction.}
	\label{fig:2D_Riemann_mesh_Tf}
\end{figure}

A good agreement is established between all the numerical methods (Figure \ref{fig:2D_Riemann_Tf}) and with the results presented in \cite{dumbser:2013}. This confirms the robustness of the computational framework presented in this work and it also shows its flexibility. Finally, we point out that the final computational meshes of the HLLC-based schemes are characterized by a slightly higher number of cells than that obtained using the Rusanov scheme (Figure \ref{fig:2D_Riemann_mesh_Tf}). To avoid any misconception about this result, we recall that the basic principle of the multiresolution is that of data compression \cite{bellotti:2022, harten:1995}. The flow fields are decomposed on a local wavelet basis, which allows to quantify their local regularity and identify regions in which the data can be compressed without loss of information --- up to a pre-defined tolerance. Hence, since the Rusanov scheme is more diffusive, it captures less flow features, thus allowing for a higher-compression rate of the mesh. In particular, the more smeared interface that we observe with the Rusanov scheme (bottom left plot of Figure \ref{fig:2D_Riemann_Tf}) is a consequence of its higher numerical diffusion and not the coarser mesh.

\section{Conclusions}
\label{sec:conclu}

We have presented a robust computational framework for the 5-equation model of Kapila \cite{kapila:2001} by means of the mixture-energy-consistent 6-equation two-phase flow model with instantaneous mechanical relaxation terms. We have shown that, despite the absence of uniquely defined shock profiles for the considered models, correct numerical solutions can be obtained for practical cases through the use of robust and accurate numerical methods. To achieve this result, we have compared different numerical methods for the discretization of the homogeneous system \eqref{eq:system_hyperbolic}, which rely on different treatments of the non-conservative terms. One of these treatments is based on a Bassi-Rebay (BR)-type approach \cite{bassi:1997a}. We have shown here for the first time that this strategy, typically employed for diffusion-type operators, can be reformulated in the framework of the path-conservative schemes, which are widely employed for the discretization of non-conservative hyperbolic systems.

First, we have compared the results obtained for the homogeneous model, using the HLLC approximate Riemann solvers with the wave-propagation method, and the Rusanov and HLLC approximate Riemann solvers in combination with the BR-type approach or the approximation proposed in \cite{crouzet:2013} for the discretization of non-conservative terms. This stage is particularly relevant and represents one of the contributions of the work. Recent works \cite{pelanti:2022, schmidmayer:2023} have equipped the 6-equation model with finite-rate mechanical relaxation, meaning that the model is no longer merely an auxiliary system for the numerical solution of the Kapila model, but becomes a standalone model in its own right. Hence, understanding its discretization is crucial. Since the model does not admit a complete set of Rankine--Hugoniot conditions uniquely defined, the construction of Riemann solvers is not straightforward. We have presented in this work some insights on the construction of approximate Riemann solvers so as to clarify the adopted underlying hypotheses. The \Crouzet\ approach \cite{crouzet:2013} does not guarantee in general an admissible solution.

Once mesh convergence is achieved, the remaining strategies do not converge in general to the same numerical solution for the homogeneous model. Moreover, even though kinetic equilibrium is prescribed, the results in general slightly differ from those obtained if one employs the full Baer-Nunziato (BN) model, as done in \cite{tokareva:2010}. In particular, spurious oscillations may arise around the contact discontinuity for the phasic pressures, see Section \ref{ssec:sonic_rarefaction}. This is a common issue for all numerical schemes, especially for HLLC-type schemes applied to the conservative part of the system and the volume fraction, combined with the discretization of the non-conservative terms in the energy equations. However, once mesh convergence is achieved, the HLLC wave-propagation scheme mitigates this issue most effectively.

Finally, we have shown in the test case of a water-air shock tube that different numerical methods converge towards different solutions. This result raises the question of the use of the 6-equation model without mechanical relaxation or when finite-rate relaxation occurs, as, e.g., in \cite{schmidmayer:2023}. Moreover, for this test case, Rusanov-based schemes need a lower Courant number to achieve a physically meaningful solution, whereas standard restrictions occur for the HLLC-based schemes. Finally, we have performed a sensitivity analysis with respect to the ``residual'' volume fraction in order to more thoroughly investigate the behaviour of the numerical schemes in the nearly pure-fluid regime. The solution obtained using the wave-propagation scheme approaches that of the Euler--Euler equations as the volume fraction tends to zero, whereas the results obtained by employing a discretization for the conservative part and a separate, simple discretization for the non-conservative terms exhibit almost no dependence on the value of the ``residual'' volume fraction. This behaviour is likely a further manifestation of the absence of a uniquely defined set of jump conditions. Overall, HLLC-based schemes, and in particular the HLLC wave-propagation scheme, are generally more robust for the numerical solution of the 6-equation model.

Next, we have focused on the case where instantaneous mechanical relaxation towards equilibrium is accounted for. In this case, for most configurations, all the methods tend in general towards the same approximate solution. However, when extreme pressure ratios are considered, different numerical schemes may tend towards different numerical solutions. A similar dependence of the numerical solution to the numerical parameters has also been observed in \cite{haegeman:2024} for this case, as different interfacial pressures lead to different shock profiles at extreme pressure ratios in the epoxy-spinel test case. This result depends on the absence of a well-defined set of Rankine--Hugoniot conditions. Nevertheless, the observed differences are moderate in comparison with the pressure ratios required to generate them, and therefore, for most applications, this shortcoming of the model might not have a huge impact, even though it should be taken into account.

In this work, we have achieved mesh convergence with first order schemes, so as to avoid any further dependence of higher order methods related to reconstruction, limiters, and large stencils. The primary goal of this work is indeed to perform a detailed comparison between different discretization strategies for the mixture-energy-consistent 6-equation two-phase model in terms of accuracy and robustness and to analyze their non-trivial interaction with the instantaneous mechanical relaxation. The numerical results have been obtained using a number of cells that is unfeasible for practical applications with the ultimate goal to highlight the possible limitations of the model when mesh convergence is really established. In future work, we aim at analyzing the high-order extension of some of the presented discretization methods. While the results of the present work suggest that the wave-propagation scheme performs better in all respects, its extension to high-order accuracy within a finite volume framework and using the classical method-of-lines approach is less straightforward than for more standard Godunov-type methods with a separate discretization of conservative and non-conservative terms, since volumetric contributions also come into play \cite{ketcheson:2013}. Alternative approaches include space-time methods, such as the Lax--Wendroff procedure \cite{pelanti:2014}, which, however, are relatively straightforward to develop only up to second-order accuracy. The outcomes of this work also show, to the best of our knowledge for the first time, that, although good results can be obtained when instantaneous mechanical relaxation is applied, without instantaneous mechanical relaxation, an unphysical dependence to the discretization method is observed. As such, for cases in which pressure non-equilibrium effects are relevant, the 6-equation model may not have the sufficient mathematical properties for practical use. For this reason, we aim at developing a numerical strategy for fully out-of-equilibrium flows, starting from the 7-equation model of Baer-Nunziato \cite{baer:1986}, or a recently derived all-topology two-fluid model \cite{haegeman:2026}, with source terms. These models are well posed, in the sense that non-conservative products are uniquely defined and they admit a mixture entropy inequality to select physically relevant weak solutions \cite{gallouet:2004,haegeman:2026}. This approach can lead to a more rigorous mathematical treatment, even though certainly increasing the computational burden. Finally, in the framework of the BN model, coupled finite-rate kinetic, mechanical, and thermal relaxation techniques have been recently developed in \cite{herard:2023}, however, the correct scaling of the different characteristic time scales involved is still an open problem.

\section*{Acknowledgements}

We gratefully acknowledge C. Le Touze for several useful discussions on related topics. This work has been supported by the CIEDS project OPEN NUM DEF (PI L. Gouarin, M. Massot, T. Pichard) as well as the HPC@Maths Initiative (PI L. Gouarin and M. Massot) of the Fondation {\'E}cole polytechnique. W.H. also acknowledges the Agence Innovation D{\'e}fense and ONERA for its support through a PhD grant.

\appendix

\section{Comments on the BR-type approach}
\label{app:BR_non_cons}

In this appendix, we further clarify, for the reader’s convenience, what we refer to when discussing the BR-type approach introduced in Section \ref{ssec:non_conservative}. One of the key contributions of \cite{bassi:1997a} (see also \cite{tumolo:2015}) is that the discrete gradient operator can be recast as the sum of two contributions: the first accounts for the elementwise gradient, while the second accounts for its jumps across element interfaces. Consider, for the sake of simplicity, the scalar case of the non-conservative product $\int_{K}B(w)\grad w\mathrm{d}\Omega$. Assume for the moment that $B(w)$ is continuous. Hence, provided that a numerical flux at the interface $\widehat{w} \equiv \widehat{w}(w_{L}, w_{R})$ is given, the discrete gradient operator $D_{h}w$ is defined by
\begin{equation}\label{eq:BR_grad}
	\int_{K}B(w)D_{h}w\mathrm{d}\Omega = \int_{K}B(w)\grad w\mathrm{d}\Omega + \int_{\partial K}B(w)\left(\widehat{w} - w\right)\bm{n}\mathrm{d}\Gamma.
\end{equation}
The above relation is obtained by first integrating by parts $\int_{K}B(w)\grad w\mathrm{d}\Omega$ which yields
\begin{equation}\label{eq:BR_grad_1}
	\int_{K}B(w)\grad w\mathrm{d}\Omega = -\int_{K}w\grad B(w)\mathrm{d}\Omega + \int_{\partial K}B(w)\widehat{w}\bm{n}\mathrm{d}\Gamma.
\end{equation}
Then, we back-integrate by parts the volumetric term $\int_{K}w\grad B(w)\mathrm{d}\Omega$ so as to obtain \eqref{eq:BR_grad}. By taking $B(w) \equiv 1$, one obtains
\begin{equation}
	\int_{K}D_{h}w\mathrm{d}\Omega = \int_{K}\grad w\mathrm{d}\Omega + \int_{\partial K}\left(\widehat{w} - w\right)\bm{n}\mathrm{d}\Gamma
\end{equation}
and, by choosing $\widehat{w} = \ave{w}$, one obtains
$$\int_{K}D_{h}w\mathrm{d}\Omega = \int_{K}\grad w\mathrm{d}\Omega + \int_{\partial K}\frac{1}{2}\left(w_{R} - w_{L}\right)\bm{n}\mathrm{d}\Gamma,$$
which is exactly the analogue of Equation (18) and the subsequent unnumbered relation in \cite{bassi:1997a}. 

The challenge when dealing with a generic non-conservative product is that $B(w)$ is, in general, not continuous. Hence, a choice must be made to approximate the boundary contribution of $B(w)$
\begin{equation}\label{eq:BR_non_cons}
	\int_{K}B(w)D_{h}w\mathrm{d}\Omega = \int_{K}B(w)\grad w\mathrm{d}\Omega + \int_{\partial K}\widehat{B(w)}\left(\widehat{w} - w\right)\bm{n}\mathrm{d}\Gamma.
\end{equation}
Taking $\widehat{B(w)} = \ave{B(w)}$ yields immediately \textit{BR-2015}. Some simple algebraic computations show that one has to consider $\widehat{B(w)} = \ave{B(w)}\left(1 - \frac{1}{2}\jump{B(w)}\right)$ to recover \textit{BR-2023}. A more flexible way to present similar ideas consists of rewriting $B(w)\grad w = \dive(B(w)w) - w\grad B(w)$ and employing two standards integration by parts and approximations, as presented in \eqref{eq:BR}. In particular, choosing $\widehat{B(w)w} = \widehat{B(w)}\widehat{w}$, one always recovers \eqref{eq:BR_non_cons}. Finally, in this work, the elementwise contribution vanishes because we consider a finite volume scheme with piecewise constant approximations inside each element. Hence, the discretization \eqref{eq:BR} falls within the BR-type approach as intended in this paper, following one of the key contributions presented in \cite{bassi:1997a}. The scheme proposed in \cite{bassi:1997a} has been developed within the framework of the discontinuous Galerkin method and is commonly referred to in the literature as BR1. In this case, the definition of $D_{h}w$ involves the introduction of a global lifting operator, as it results from \eqref{eq:BR_grad} and \eqref{eq:BR_grad_1}. A variant of the BR1 approach, commonly known as BR2, based on a local lifting operator has been later proposed in \cite{bassi:1997b} in order to reduce the stencil as well as to enhance stability and accuracy of the method. However, these aspects lie far beyond the scope of the present work and do not affect the above discussion.

\printbibliography[category=cited]
\printbibliography[title={Further Reading}, notcategory=cited, resetnumbers=true]

\end{document}